\documentclass[12pt]{amsart}
\usepackage{latexsym}
\usepackage{float}
\usepackage{amsmath,amscd}
\usepackage{amssymb}
\usepackage{indentfirst}
\usepackage{amsthm}
\usepackage{thmtools}
\usepackage[capitalise]{cleveref}
\usepackage{graphicx}
\usepackage{amsfonts}
\usepackage{mathrsfs}
\usepackage{graphicx}
\usepackage{booktabs}
\usepackage{latexsym}
\usepackage{bm}
\usepackage{dsfont}
\usepackage{array}
\usepackage{color}
\usepackage{epstopdf}
\usepackage{subfigure}
\usepackage[top=1in, bottom=1in, left=1.25in, right=1.25in]{geometry}
\usepackage{mathtools}
\usepackage{leftidx}
\usepackage{tensor}
\usepackage{tikz}
\usepackage{tikz-cd}
\usetikzlibrary{matrix,arrows,decorations.pathmorphing}

\DeclareMathOperator{\Spec}{Spec}

\DeclareMathOperator{\Sym}{Sym}
\DeclareMathOperator{\Div}{Div}
\DeclareMathOperator{\Pic}{Pic}
\DeclareMathOperator{\Res}{Res}
\DeclareMathOperator{\can}{can}
\DeclareMathOperator{\Dim}{Dim}

\DeclareMathOperator{\im}{im}

\DeclareMathOperator{\rank}{rank}
\DeclareMathOperator{\supp}{supp}
\DeclareMathOperator{\Fr}{Fr}

\DeclareMathOperator{\Mat}{Mat}
\DeclareMathOperator{\Sht}{Sht}

\DeclareMathOperator{\TrSht}{TrSht}
\DeclareMathOperator{\RedSht}{RedSht}
\DeclareMathOperator{\all}{all}
\DeclareMathOperator{\ToySht}{ToySht}
\DeclareMathOperator{\RToySht}{RToySht}
\DeclareMathOperator{\LToySht}{LToySht}
\DeclareMathOperator{\ooSht}{\tensor*[^{\circ\circ}]{\Sht}{}}
\DeclareMathOperator{\oToySht}{\tensor*[^\circ]{\ToySht}{}}
\DeclareMathOperator{\oRToySht}{\tensor*[^\circ]{\RToySht}{}}
\DeclareMathOperator{\oLToySht}{\tensor*[^\circ]{\LToySht}{}}

\DeclareMathOperator{\ooToySht}{\tensor*[^{\circ\circ}]{\ToySht}{}}
\DeclareMathOperator{\Grass}{Grass}
\DeclareMathOperator{\oGrass}{\tensor*[^\circ]{\Grass}{}}
\DeclareMathOperator{\Schub}{Schub}
\DeclareMathOperator{\Flag}{Flag}
\DeclareMathOperator{\I}{I}
\DeclareMathOperator{\II}{II}
\DeclareMathOperator{\Id}{Id}
\DeclareMathOperator{\coker}{coker}

\DeclareMathOperator{\Hom}{Hom}

\DeclareMathOperator{\Maps}{Maps}

\DeclareMathOperator{\Aut}{Aut}
\DeclareMathOperator{\Four}{Four}
\DeclareMathOperator{\Av}{Av}
\DeclareMathOperator{\SHom}{\mathscr{H}\kern -2pt \textit{om}}
\DeclareMathOperator{\Art}{Art}

\DeclareMathOperator{\pr}{\mathrm{pr}}
\DeclareMathOperator{\bcirc}{\mathbin{\circ}}
\DeclareMathOperator{\simto}{\xrightarrow{\sim}}

\theoremstyle{definition}
\declaretheorem[parent=subsection]{Definition}

\theoremstyle{remark}
\declaretheorem[sibling=Definition]{Remark}
\theoremstyle{plain}
\declaretheorem[sibling=Definition]{Theorem}
\declaretheorem[sibling=Definition]{Lemma}
\declaretheorem[sibling=Definition]{Proposition}
\declaretheorem[sibling=Definition]{Corollary}

\crefname{Definition}{Definition}{Definitions}
\crefname{Remark}{Remark}{Remarks}
\crefname{Example}{Example}{Examples}
\crefname{Lemma}{Lemma}{Lemmas}
\crefname{Theorem}{Theorem}{Theorems}
\crefname{Proposition}{Proposition}{Propositions}
\crefname{Corollary}{Corollary}{Corollaries}
\crefname{Conjecture}{Conjecture}{Conjectures}
\crefname{Claim}{Claim}{Claims}
\crefname{Section}{Section}{Sections}
\crefname{Chapter}{Chapter}{Chapters}

\numberwithin{equation}{section}

\title{A Toy Model of Shtukas}
\author{Zhiyuan Ding}

\begin{document}
\begin{abstract}
  Motivated by the question of constructing certain rational functions (modular units) on the moduli stack of Drinfeld shtukas, we introduce the notion of toy shtukas. We prove basic properties of the moduli scheme of toy shtukas. Analogously to horospherical divisors on the moduli stack of Drinfeld shtukas, there are toy horospherical divisors on the moduli scheme of toy shtukas. We describe the space of principal toy horospherical divisors. There is a canonical morphism from the moduli stack of Drinfeld shtukas to the moduli scheme of toy shtukas. Our main result is a description of the space of principal horospherical divisors obtained from the pullback.
\end{abstract}

\maketitle

\section{Introduction}
\subsection{Notation and conventions}\label{(Section)notation and conventions}
The following notation and conventions will be used throughout the article.

Fix a prime number $p$ and fix a finite field $\mathbb{F}_q$ of characteristic $p$ with $q$ elements.

For any vector space $V$ over $\mathbb{F}_q$, we denote $\mathbf{P}_{V^*}$ to be the set of codimension 1 subspaces of $V$, and we denote $\mathbf{P}_V$ to be the set of dimension 1 subspaces of $V$.

For any scheme $S$ over $\mathbb{F}_q$, we denote $\Fr_S$ to be its Frobenius endomorphism relative to $\mathbb{F}_q$. For two schemes $S_1$ and $S_2$ over $\mathbb{F}_q$, $S_1\times S_2$ denotes the product of $S_1$ and $S_2$ over $\mathbb{F}_q$, and by a morphism $S_1\to S_2$ we mean a morphism over $\mathbb{F}_q$.

\subsection{The notion of toy shtukas}
Let $V$ be a finite dimensional vector space over $\mathbb{F}_q$. Let $\Grass_V^n$ be the Grassmannian for subspaces of $V$ of dimension $n$. We define a closed subscheme $\ToySht_V^n\subset\Grass_V^n$ whose $\overline{\mathbb{F}_q}$-points consist of subspaces $L\subset V\otimes\overline{\mathbb{F}_q}$ such that $L\cap\Fr^*L$ has codimension at most 1 in $L$. (See \cref{definition of a toy shtuka} for the description of $\ToySht_V^n(S)$ for any scheme $S$ over $\mathbb{F}_q$.)

A point of $\ToySht_V^n$ is called a $\emph{toy shtuka}$ for $V$ of dimension $n$.

Analogously to Drinfeld shtukas, we have the notion of left and right toy shtukas, and there are partial Frobeniuses relating them. See \cref{(Section)definition of toy shtukas} and \cref{(Section)partial Frobenius} for more details.

\subsection{Motivation}
The motivation of our work is to construct certain rational functions on the moduli stack of Drinfeld shtukas and use them to generate linear relations between the classes of horospherical divisors to give an explicit version of Manin-Drinfeld theorem for shtukas. These rational functions are somewhat similar to modular units. In~\cite{Ding17}, some rational functions are constructed and the their divisors are calculated. The notion of toy shtukas introduced in this article not only allows us to construct more rational functions, but also gives more insight into the problem. We hope that \cref{a subspace of principal horospherical Q-divisors} gives all principal horospherical $\mathbb{Z}[\frac{1}{p}]$-divisors.

\subsection{Contents}
The definition of toy shtukas is given in \cref{(Section)definition of toy shtukas}.

In \cref{(Section)the scheme of toy shtukas}, we prove basic properties of $\ToySht_V^n$.

Analogously to horospherical divisors on the stack of shtukas, we have toy horospherical divisors on the scheme of toy shtukas. We give more details here. Let $V$ be a finite dimensional vector space over $\mathbb{F}_q$ with $\dim_{\mathbb{F}_q}V\ge3$. Let $1\le n\le N-1$. For $J\in \mathbf{P}_V$, those toy shtukas for $V$ of dimension $n$ which contain $J$ form a codimension 1 subscheme of $\ToySht_V^n$. There is a similar statement for any $H\in\mathbf{P}_{V^*}$. When restricted to the smooth locus of $\ToySht_V^n$, these codimension 1 subschemes are called \emph{toy horospherical divisors}. We give  $\mathbb{Z}[\frac{1}{p}]$-linear relations between the classes of toy horospherical divisors in \cref{space of principal rational toy horospherical divisors}.

Restrictions of Schubert divisors of $\Grass_V^n$ give divisors on the nonsingular locus of $\ToySht_V^n$. In \cref{(Section)Schubert divisors of the schemes of toy shtukas}, we show that these divisors are toy horospherical divisors, and we give an explicit description of them. These divisors are related to the Knudsen-Mumford divisors on the stack of Drinfeld shtukas.

When dealing with shtukas, to get the action of the adelic group, we consider shtukas with structures of all levels. Correspondingly, we consider toy shtukas for Tate spaces. In \cref{(Section)Tate linear algebra}, we give a brief survey of Tate spaces. Eventually the Tate space will be $\mathbb{A}^d$, where $\mathbb{A}$ is the ring of adeles of some function field, and $d$ is the rank of the shtukas.

Let $T$ be a Tate space over $\mathbb{F}_q$. Roughly speaking, a Tate toy shtuka for $T$ is a discrete lattice of $T$ which is almost preserved by Frobenius. The moduli problem of Tate toy shtukas for $T$ is representable by a scheme $\ToySht_T^n$, which is a closed subschemes of the Sato Grassmannian $\Grass_T^n$. Let $\oToySht_T^n$ be the nontrivial locus of $\ToySht_T^n$, i.e., the locus where the discrete lattice is not preserved by Frobenius. When $T$ is nondiscrete and noncompact, $\oToySht_T^n$ is not locally Noetherian. Therefore, one cannot hope to write a Cartier divisor of $\oToySht_T^n$ as a sum of irreducible ones.

In \cref{space of Tate toy horospherical subschemes}, we identify the (partially) ordered abelian group of Tate toy horospherical divisors of $\oToySht_T^n$ as a certain class of locally constant functions on the totally disconnected topological space $(T-\{0\})\coprod (T^*-\{0\})$. Here we are seriously working with Cartier divisors of schemes which are not locally Noetherian. This seems to be a novel feature of our work.

Our first main result is \cref{space of principal rational Tate toy horospherical divisors}, which describes of the (partially) ordered abelian group of those principal $\mathbb{Z}[\frac{1}{p}]$-divisors of $\oToySht_T^n$ supported on the union of Tate toy horospherical subschemes. The appearance of Fourier transform there gives a link to the intertwining operator for Eisenstein series, which plays an important role in the theory of horospherical divisors on the stack of Drinfeld shtukas.

We relate Drinfeld shtukas with toy shtukas in \cref{(Section)relation between shtukas and toy shtukas}. In Propositions \ref{from a right shtuka to a right toy shtuka}, \ref{from a left shtuka to a left toy shtuka} and \ref{from a shtuka to a toy shtuka}, to a shtuka equipped with level structures satisfying certain vanishing conditions on its cohomology, we functorially associate a toy shtuka. In \cref{(Section)from a shtuka with all level structures to a Tate toy shtuka}, we obtain the canonical morphism $\theta$ from the moduli scheme of shtukas with structures of all levels to the moduli scheme of Tate toy shtukas.

In \cref{(Section)pullback of toy horospherical divisors under theta}, we calculate the pullback of a Tate toy horospherical divisor under the morphism $\theta$. The pullback turns out to be an averaging operator, as shown in \cref{formula for pullback of toy horospherical divisors under theta}. Our second main result is \cref{a subspace of principal horospherical Q-divisors}, which gives a subspace of principal horospherical $\mathbb{Z}[\frac{1}{p}]$-divisors of the moduli scheme of shtukas with structures of all levels. We hope that this subspace actually equals the whole space.
\subsection{Acknowledgements}
The research was partially supported by NSF grant DMS-1303100. I would like to express my deep gratitude to my advisor Vladimir Drinfeld for his continual guidance and support.

\section{A toy model of shtukas}\label{(Section)definition of toy shtukas}
Fix a vector space $V$ over $\mathbb{F}_q$ with $\dim_{\mathbb{F}_q}V=N<\infty$. Let $S$ be a scheme over $\mathbb{F}_q$. For an $\mathscr{O}_S$-module $\mathscr{F}$ and a point $s\in S$, we denote $\mathscr{F}_s$ to be the pullback of $\mathscr{F}$ to $s$.

\subsection{Definition of toy shtukas}
\begin{Definition}\label{definition of a toy shtuka}
A \emph{toy shtuka} for $V$ over $S$ is an $\mathscr{O}_S$-submodule $\mathscr{L}\subset V\otimes\mathscr{O}_S$ such that $(V\otimes\mathscr{O}_S)/\mathscr{L}$ is locally free and the composition
\[\Fr_S^*\mathscr{L}\hookrightarrow(\Fr_S^*V)\otimes\mathscr{O}_S=V\otimes\mathscr{O}_S \twoheadrightarrow(V\otimes\mathscr{O}_S)/\mathscr{L}\]
has rank at most 1. (In other words, the corresponding morphism $\bigwedge^2 \Fr_S^*\mathscr{L}\to\bigwedge^2((V\otimes\mathscr{O}_S)/\mathscr{L})$ is zero.)
\end{Definition}

\begin{Definition}
A \emph{right toy shtuka} for $V$ over $S$ is a pair of $\mathscr{O}_S$-submodules $\mathscr{L},\mathscr{L}'\subset V\otimes\mathscr{O}_S$ such that $(V\otimes\mathscr{O}_S)/\mathscr{L}$ and $(V\otimes\mathscr{O}_S)/\mathscr{L}'$ are locally free, $\mathscr{L}\subset\mathscr{L}'$, $\Fr_S^*\mathscr{L}\subset\mathscr{L}'$, and $\dim\mathscr{L}'_s-\dim\mathscr{L}_s=1$ for every $s\in S$.
\end{Definition}

\begin{Remark}
In the exact sequence
\[\begin{tikzcd}
  0\arrow[r]&\mathscr{L}'/\mathscr{L}\arrow[r]&(V\otimes\mathscr{O}_S)/\mathscr{L}\arrow[r]&(V\otimes\mathscr{O}_S)/\mathscr{L}'\arrow[r]&0,
\end{tikzcd}\]
both $(V\otimes\mathscr{O}_S)/\mathscr{L}$ and $(V\otimes\mathscr{O}_S)/\mathscr{L}'$ are locally free. So $\mathscr{L}'/\mathscr{L}$ is also locally free. Hence the condition $\dim\mathscr{L}'_s-\dim\mathscr{L}_s=1$ for every $s\in S$ implies that $\mathscr{L}'/\mathscr{L}$ is invertible. Similarly $\mathscr{L}'/\Fr_S^*\mathscr{L}$ is invertible.
\end{Remark}

\begin{Definition}
A \emph{left toy shtuka} for $V$ over $S$ is a pair of $\mathscr{O}_S$-submodules $\mathscr{L},\mathscr{L}'\subset V\otimes\mathscr{O}_S$ such that $(V\otimes\mathscr{O}_S)/\mathscr{L}$ and $(V\otimes\mathscr{O}_S)/\mathscr{L}'$ are locally free, $\mathscr{L}'\subset\mathscr{L}$, $\mathscr{L}'\subset\Fr_S^*\mathscr{L}$, and $\dim\mathscr{L}_s-\dim\mathscr{L}'_s=1$ for every $s\in S$.
\end{Definition}

\begin{Remark}
Similarly to the above remark, $\mathscr{L}/\mathscr{L}'$ and $(\Fr_S^*\mathscr{L})/\mathscr{L}'$ are invertible.
\end{Remark}

\section{The schemes of toy shtukas}\label{(Section)the scheme of toy shtukas}
In this section, schemes (e.g. Grassmannians, the schemes of matrices) are defined over $\mathbb{F}_q$.

Fix an vector space $V$ over $\mathbb{F}_q$ with $\dim_{\mathbb{F}_q}V=N<\infty$.

For $0\le n \le N$, let $\Grass_V^n$ denote the Grassmannian of $n$-dimensional subspaces of $V$.

\subsection{Definitions of the schemes of toy shtukas}
\begin{Definition}
Let $\ToySht_V^n$ (resp. $\LToySht_V^n$, resp. $\RToySht_V^n$) be the functor which associates to each $\mathbb{F}_q$-scheme $S$ the set of isomorphism classes of toy shtukas $\mathscr{L}\subset V\otimes\mathscr{O}_S$ (resp. left toy shtukas $\mathscr{L}'\subset\mathscr{L}\subset V\otimes\mathscr{O}_S$, resp. right toy shtukas $\mathscr{L}\subset\mathscr{L}'\subset V\otimes\mathscr{O}_S$)  such that $\rank\mathscr{L}=n$.
\end{Definition}

\begin{Remark}
From the definition of toy shtukas, we see that $\ToySht_V^n$ (resp. $\LToySht_V^n$, resp. $\RToySht_V^n$) is representable by a subschemes of $\Grass_V^n$ (resp. $\Grass_V^{n-1}\times\Grass_V^n$, resp. $\Grass_V^n\times\Grass_V^{n+1}$). For explicit local description of these schemes, see \cref{(Section)explicit local description of ToySht,(Section)explicit local description of LToySht,(Section)explicit local description of RToySht}.
\end{Remark}

\begin{Remark}\label{forgetting morphism to the scheme of toy shtukas}
There is a natural morphism $\LToySht_V^n\to\ToySht_V^n$ which maps a left toy shtuka $\mathscr{L}'\subset\mathscr{L}$ to $\mathscr{L}$, and there is a natural morphism $\RToySht_V^n\to\ToySht_V^n$ which maps a right toy shtuka $\mathscr{L}\subset\mathscr{L}'$ to $\mathscr{L}$.
\end{Remark}

\subsection{Explicit local description of the scheme of toy shtukas}\label{(Section)explicit local discription of the scheme of toy shtukas}
\subsubsection{Affine open charts of Grassmannians}\label{(Section)affine open charts of Grassmannians}
For a finite dimensional vector space $M$ over $\mathbb{F}_q$, we denote $\underline{M}=\Spec\Sym M^*$. If $M=\Hom(W',W)$, we write  $\underline{\Hom}(W',W)$ instead of $\underline{\Hom(W',W)}$.

Fix an $(N-n)$-dimensional subspace $W$ of $V$. Denote $U_W$ to be the open subscheme of $\Grass_V^n$ parameterizing those $n$-dimensional subspaces of $V$ which are transversal to $W$. Fix an $n$-dimensional subspace $W'$ of $V$ such that $V=W\oplus W'$. Then we have an identification $U_W=\underline{\Hom}(W',W)$. For an $\mathbb{F}_q$-algebra $R$, an $R$-point of $U_W$ is the graph of an $R$-linear map $W'\otimes R\to W\otimes R$.

Define the Artin-Schreier morphism $AS_{W',W}: \underline{\Hom}(W',W)\to \underline{\Hom}(W',W)$ by $AS_{W',W}=\Id-\Fr$. We know that $AS_{W',W}$ is surjective finite \'etale of degree $q^{n(N-n)}$.

Define a closed subscheme $U_W^{\le 1}:=\underline{\Hom}(W',W)^{\rank\le 1}\subset\underline{\Hom}(W',W)=U_W$ whose $R$-points are $R$-linear maps $A:W'\otimes R\to W\otimes R$ such that $\bigwedge^2 A=0$.

\subsubsection{Explicit local description of $\ToySht$}\label{(Section)explicit local description of ToySht}
We consider $\ToySht_V^n$ as a subscheme of $\Grass_V^n$.
\begin{Lemma}\label{explicit local description of ToySht}
When $1\le n\le N-1$, $\ToySht_V^n\cap U_W$ is the inverse image of $U_W^{\le 1}$ under $AS_{W',W}$. In particular $\ToySht_V^n$ is a closed subscheme of $\Grass_V^n$.
\end{Lemma}

\begin{proof}
The statement follows from the definition of toy shtukas.
\end{proof}

\subsubsection{Explicit local description of $\LToySht$}\label{(Section)explicit local description of LToySht}
For $0\le i\le j\le N$, let $\Flag_V^{i,j}$ be the closed subscheme of $\Grass_V^i\times\Grass_V^j$ which consists of pairs $(M_1,M_2)$ such that $M_1\subset M_2$.

We consider $\LToySht_V^n$ as a subscheme of $\Flag_V^{n-1,n}$.

Let $R$ be an $\mathbb{F}_q$-algebra. Let $A:W'\otimes R\to W\otimes R$ be an $R$-linear morphism and $\Gamma_A$ be its graph. Projection from $\Gamma_A$ to $W'\otimes R$ induces a bijective correspondence between those submodules of $\Gamma_A$ whose quotient is locally free and those submodules of $W'\otimes R$ whose quotient is locally free. This bijection is functorial in $R$. In this way we get an identification $\Flag_V^{n-1,n}\cap(\Grass_V^{n-1}\times U_W)=\Grass_{W'}^{n-1}\times U_W$.

Define a closed subscheme $C_{W,W'}^\flat\subset\Grass_{W'}^{n-1}\times U_W=\Grass_{W'}^{n-1}\times \underline{\Hom}(W',W)$ consisting of pairs $(H,A)$ such that $H\subset\ker A$.

\begin{Lemma}\label{explicit local description of LToySht}
When $1\le n\le N-1$, $\LToySht_V^n\cap(\Grass_V^{n-1}\times U_W)$ is the inverse image of $C_{W,W'}^\flat$ under $\Id_{\Grass_{W'}^{n-1}}\times AS_{W',W}$. In particular, $\LToySht_V^n$ is a closed subscheme of $\Flag_V^{n-1,n}$.
\end{Lemma}
\begin{proof}
The statement follows from the definition of left toy shtukas.
\end{proof}

\begin{Lemma}\label{smoothness and dimension of LToySht}
When $1\le n\le N-1$, $\LToySht_V^n$ is smooth of pure dimension $N-1$ over $\mathbb{F}_q$.
\end{Lemma}
\begin{proof}
We observe that $C_{W,W'}^\flat$ is a $(N-n)$-dimensional vector bundle over $\Grass_{W'}^{n-1}$. Hence it is smooth of pure dimension $N-1$. Since  $\Id_{\Grass_{W'}^{n-1}}\times AS_{W',W}$ is \'etale, the statement follows from \cref{explicit local description of LToySht}.
\end{proof}
\subsubsection{Explicit local description of $\RToySht$}\label{(Section)explicit local description of RToySht}
We consider $\RToySht_V^n$ as a subscheme of $\Flag_V^{n,n+1}$.

Let $R$ be an $\mathbb{F}_q$-algebra. Let $A:W'\otimes R\to W\otimes R$ be an $R$-linear morphism and $\Gamma_A$ be its graph. Intersection with $W\otimes R$ induces a bijective correspondence between those submodules $M\subset V\otimes R$ containing $\Gamma_A$ such that $(V\otimes R)/M$ is locally free and those submodules of $W\otimes R$ whose quotient is locally free. This bijection is functorial in $R$. In this way we get an identification $\Flag_V^{n,n+1}\cap(U_W\times\Grass_V^{n+1})=U_W\times\Grass_W^1$.

Define a closed subscheme $C_{W,W'}^\sharp\subset U_W\times \Grass_{W}^1=\underline{\Hom}(W',W)\times\Grass_W^1$ consisting of pairs $(A,H)$ such that $\im A\subset H$.

\begin{Lemma}\label{explicit local description of RToySht}
When $1\le n\le N-1$, $\RToySht_V^n\cap(U_W\times\Grass_V^{n+1})$ is the inverse image of $C_{W,W'}^\sharp$ under $AS_{W',W}\times\Id_{\Grass_W^1}$. In particular, $\RToySht_V^n$ is a closed subscheme of $\Flag_V^{n,n+1}$.
\end{Lemma}
\begin{proof}
The statement follows from the definition of right toy shtukas.
\end{proof}

\begin{Lemma}\label{smoothness and dimension of RToySht}
When $1\le n\le N-1$, $\RToySht_V^n$ is smooth of pure dimension $N-1$ over $\mathbb{F}_q$.
\end{Lemma}
\begin{proof}
We observe that $C_{W,W'}^\sharp$ is an $n$-dimensional vector bundle over $\Grass_W^1$. Hence it is smooth of pure dimension $N-1$. Since  $AS_{W',W}\times\Id_{\Grass_W^1}$ is \'etale, the statement follows from \cref{explicit local description of RToySht}.
\end{proof}

\subsection{Basic properties of $\ToySht$}
Determinantal varieties are proved to be Cohen-Macaulay in~\cite{HE}. See Section 3 of Chapter 2 of~\cite{ACGH} for a review of basic properties of determinantal varieties.

We know that for $1\le n\le N-1$, $\Mat_{n\times(N-n)}^{\rank\le 1}$ is the affine cone over $(\mathbb{P}^\vee)^{N-n-1}\times \mathbb{P}^{n-1}$. In particular, it has pure dimension $N-1$ and it is reduced.

\begin{Lemma}\label{basic properties of the scheme of toy shtukas}
For $1\le n\le N-1$, $\ToySht_V^n$ is reduced and has pure dimension $N-1$. We have $\ToySht_V^1=\Grass_V^1$, $\ToySht_V^{N-1}=\Grass_V^{N-1}$. For $n\in\{1,N-1\}$, $\ToySht_V^n$ is smooth. For $2\le n\le N-2$, $\ToySht_V^n$ is smooth outside $\Grass_V^n(\mathbb{F}_q)$ and the singularity at each point of $\Grass_V^n(\mathbb{F}_q)$ is \'etale locally the vertex of the cone over $(\mathbb{P}^\vee)^{N-n-1}\times \mathbb{P}^{n-1}$.
\end{Lemma}
\begin{proof}
A choice of bases of $W$ and $W'$ identifies $\underline{\Hom}(W',W)^{\rank\le 1}$ with $\Mat_{n\times(N-n)}^{\rank\le 1}$. Since the Artin-Schreier morphism is finite \'etale, the statement follows from the corresponding properties of $\Mat_{n\times(N-n)}^{\rank\le 1}$.
\end{proof}

\subsection{The scheme of nontrivial toy shtukas}
\begin{Definition}
For a toy shtuka (resp. a right toy shtuka, resp. a left toy shtuka) $\mathscr{L}$ over $S$, we say that it is \emph{nontrivial at $s\in S$} if $\Fr_s^*\mathscr{L}_s\ne\mathscr{L}_s$, where $\mathscr{L}_s$ is the pullback of $\mathscr{L}$ to $s$.
\end{Definition}

\begin{Remark}
A left/right toy shtuka is nontrivial at $s$ if and only if it is as a toy shtuka.
\end{Remark}

\begin{Remark}
A pointwise nontrivial left/right toy shtuka is the same as a pointwise nontrivial toy shtuka. See \cref{identification of nontrivial left/right toy shtukas} for a more precise statement.
\end{Remark}

\begin{Remark}\label{description of the trivial locus}
We know that the trivial locus of $\ToySht_V^n$ is equal to the Frobenius fixed points, or equivalently it is the intersection of the graph of Frobenius morphism and the diagonal. So the trivial locus of $\ToySht_V^n$ is $\Grass_V^n(\mathbb{F}_q)$, a reduced 0-dimensional closed subscheme of $\Grass_V^n$.
\end{Remark}

\begin{Definition}
Let $\oToySht_V^n$ (resp. $\oLToySht_V^n$, resp. $\oRToySht_V^n$) be the nontrivial locus of $\ToySht_V^n$ (resp. $\RToySht_V^n$, resp. $\LToySht_V^n$).
\end{Definition}

\begin{Remark}\label{description of the nontrivial locus}
Since nontrivialness is an open condition for toy shtukas (resp. left toy shtukas, resp. right toy shtukas), $\oToySht_V^n$ (resp. $\oLToySht_V^n$, resp. $\oRToySht_V^n$) is an open subscheme of $\ToySht_V^n$ (resp. $\RToySht_V^n$, resp. $\LToySht_V^n$). We have $\oToySht_V^n=\ToySht_V^n-\Grass_V^n(\mathbb{F}_q)$ from \cref{description of the trivial locus}. In particular, we know that $\oToySht_V^n$ is smooth over $\mathbb{F}_q$. We also see that $\oToySht_V^n$ is nonempty if and only if $1\le n\le N-1$.
\end{Remark}

\begin{Lemma}\label{intersection and sum are subbundles for nontrivial toy shtukas}
Let $\mathscr{L}$ be a toy shtuka over $S$ which is nontrivial at each $s\in S$. Then $\mathscr{L}+\Fr_S^*\mathscr{L}$ and $\mathscr{L}\cap\Fr_S^*\mathscr{L}$ are subbundles of $\mathsf{p}_S^*V$. Also, $(\mathscr{L}+\Fr_S^*\mathscr{L})/\mathscr{L}$ and $\mathscr{L}/(\mathscr{L}\cap\Fr_S^*\mathscr{L})$ are invertible sheaves.
\end{Lemma}
\begin{proof}
By the definition of toy shtukas, the morphism $\Fr_S^*\mathscr{L}\to\mathsf{p}_S^*V/\mathscr{L}$ has rank at most 1. It has strictly constant rank 1 since $\mathscr{L}$ is nontrivial at each $s\in S$. The statements now follow from \cref{equivalent conditions of the relative position of two subbundles}.
\end{proof}

\begin{Corollary}\label{identification of nontrivial left/right toy shtukas}
For $1\le n\le N$, the natural morphisms $\LToySht_V^n\to\ToySht_V^n$ and $\RToySht_V^n\to\ToySht_V^n$ in \cref{forgetting morphism to the scheme of toy shtukas} induce isomorphisms $\oLToySht_V^n\to\oToySht_V^n$ and $\oRToySht_V^n\to\oToySht_V^n$.
\end{Corollary}

\subsection{Partial Frobeniuses}\label{(Section)partial Frobenius}
We have the following constructions for left/right toy shtukas:

(i) For a left toy shtuka $\mathscr{L}'\subset\mathscr{L}$ over an $\mathbb{F}_q$-scheme $S$, the pair $\mathscr{L}'\subset\Fr_S^*\mathscr{L}$ forms a right toy shtuka over $S$.

(ii) For a right toy shtuka $\mathscr{L}\subset\mathscr{L}'$ over an $\mathbb{F}_q$-scheme $S$, the pair $\Fr_S^*\mathscr{L}\subset\mathscr{L}'$ forms a left toy shtuka over $S$.

\begin{Definition}\label{definition of partial Frobeniuses}
We define partial Frobeniuses $F_{V,n}^-:\LToySht_V^{n}\to\RToySht_V^{n-1}, (1\le n\le N)$ and $F_{V,n}^+:\RToySht_V^n\to\LToySht_V^{n+1}, (0\le n\le N-1)$ induced by the above constructions.
\end{Definition}

\begin{Lemma}\label{composition is universal homeomorphism implies second morphism is}
Let $f_1:Y_1\to Y_2,f_2:Y_2\to Y_3$ be two morphisms of schemes. If $f_1$ is surjective and $f_2\bcirc f_1$ is a universal homeomorphism, then $f_2$ is also a universal homeomorphism.
\end{Lemma}
\begin{proof}
This is Proposition 3.8.2(iv) of~\cite{EGA_I_2nd}.
\end{proof}

\begin{Lemma}\label{mutual inverses up to universal homeomorphism}
For two morphisms $f_1:Y_1\to Y_2,f_2:Y_2\to Y_1$ between two schemes $Y_1$ and $Y_2$, if $f_2\bcirc f_1$ and $f_1\bcirc f_2$ are universal homeomorphisms, then $f_1$ and $f_2$ are universal homeomorphisms.
\end{Lemma}
\begin{proof}
Since $f_1\bcirc f_2$ is a homeomorphism, $f_1$ is surjective. Thus $f_2$ is a universal homeomorphism by \cref{composition is universal homeomorphism implies second morphism is}. Similarly, $f_1$ is also a universal homeomorphism.
\end{proof}

\begin{Lemma}\label{universal homeomorphisms between schemes of left and right toy shtukas}
Assume $0\le n\le N-1$. We have $F_{V,n+1}^-\bcirc F_{V,n}^+=\Fr_{\RToySht_V^n}$, $F_{V,n}^+\bcirc F_{V,n+1}^-=\Fr_{\LToySht_V^{n+1}}$. In particular, $F_{V,n+1}^-$ and $F_{V,n}^+$ induce universal homeomorphisms between $\RToySht_V^{n}$ and $\LToySht_V^{n+1}$.
\end{Lemma}
\begin{proof}
The first statement follows from definition of partial Frobeniuses. The second statement follows from \cref{mutual inverses up to universal homeomorphism}.
\end{proof}

\subsection{Irreducibility of the scheme of toy shtukas}
Recall that we have $\oToySht_V^n=\oLToySht_V^n=\oRToySht_V^n$ for $1\le n\le N$ by \cref{identification of nontrivial left/right toy shtukas}.

\begin{Lemma}\label{nontrivial shtukas are dense}
For $1\le n\le N-1$, $\oToySht_V^n$ is dense in $\ToySht_V^n$, $\LToySht_V^n$ and $\RToySht_V^n$.
\end{Lemma}
\begin{proof}
Since $\ToySht_V^n$ has pure dimension $N-1$ by \cref{basic properties of the scheme of toy shtukas}, we see from \cref{description of the nontrivial locus} that $\oToySht_V^n=\ToySht_V^n-\Grass_V^n(\mathbb{F}_q)$ is dense in $\ToySht_V^n$.

Consider the morphism $\LToySht_V^n\to \ToySht_V^n$. It induces an isomorphism $\oLToySht_V^n\to\oToySht_V^n$. For any $M\in \Grass_V^n(\mathbb{F}_q)$, wee see that the inverse image of $M$ in $\LToySht_V^n$ is isomorphic to $\Grass_M^{n-1}$, and we have $\dim\Grass_M^{n-1}=n-1<N-1$. Since $\LToySht_V^n$ has pure dimension $N-1$ by \cref{smoothness and dimension of LToySht}, $\oLToySht_V^n$ is dense in $\LToySht_V^n$.

The proof for $\RToySht_V^n$ is similar, with $\Grass_M^{n-1}$ replaced by $\Grass_{V/M}^1$.
\end{proof}

\begin{Proposition}\label{irreducibility of the scheme of toy shtukas}
For $1\le n\le N-1$, $\ToySht_V^n$, $\LToySht_V^n$ and $\RToySht_V^n$ are geometrically irreducible.
\end{Proposition}
\begin{proof}
The three schemes in question all contain a dense subscheme $\oToySht_V^n$ by \cref{nontrivial shtukas are dense}. This is still true after base change to $\overline{\mathbb{F}}_q$ since the morphism $\Spec\overline{\mathbb{F}}_q\to \Spec \mathbb{F}_q$ is universally open. Now \cref{universal homeomorphisms between schemes of left and right toy shtukas} and \cref{identification of nontrivial left/right toy shtukas} reduce the question to the case $n=1$. We know that $\ToySht_V^1=\Grass_V^1$ is geometrically irreducible.
\end{proof}

\section{A lemma for toy shtukas inspired by a fact about reducible shtukas}
\subsection{A fact about reducible shtukas}
We have the following fact about reducible shtukas over a field from the paragraph before Proposition 4.2 of~\cite{Drinfeld1}.

Let $X$ be a smooth projective geometrically connected curve over $\mathbb{F}_q$. Let $E$ be a field over $\mathbb{F}_q$. Denote $\Phi_E=\Id_X\otimes\Fr_E:X\otimes E\to X\otimes E$. Let $\alpha,\beta:\Spec E\to X$ be two morphisms such that $\Gamma_\alpha\cap\Gamma_\beta=\emptyset$. Let $\Phi_E^*\mathscr{F}\hookrightarrow\mathscr{F}'\hookleftarrow\mathscr{F}$ be a right shtuka of rank $d$ over $\Spec E$ with zero $\alpha$ and pole $\beta$. Suppose there exists a nonzero subsheaf $\mathscr{E}\subset\mathscr{F}$ of rank $r<d$ such that  $\Phi_E^*\mathscr{E}\subset\mathscr{E}(\Gamma_\beta)$. Let $\mathscr{A}$ be the saturation of $\mathscr{E}$ in $\mathscr{F}$ and let $\mathscr{B}=\mathscr{F}/\mathscr{A}$. Then one of the two possibilities hold:

(i) $\mathscr{A}$ is a right shtuka with zero $\alpha$ and pole $\beta$, and the image of the morphism $\Phi_E^*\mathscr{B}\to\mathscr{B}(\Gamma_\beta)$ is equal to $\mathscr{B}$.

(ii) The image of the morphism $\Phi_E^*\mathscr{A}\to\mathscr{A}(\Gamma_\beta)$ is equal to $\mathscr{A}$ and $\mathscr{B}$ is a right shtuka with zero $\alpha$ and pole $\beta$.

\subsection{A lemma for toy shtukas}
Let $E$ be a field over $\mathbb{F}_q$ and let $L\in\ToySht_V(E)$; in other words, $L$ is a subspace of $V\otimes E$ is a subspace such that $\dim L-\dim((\Fr_E^*L)\cap L)=\dim \Fr_E^*L-\dim((\Fr_E^*L)\cap L)\le 1$.

Let $W$ be a subspace of $V$.

Denote $L'=L\cap (W\otimes E),L''=\im(L\to (V/W)\otimes E)$. We have $\Fr_E^*L'=(\Fr_E^*L)\cap(W\otimes E),\Fr_E^*L''=\im(\Fr_E^*L\to(V/W)\otimes E)$.

\begin{Lemma}\label{triviality of intersection or quotient of intersection}
At least one of the following two statements holds.

(i) $\Fr_E^*L'=L'$;

(ii) $\Fr_E^*L''=L''$.
\end{Lemma}
\begin{proof}
First note that $\dim L'=\dim\Fr_E^*L'$, $\dim L''=\dim\Fr_E^*L''$. Suppose there exists $v\in \Fr_E^*L'$ such that $v\notin L'$. Then $\Fr_E^*L$ is contained in the linear span of $L$ and $v$. Since $v$ maps to $0$ in $(V/W)\otimes E$, we have $\Fr_E^*L''\subset L''$. Hence $\Fr_E^*L''=L''$ by dimension comparison.
\end{proof}

\section{Schubert divisors of the schemes of toy shtukas}\label{(Section)Schubert divisors of the schemes of toy shtukas}
Fix a finite dimensional vector space $V$ over $\mathbb{F}_q$ with $\dim_{\mathbb{F}_q}V=N\ge 3$.

Fix an integer $n$ such that $1\le n\le N-1$, and a subspace $W\subset V$ of codimension $n$.

As in \cref{(Section)affine open charts of Grassmannians}, for a codimension $n$ subspace $M$ in $V$, let $U_{M}$ be the affine chart of $\Grass_V^n$ parameterizing those $n$-dimensional subspaces of $V$ which are transversal to $M$.
\subsection{Definition of Schubert divisors}
Define the closed subscheme $\Schub_V^W\subset\Grass^n_V$ by the following equation:
\[\Schub_V^W:=\{L\in\Grass^n_V|\det(L\to V/W)=0\}.\]
It is easy to prove and well-known that $\Schub_V^W\subset\Grass^n_V$ is an irreducible divisor (a Schubert variety).

\begin{Lemma}
$\Schub_V^W\cap\ToySht_V^n$ is a Cartier divisor in $\ToySht^n_V$.
\end{Lemma}
\begin{proof}
We know that $\ToySht_V^n$ is irreducible and reduced, and $\Schub_V^W$ is a Cartier divisor of $\Grass_V^n$. To prove the statement, it suffices to show that $\ToySht_V^n$ is not contained in $\Schub_V^W$. We have $U_W=\Grass_V^n-\Schub_V^W$, and we see in \cref{(Section)explicit local description of ToySht} that $\ToySht_V^n\cap U_W$ is nonempty.
\end{proof}

We have a perfect complex
\[\mathscr{S}_{V,W}^\bullet=(\mathscr{S}_{V,W}^{-1}\to\mathscr{S}_{V,W}^0)\]
on $\Grass_V^n$, where $\mathscr{S}_{V,W}^{-1}\subset V\otimes\mathscr{O}_{\Grass_V^n}$ is the universal locally free sheaf on $\Grass_V^n$, $\mathscr{S}_{V,W}^0=(V/W)\otimes\mathscr{O}_{\Grass_V^n}$, and the morphism $\mathscr{S}_{V,W}^{-1}\to\mathscr{S}_{V,W}^0$ is the natural one.

Since $W$ has codimension $n$ in $V$, the morphism $\mathscr{S}_{V,W}^{-1}\to\mathscr{S}_{V,W}^0$ is an isomorphism on the dense open subscheme $U_W\subset\Grass_V^n$. So the complex $\mathscr{S}_{V,W}^\bullet$ is good in the sense of Knudsen-Mumford.

\begin{Remark}\label{equivalence of two definitions of the Schubert divisor}
  We have $\Schub_V^W=\Div(\mathscr{S}_W^\bullet)$. (See Chapter II of~\cite{KM} for the definition of $\Div$.)
\end{Remark}

We call $\Schub_V^W\cap\ToySht_V^n$ (resp. $\Schub_V^W\cap\oToySht_V^n$) the \emph{Schubert divisor} of $\ToySht_V^n$ (resp. $\oToySht_V^n$) for $W$.

\subsection{Description of Schubert divisors}\label{(Section)description of Schubert divisors}
For any subspace $H\subset V$, we have a natural closed immersion $\ToySht_H^n\to \ToySht_V^n$. We consider $\ToySht_H^n$ as a subscheme of $\ToySht_V^n$. It is a Weil divisor if $H$ has codimension 1 in $V$. If $H$ contains $W$, then $\ToySht_H^n\subset \Schub_V^W\cap\ToySht_V^n$.

For a subspace $J\subset V$, there is a natural bijection between subspaces of $V/J$ and subspaces of $V$ containing $J$. So we get a natural closed immersion $\ToySht_{V/J}^{n-\dim J}\to\ToySht_V^n$. We consider $\ToySht_{V/J}^{n-\dim J}$ as a subscheme of $\ToySht_V^n$. It is a Weil divisor if $J$ has dimension 1. If $J$ is contained in $W$, then $\ToySht_{V/J}^{n-\dim J}\subset\Schub_V^W\cap\ToySht_V^n$.

The following statement describes the Schubert divisors. Note that by \cref{irreducibility of the scheme of toy shtukas},  $\ToySht_H^n$ is irreducible when $n<\dim H$, and $\ToySht_{V/J}^{n-\dim J}$ is irreducible when $n>\dim J$. Also note that by \cref{description of the nontrivial locus}, $\oToySht_H^n$ is empty when $n=\dim H$, and $\oToySht_{V/J}^{n-\dim J}$ is empty when $n=\dim J$.

Recall that we denote $\mathbf{P}_{V^*}$ to be the set of codimension 1 subspaces of $V$ and we denote $\mathbf{P}_V$ to be the set of dimension 1 subspaces of $V$.

\begin{Theorem}\label{description of Schubert divisors of the scheme of toy shtukas}
We have an equality of Cartier divisors of $\oToySht_V^n$
\[\Schub_V^W\cap\oToySht_V^n=\sum_{\substack{H\in\mathbf{P}_{V^*}\\H\supset W}}\oToySht_H^n+\sum_{\substack{J\in\mathbf{P}_{V}\\J\subset W}}\oToySht_{V/J}^{n-1}.\]
\end{Theorem}
\begin{proof}
Put $Z_V^W=\Schub_V^W\cap\oToySht_V^n$.

From the above discussion we know that $\oToySht_H^n$ and $\oToySht_{V/J}^{n-1}$ are Weil divisors of $\oToySht_V^n$, hence Cartier divisors of $\oToySht_V^n$ since $\oToySht_V^n$ is smooth. We also see that $\oToySht_H^n\subset Z_V^W$ if $H\supset W$ and $\oToySht_{V/J}^{n-1}\subset Z_V^W$ if $J\subset W$.

On the other hand, if $s\in Z_V^W$, then \cref{triviality of intersection or quotient of intersection} implies that $s\in \ToySht_H^n$ for some $H\in\mathbf{P}_{V^*}, H\supset W$ or $s\in \ToySht_{V/J}^{n-1}$ for some $J\in\mathbf{P}_V, J\subset W$ (or both). This shows that the statement is set-theoretically true.

It remains to show that the multiplicity of $Z_V^W$ at each $\oToySht_H^n$ or $\oToySht_{V/J}^{n-1}$ is one. This follows from \cref{multiplicity freeness of Schubert divisor}.
\end{proof}

\subsection{An open subscheme of $\Schub_V^W\cap\ToySht_V^n$}
Define a closed subset $(\Schub_V^W)^{\ge 2}\subset \Schub_V^W$ by the condition
\[(\Schub_V^W)^{\ge 2}:=\{L\in\Grass^n_V|\rank(L\to V/W)\le n-2\}.\]

Let $(\Schub_V^W)^1=\Schub_{V}^{W}-(\Schub_V^W)^{\ge 2}$ be the open subscheme of $\Schub_V^W$.

\begin{Remark}\label{open cover of Schub_W^1}
We know that $(\Schub_V^W)^1=\Schub_V^W\cap(\cup U_{M})$ where $M$ runs through all codimension $n$ subspaces of $V$ whose intersection with $W$ has codimension $n+1$ in $V$.
\end{Remark}

\begin{Remark}\label{equation for a Schubert cell}
Fix one such $M$ as in \cref{open cover of Schub_W^1} and choose an $n$-dimensional subspace $M'$ of $V$ such that $V=M\oplus M'$ and $M'\cap W\ne 0$. Then $\dim (M'\cap W)=1$. We identify $U_{M}$ with $\underline{\Hom}(M',M)$ as in \cref{(Section)affine open charts of Grassmannians}. We see that $\Schub_V^W\cap U_{M}=(\Schub_V^W)^1\cap U_{M}$ as a subscheme of $U_{M}$ is defined by the condition that the induced morphism $M'\cap W\to M/(M\cap W)$ is zero.
\end{Remark}

\begin{Lemma}\label{codimension of higher Schubert cycle}
$(\Schub_V^W)^{\ge 2}\cap\ToySht_V^n$ has codimension at least 2 in $\ToySht_V^n$.
\end{Lemma}
\begin{proof}
From \cref{triviality of intersection or quotient of intersection} we see that $s\in (\Schub_V^W)^{\ge 2}\cap\ToySht_V^n$ implies $s\in \ToySht_{V/J}^{n-2}$ for some 2-dimensional subspace $J\subset V$ or $s\in \ToySht_H^n$ for some codimension 2 subspace $H\subset V$ (or both). The statement follows from the result about dimensions of $\ToySht$ in \cref{basic properties of the scheme of toy shtukas}.
\end{proof}
\subsection{Transversality}
For $1\le a\le s,1\le b\le t$, let $\Mat_{s\times t}^{(a,b)=0}$ denote the subscheme of $\Mat_{s\times t}$ consisting of matrices whose $(a,b)$-entry is zero.

\begin{Lemma}\label{coordinate hyperplane meets determinantal variety transversally}
The locus of $\Mat_{s\times t}^{\rank\le 1}$ where it does not meet $\Mat_{s\times t}^{(a,b)=0}$ transversally is
\[\Mat_{s\times t}^{\rank\le 1}\cap(\bigcap_{i=1}^s\Mat_{s\times t}^{(i,b)=0})\cap(\bigcap_{j=1}^t\Mat_{s\times t}^{(a,j)=0}),\]
i.e., the locus where the $a$-th row and the $b$-th column is zero. In particular, this locus is isomorphic to $\Mat_{(s-1)\times(t-1)}^{\rank\le 1}$ and it has codimension 2 in $\Mat_{s\times t}^{\rank\le 1}$. \qed
\end{Lemma}

\begin{Lemma}\label{coordinate hyperplane meets ToySht transversally}
Let $AS=\Id-\Fr:\Mat_{s\times t}\to\Mat_{s\times t}$ be the Artin-Schreier morphism.  Then $\Mat_{s\times t}^{(a,b)=0}\cap AS^{-1}(\Mat_{s\times t}^{\rank\le 1}-\{0\})$ as a divisor of $AS^{-1}(\Mat_{s\times t}^{\rank\le 1}-\{0\})$ has multiplicity one at each irreducible component.
\end{Lemma}
\begin{proof}
If we identify the tangent spaces at each point of $\Mat_{s\times t}$ with the vector space $\Mat_{s\times t}$ itself, then $AS$ induces the identity on the tangent spaces. Now the statement follows from \cref{coordinate hyperplane meets determinantal variety transversally} and the fact that $AS$ is finite.
\end{proof}

\begin{Lemma}\label{multiplicity freeness of Schubert divisor}
$\Schub_V^W\cap\oToySht_V^n$ is a reduced scheme.
\end{Lemma}
\begin{proof}
By \cref{codimension of higher Schubert cycle}, it suffices to show that $(\Schub_V^W)^1\cap\oToySht_V^n$ as a divisor of $\oToySht_V^n-(\Schub_V^W)^{\ge 2}$ has multiplicity one at each irreducible component. By \cref{open cover of Schub_W^1}, it suffices to show that $(\Schub_V^W)^1\cap\oToySht_V^n\cap U_{M}$ as a divisor of $\oToySht_V^n\cap U_{M}$ has multiplicity one at each irreducible component, where $M$ is any codimension $n$ subspace of $V$ whose intersection with $W$ has codimension $n+1$ in $V$. From the description of $(\Schub_V^W)^1\cap U_{M}$ in \cref{equation for a Schubert cell}, the statement follows from \cref{explicit local description of ToySht} and \cref{coordinate hyperplane meets ToySht transversally}.
\end{proof}

\section{Toy Horospherical Divisors}
\subsection{Notation}
Fix a finite dimensional vector space $V$ over $\mathbb{F}_q$ with $\dim_{\mathbb{F}_q}V=N\ge 3$. Fix $n\in\mathbb{Z}$ such that $1\le n\le N-1$.

As in \cref{(Section)description of Schubert divisors}, for any $H\in\mathbf{P}_{V^*}$, we consider $\oToySht_H^n$ (resp. $\LToySht_H^n$, resp.  $\RToySht_H^n$) as a Cartier divisor of $\oToySht_V^n$ (resp. $\LToySht_V^n$, resp.  $\RToySht_V^n$). For any $J\in\mathbf{P}_{V}$, we consider $\oToySht_{V/J}^{n-1}$ (resp. $\LToySht_{V/J}^{n-1}$, resp. $\RToySht_{V/J}^{n-1}$) as a Cartier divisor of $\oToySht_V^n$ (resp. $\LToySht_V^n$,  resp. $\RToySht_V^n$). These divisors are called \emph{toy horospherical divisors} of $\oToySht_V^n$ (resp. $\LToySht_V^n$, resp. $\RToySht_V^n$).

As in \cref{definition of partial Frobeniuses}, we have partial Frobeniuses $F_{V,n}^-:\LToySht_V^{n}\to\RToySht_V^{n-1}, (1\le n\le N)$ and $F_{V,n}^+:\RToySht_V^n\to\LToySht_V^{n+1}, (0\le n\le N-1)$.

\subsection{Short complexes related with principal toy horospherical divisors}\label{(Section)short complexes related with principal toy horospherical divisors}
Let $Y$ be a smooth scheme over $\mathbb{F}_q$ and $U\subset Y$ be a dense open subscheme. Then we have a short complex
\[\begin{tikzcd}
  0\arrow[r]&C^0(Y,U)\arrow[r,"d"]&C^1(Y,U)\arrow[r]&0,
\end{tikzcd}\]
where $C^0(Y,U):=H^0(U,\mathscr{O}_U^\times)$, $C^1(Y,U):=\{ \text{divisors on $Y$ with zero restriction to $U$}\}$

The complex $C^\bullet(Y,U)$ is a contravariant functor in $(Y,U)$.

\begin{Remark}\label{removing codimension 2 subsets does not change divisors with specified support}
Let $Z,Z'\subset Y$ be closed subsets of codimension at least 2 such that $U-Z'\subset Y-Z'$. Then the map $C^\bullet(Y,U)\to C^\bullet(Y-Z,U-Z')$ is an isomorphism.
\end{Remark}

\begin{Remark}\label{Frobenius induces multiplication by q on the short complex}
The endomorphism of $C^\bullet(Y,U)$ corresponding to $\Fr:(Y,U)\to(Y,U)$ is multiplication by $q$.
\end{Remark}

We denote
\[\ooToySht_V^n:=\ToySht_V^n-(\bigcup_{H\in\mathbf{P}_{V^*}}\ToySht_H^n)\cup(\bigcup_{J\in\mathbf{P}_{V}}\ToySht_{V/J}^{n-1})\]
to be the open subscheme of $\ToySht_V^n$. We have $\ooToySht_V^n\subset\oToySht_V^n$.

When $2\le n\le N-2$, the complement of $\oToySht_V^n$ has codimension at least 2 in $\LToySht_V^n$ or $\RToySht_V^n$. \cref{removing codimension 2 subsets does not change divisors with specified support} implies that the natural maps of complexes
\[C^\bullet(\LToySht_V^n,\ooToySht_V^n)\to C^\bullet(\oToySht_V^n,\ooToySht_V^n),\]
\[C^\bullet(\RToySht_V^n,\ooToySht_V^n)\to C^\bullet(\oToySht_V^n,\ooToySht_V^n)\]
are isomorphisms. We denote
\[C_{V,n}^\bullet=C^\bullet(Y,\ooToySht_V^n)\]
for $Y=\oToySht_V^n, \LToySht_V^n, \RToySht_V^n$.

We define
\[C_{V,1}^\bullet=C^\bullet(\RToySht_V^1,\ooToySht_V^n),\]
\[C_{V,N-1}^\bullet=C^\bullet(\LToySht_V^{N-1},\ooToySht_V^n).\]

We see that for $1\le n\le N-2$, $C^1_{V,n}=C^1(\RToySht_V^n,\ooToySht_V^n)$ is freely generated by $\RToySht_H^n(H\in\mathbf{P}_{V^*})$ and $\RToySht_{V/J}^{n-1}(J\in\mathbf{P}_{V})$. For $2\le n\le N-1$, $C^1_{V,n}=C^1(\LToySht_V^n,\ooToySht_V^n)$ is freely generated by $\LToySht_H^n(H\in \mathbf{P}_{V^*})$ and $\LToySht_{V/J}^{n-1}(J\in\mathbf{P}_{V})$.

\subsection{Partial Frobeniuses and toy horospherical divisors}
Partial Frobeniuses induce morphisms of complexes
\[\begin{tikzcd}
C^\bullet_{V,1}\arrow[r,shift left=0.5ex,"(F_{V,2}^-)^*"]&C^\bullet_{V,2}\arrow[r,shift left=0.5ex,"(F_{V,3}^-)^*"]\arrow[l,shift left=0.5ex,"(F_{V,1}^+)^*"]&\dots\arrow[r,shift left=0.5ex]\arrow[l,shift left=0.5ex,"(F_{V,2}^+)^*"]&C^\bullet_{V,N-1}\arrow[l,shift left=0.5ex]
\end{tikzcd}\]

\begin{Lemma}\label{multiplicity 1 of pullback of toy horospherical divisors by partial Frobeniuses}
For $1\le n\le N-2$ and $H\in\mathbf{P}_{V^*}$, $(F_{V,n}^+)^*\LToySht_H^{n+1}=\RToySht_H^n$. For $2\le n\le N-1$ and $J\in\mathbf{P}_{V}$, $(F_{V,n}^-)^*\RToySht_{V/J}^{n-2}=\LToySht_{V/J}^{n-1}$.
\end{Lemma}
\begin{proof}
Since $\LToySht_H^{n+1}$ and $\RToySht_H^n$ are reduced and irreducible by \cref{smoothness and dimension of LToySht} and \cref{smoothness and dimension of RToySht}, to prove the first statement, it suffices to show that $(F_{V,n}^+)^{-1}\LToySht_H^{n+1}\subset\RToySht_H^n$.

We have a Cartesian diagram
\[\begin{tikzcd}
(F_{V,n}^+)^{-1}\LToySht_H^{n+1}\arrow[r]\arrow[d]&\RToySht_V^n\arrow[d,"F_{V,n}^+"]\\
\LToySht_H^{n+1}\arrow[r]&\LToySht_V^{n+1}
\end{tikzcd}\]
Suppose we have a morphism $S\to (F_{V,n}^+)^{-1}\LToySht_H^{n+1}$. From the above Cartesian diagram, we get a right toy shtuka $\mathscr{L}\subset\mathscr{L}'$ of rank $n$ over $S$, such that $\mathscr{L}'\subset H\otimes\mathscr{O}_S$. This shows that the morphism $S\to (F_{V,n}^+)^{-1}\LToySht_H^{n+1}$ factors through $\RToySht_H^n$. Hence $(F_{V,n}^+)^{-1}\LToySht_H^{n+1}\subset\RToySht_H^n$.

The proof of the second statement is similar.
\end{proof}

\begin{Lemma}\label{multiplicity q of pullback of toy horospherical divisors by partial Frobeniuses}
For $1\le n\le N-2$ and $J\in\mathbf{P}_{V}$, $(F_{V,n}^+)^*\LToySht_{V/J}^{n}=q\cdot\RToySht_{V/J}^{n-1}$. For $2\le n\le N-1$ and $H\in\mathbf{P}_{V^*}$, $(F_{V,n}^-)^*\RToySht_{H}^{n-1}=q\cdot\LToySht_{H}^{n}$.
\end{Lemma}
\begin{proof}
Since the $F_{V,n}^+:\RToySht_V^n\to\LToySht_V^{n+1}$ and $F_{V,n+1}^-:\LToySht_V^{n+1}\to\RToySht_V^n$ are universal homeomorphisms and their composition is the Frobenius morphism, the statements follow from \cref{multiplicity 1 of pullback of toy horospherical divisors by partial Frobeniuses}.
\end{proof}

The following statement about finite Radon transform is standard.
\begin{Lemma}\label{equivalent conditions for coefficients of toy horospherical divisors}
For $1\le n\le N-1$, The following two conditions for the two sets of numbers $\{\lambda_H\}_{H\in\mathbf{P}_{V^*}}, \{\mu_J\}_{J\in\mathbf{P}_{V}} \subset \mathbb{Z}[\frac{1}{p}]$ are equivalent.

($i_n$) $\sum\limits_{J\in\mathbf{P}_{V}}\mu_J=0$ and $\lambda_H=q^{n-(N-1)}\sum\limits_{\substack{J\in\mathbf{P}_{V}\\J\subset H}}\mu_J$ for any $H\in\mathbf{P}_{V^*}$.

($ii_n$) $\sum\limits_{H\in\mathbf{P}_{V^*}}\lambda_H=0$ and $\mu_J=q^{1-n}\sum\limits_{\substack{H\in\mathbf{P}_{V^*}\\H\supset J}}\lambda_H$ for any $J\in\mathbf{P}_{V}$. \qed
\end{Lemma}

For $J\in\mathbf{P}_V,H\in\mathbf{P}_{V^*}$, we denote
\[Z_{V,1,H}=\RToySht_H^1\in C_{V,1}^1, \quad Z_{V,1,J}=\RToySht_{V/J}^0\in C_{V,1}^1,\]
\[Z_{V,N-1,H}=\LToySht_H^{N-1}\in C_{V,N-1}^1, \quad Z_{V,N-1,J}=\LToySht_{V/J}^{N-2}\in C_{V,N-1}^1.\]
When $2\le n\le N-2$, we denote
\[Z_{V,n,H}=\RToySht_H^n\in C_{V,n}^1, \text{ or equivalently, } Z_{V,n,H}=\LToySht_H^n\in C_{V,n}^1,\]
\[Z_{V,n,J}=\RToySht_{V/J}^{n-1}\in C_{V,n}^1,  \text{ or equivalently, } Z_{V,n,J}=\LToySht_{V/J}^{n-1}\in C_{V,n}^1.\]

\begin{Theorem}\label{space of principal rational toy horospherical divisors}
The element $\sum_{H\in\mathbf{P}_{V^*}}\lambda_H\cdot Z_{V,n,H}+\sum_{J\in\mathbf{P}_{V}}\mu_J\cdot Z_{V,n,J}$ belongs to $\im(C^0_{V,n}\to C^1_{V,n})\otimes\mathbb{Z}[\frac{1}{p}]$ if and only if $\{\lambda_H\}_{H\in\mathbf{P}_{V^*}},\{\mu_J\}_{J\in\mathbf{P}_V}$ satisfy condition ($i_n$) (or equivalently condition ($ii_n$)) in \cref{equivalent conditions for coefficients of toy horospherical divisors}.
\end{Theorem}
\begin{proof}
\cref{multiplicity 1 of pullback of toy horospherical divisors by partial Frobeniuses,multiplicity q of pullback of toy horospherical divisors by partial Frobeniuses} show that $(F_{V,n}^+)^*Z_{V,n+1,H}=Z_{V,n,H}$, $(F_{V,n}^+)^*Z_{V,n+1,J}=q\cdot Z_{V,n,J}$ for $1\le n\le N-2, H\in\mathbf{P}_{V^*},J\in\mathbf{P}_{V}$. In view of \cref{Frobenius induces multiplication by q on the short complex}, it suffices to prove the statement in the case $n=N-1$.

We know that the morphism $\pi_{L,N-1}:\LToySht_V^{N-1}\to \ToySht_V^{N-1}$ is the blow-up at the points of $\ToySht_V^{N-1}(\mathbb{F}_q)$. For $H\in \ToySht_V^{N-1}(\mathbb{F}_q)=\mathbf{P}_{V^*}$, the exceptional divisor with center $H$ is $\LToySht_H^{N-1}$. Thus for $J\in\mathbf{P}_{V}$, we have an identity of divisors of $\LToySht_V^{N-1}$
\[(\pi_{L,N-1})^*\ToySht_{V/J}^{N-2}=\LToySht_{V/J}^{N-2}+\sum_{\substack{H\in\mathbf{P}_{V^*}\\H\supset J}}\LToySht_H^{N-1}.\]
Also note that $C^0_{V,N-1}$ consists of rational functions on $\ToySht_V^{N-1}=\mathbb{P}^\vee(V)$ with zeros and poles supported on the hyperplanes $\ToySht_{V/J}^{N-2}(J\in\mathbf{P}_{V})$. Hence the statement is true in the case $n=N-1$.
\end{proof}

\subsection{Degree of partial Frobeniuses}
Recall that we denote $\Flag_V^{i,j}$ be the closed subscheme of $\Grass_V^i\times\Grass_V^j$ which consists of pairs $(M,M')$ such that $M\subset M'$.

For $0\le n\le N-1$, we define a scheme $\Flag_{V,\{2\}}^{n,n+1}$ by the following Cartesian diagram
\begin{equation}\label{Cartesian diagram defining Frobenius twisted partial flag variety}
\begin{tikzcd}
\Flag_{V,\{2\}}^{n,n+1}\arrow[r]\arrow[d]&\Flag_V^{n,n+1}\arrow[d]\\
\Grass_V^{n+1}\arrow[r,"\Fr"]&\Grass_V^{n+1}
\end{tikzcd}
\end{equation}
where the right vertical arrow is the projection. So $\Flag_{V,\{2\}}^{n,n+1}$ parameterizes pairs $(M,M')$ such that $M\subset\Fr^*M'$.

The commutative diagram
\[\begin{tikzcd}
\Flag_V^{n,n+1}\arrow[r,"\Fr"]\arrow[d]&\Flag_V^{n,n+1}\arrow[d]\\
\Grass_V^{n+1}\arrow[r,"\Fr"]&\Grass_V^{n+1}
\end{tikzcd}\]
induces a morphism $f_n:\Flag_V^{n,n+1}\to\Flag_{V,\{2\}}^{n,n+1}$, which is finite.

\begin{Lemma}\label{degree of the morphism from flag variety to Frobenius twisted flag variety}
The morphism $f_n:\Flag_V^{n,n+1}\to\Flag_{V,\{2\}}^{n,n+1}$ has degree $q^n$
\end{Lemma}
\begin{proof}
In commutative diagram (\ref{Cartesian diagram defining Frobenius twisted partial flag variety}), denote $g:\Flag_{V,\{2\}}^{n,n+1}\to\Flag_V^{n,n+1}$. We have $g\bcirc f_n=\Fr_{\Flag_V^{n,n+1}}$. Since $\Fr_{\Grass_V^{n+1}}$ is finite and flat, we have $\deg g=\deg \Fr_{\Grass_V^{n+1}}$. Thus
\[\deg f_n=\frac{\deg \Fr_{\Flag_V^{n,n+1}}}{\deg g}=\frac{\deg \Fr_{\Flag_V^{n,n+1}}}{\deg \Fr_{\Grass_V^{n+1}}}=q^{\dim\Flag_V^{n,n+1}-\dim\Grass_V^{n+1}}=q^n\]
\end{proof}

\begin{Proposition}\label{degree of partial Frobeniuses for toy shtukas}
For $0\le n\le N-1$, we have $\deg F_{V,n}^+=q^n$. For $1\le n\le N$, we have $\deg F_{V,n}^-=q^{N-n}$.
\end{Proposition}
\begin{proof}
The first statement follows from \cref{degree of the morphism from flag variety to Frobenius twisted flag variety} and the Cartesian diagram
\[\begin{tikzcd}
\RToySht_V^n\arrow[d,"F_{V,n}^+"]\arrow[r]&\Flag_V^{n,n+1}\arrow[d,"f_n"]\\
\LToySht_V^{n+1}\arrow[r]&\Flag_{V,\{2\}}^{n,n+1}
\end{tikzcd}\]
where the morphism $\RToySht_V^n\to\Flag_V^{n,n+1}$ sends a right toy shtuka $\Fr^*L\subset L'\supset L$ to the pair $(L,L')$, and the morphism $\LToySht_V^{n+1}\to\Flag_{V,\{2\}}^{n,n+1}$ sends a left toy shtuka $\Fr^*M\supset M'\subset M$ to the pair $(M',M)$.

The second statement follows from duality.
\end{proof}

\section{Tate toy shtukas}
\subsection{Tate linear algebra}\label{(Section)Tate linear algebra}
We recall some basic definitions and facts in Tate linear algebra. See Section 6 of Chapter 2 of~\cite{Lef}, Section 2.7.7 of~\cite{BD04}, Section 3.1 of~\cite{Dr06} and Section 1 of~\cite{Kap} for more details.

We consider topological vector spaces over a discrete field $E$.

For a topological vector space $M$, its dual $M^*$ is by definition the space of all continuous linear functionals $M\to E$.

For a discrete topological vector space $Q$, the topology on its dual $Q^*$ is the weakest one such that the linear functional $\left<v,-\right>:Q^*\to E$ is continuous for all $v\in Q$.

\begin{Definition}
A \emph{Tate space} is a topological vector space isomorphic to $P\oplus Q^*$, where $P$ and $Q$ are discrete.
\end{Definition}

\begin{Remark}\label{Hausdorff property of a Tate space}
The topology on a Tate space is Hausdorff.
\end{Remark}

\begin{Definition}
Let $T$ be a Tate space. A linear subspace $\Lambda\subset T$ is said to be \emph{linearly compact} (resp. \emph{linearly cocompact}) if it is closed and for any open vector subspace $U\subset T$ one has $\dim \Lambda/(\Lambda\cap U)<\infty$ (resp. $\dim T/(\Lambda+U)<\infty$).
\end{Definition}

For a Tate space $T$, we equip its dual $T^*$ with a topology as follows. We require the topology on $T^*$ to be linear, i.e., its open linear subspaces form a basis of open neighborhoods of $0$. A linear subspace of $T^*$ is open if and only if it is the orthogonal complement of a linearly compact linear subspace of $T$.

Note that when $T$ is discrete, the topology on $T^*$ agrees with the one given before.

\begin{Remark}
For a Tate space $T$, its dual $T^*$ is again a Tate space. The canonical map $T\to T^{**}$ is an isomorphism. Any discrete vector space is a Tate space. Any linearly compact topological vector space is a Tate space. Duality interchanges discrete and linearly compact Tate spaces.
\end{Remark}

\begin{Remark}
A Tate space over a finite field is locally compact. A linearly compact Tate space over a finite field is compact.
\end{Remark}

\begin{Definition}
A linear subspace of a Tate space is called a $\emph{c-lattice}$ if it is open and linearly compact. A linear subspace of a Tate space is called a $\emph{d-lattice}$ if it is discrete and linearly cocompact.
\end{Definition}

\begin{Remark}
Suppose $T=P\oplus Q^*$ is a Tate space, where $P$ and $Q$ are discrete. Then $P$ is a d-lattice of $T$, and $Q^*$ is a c-lattice of $T$. Thus there exist c-lattices and d-lattices in every Tate space.
\end{Remark}

\begin{Remark}\label{c-lattices form a basis of neighborhood of 0}
c-lattices of a Tate space form a basis of neighborhood of $0$.
\end{Remark}

\begin{Definition}
For a Tate space $T$ and two c-lattices $L_1,L_2$ of $T$, we define the \emph{relative dimension} $d_{L_1}^{L_2}:=\dim(L_2/L_1\cap L_2)-\dim(L_1/L_1\cap L_2)\in\mathbb{Z}$.
\end{Definition}

\begin{Definition}
A \emph{dimension theory} on a Tate space $T$ is a function
\[d:\{\text{c-lattices of $T$}\}\to\mathbb{Z}\]
such that $d(L_2)-d(L_1)=d_{L_1}^{L_2}$ for any $L_1,L_2$.
\end{Definition}

A dimension theory exists and is unique up to adding an integer. So dimension theories on a Tate space $T$ form a $\mathbb{Z}$-torsor, denoted by $\Dim_T$.

For a d-lattice $I$ and a c-lattice $L$ of a Tate space $T$, we denote
\[\chi(I,L)=\dim(I\cap L)-\dim(T/(I+L)).\]
For any d-lattice $I\subset T$, the function $L\mapsto\chi(I,L)$ is a dimension theory.

\begin{Definition}\label{definition of the dimension of a d-lattice}
For a Tate space $T$, we define a function
\begin{align*}
\dim:\{\text{d-lattices of $T$}\}&\to\Dim_T\\
I&\mapsto(L\mapsto\chi(I,L))
\end{align*}
We call $\dim(I)$ the \emph{dimension} of $I$.
\end{Definition}

We borrow the definition of relative determinant from Section 4 of~\cite{AdCK}.

\begin{Definition}\label{definition of relative determinant}
For a Tate space $T$ and two c-lattices $L_1,L_2$ of $T$, we define their \emph{relative determinant} to be the one-dimensional vector space
\[\det\nolimits_{L_1}^{L_2}:=\det(L_1/L_1\cap L_2)^*\otimes\det(L_2/L_1\cap L_2).\]
\end{Definition}

\subsection{Sato Grassmannians}
Let $E$ be a discrete field.

For an $E$-vector space $V$ and a Tate space $M$ over $E$, we denote
\[M\widehat{\otimes}V:=\varprojlim_{\Lambda}(M/\Lambda)\otimes V,\]
where the projective limit is taken over all c-lattices $\Lambda$ of $M$.

For an $E$-scheme $S$ and a quasi-coherent sheaf $\mathscr{G}$ on $S$, we denote $M\widehat{\otimes}\mathscr{G}$ to be the sheaf of $\mathscr{O}_S$-modules such that
\[(M\widehat{\otimes}\mathscr{G})(U)=\varprojlim_{\Lambda}(M/\Lambda)\otimes\mathscr{G}(U)\]
for all open subset $U\subset S$, where the projective limit is taken over all c-lattices $\Lambda$ of $M$.

Let $T$ be a Tate space over $E$ and fix $n\in \Dim_T$.

\begin{Definition}\label{definition of Sato Grassmannian by Drinfeld}
We define a functor $\Grass_T^n$ from the category of $E$-schemes to the category of sets. For an $E$-algebra $R$, an $R$-point of $\Grass_T^n$ is an $R$-submodule $L\subset T\widehat{\otimes}R$ such that there exists a c-lattice $\Lambda'$ of $T$ such that the morphism $L\to (T/\Lambda')\otimes R$ is injective and its cokernel is a finitely generated projective $R$-module of rank $-n(\Lambda')$.
\end{Definition}

The definition below is due to Kashiwara. (Cf. Section 2 of~\cite{Kash}.) \cref{definition of Sato Grassmannian by Drinfeld} and \cref{definition of Sato Grassmannian by Kashiwara} are equivalent by \cref{equivalence of definitions of Sato Grassmannian by Drinfeld and Kashiwara}.

\begin{Definition}\label{definition of Sato Grassmannian by Kashiwara}
For an $E$-scheme $S$, an $S$-point of $\Grass_T^n$ is an $\mathscr{O}_S$-submodule $\mathscr{F}\subset T\widehat{\otimes}\mathscr{O}_S$ such that locally in the Zariski topology there exists a c-lattice $\Lambda\subset T$ such that $n(\Lambda)=0$ and $\mathscr{F}\to (T/\Lambda)\otimes \mathscr{O}_S$ is an isomorphism.
\end{Definition}

\begin{Lemma}\label{equivalence of definitions of Sato Grassmannian by Drinfeld and Kashiwara}
Let $S=\Spec R$, where $R$ is an $E$-algebra. Let $G_1$ (resp. $G_2$) be the set $\Grass_T^n(S)$, where $\Grass_T^n$ is the functor in \cref{definition of Sato Grassmannian by Drinfeld} (resp. \cref{definition of Sato Grassmannian by Kashiwara}). Then we have a canonical bijection $G_2\to G_1$.
\end{Lemma}
\begin{proof}
Since $S$ is quasi-compact, we get a map $G_2\to G_1$.

For $L\in G_1$ and $s\in S$, we can find a c-lattice $\Lambda$ containing $\Lambda'$ such that $n(\Lambda)=0$ and $L\otimes k(s)\to(T/\Lambda)\otimes k(s)$ is an isomorphism of vector spaces, where $k(s)$ is the residue field of $s$. Let $\mathscr{F}=L\widehat{\otimes}_R\mathscr{O}_S$. Since $((T/\Lambda')\otimes R)/L$ is projective, we see that for the fixed c-lattice $\Lambda\subset T$ and the fixed $\mathscr{O}_S$-submodule $\mathscr{F}$ the condition in \cref{definition of Sato Grassmannian by Kashiwara} is open on $S$. Hence there exists a Zariski neighborhood of $s$ in which $\mathscr{F}\to(T/\Lambda)\otimes\mathscr{O}_S$ is an isomorphism. This shows that $\mathscr{F}\in G_1$. So we get a map $G_1\to G_2$.

It is easy to see that the above two maps are inverse of each other.
\end{proof}

\begin{Proposition}
$\Grass_T^n$ is representable by a separated scheme over $E$.
\end{Proposition}
\begin{proof}
This is Proposition 2.2.1 of~\cite{Kash}.
\end{proof}

\begin{Remark}
If $R$ is a field over $E$, then $\Grass_T^n(R)$ is the set of d-lattices of dimension $n$ of the Tate space $T\widehat{\otimes}R$ over $R$.
\end{Remark}

\subsection{Determinant of a family of d-lattices relative to a c-lattice}
Let $T$ be a Tate space over a field $E$.

We denote
\[\Grass_T=\coprod_{n\in\Dim_T}\Grass_T^n.\]
We see that $\Grass_T$ parameterizes d-lattices of $T$.

\begin{Definition}\label{definition of the determinant of a family of d-lattices relative to a c-lattice}
Let $S$ be a scheme over $E$. For $\mathscr{L}\in\Grass_T(S)$ and a c-lattice $W$ of $T$, the \emph{determinant} of $\mathscr{L}$ relative to $W$, denoted by $\det(\mathscr{L},W)$, is defined to be the invertible sheaf $\det(\mathscr{L}\to(T/W)\otimes\mathscr{O}_S)$ on $S$, where $\mathscr{L}$ is in degree 0.
\end{Definition}

\begin{Remark}\label{base change for relative determinant of a family of d-lattices}
The above definition commutes with base change, i.e., for a morphism of schemes $f:S_1\to S_2$ over $E$, an element $\mathscr{L}\in\Grass_T(S_2)$ and a c-lattice $W$ of $T$, we have a canonical isomorphism $\det(f^*\mathscr{L},W)\cong f^*\det(\mathscr{L},W)$
\end{Remark}

\begin{Remark}\label{change of c-lattices for relative determinant of a family of d-lattices}
Let $S$ be a scheme over $E$. For $\mathscr{L}\in\Grass_T(S)$ and two c-lattice $W_1,W_2$ of $T$, we have a canonical isomorphism $\det(\mathscr{L},W_1)\otimes\det_{W_1}^{W_2}\cong\det(\mathscr{L},W_2)$. In particular, the two invertible sheaves $\det(\mathscr{L},W_1)$ and $\det(\mathscr{L},W_2)$ are isomorphic.
\end{Remark}

\begin{Lemma}\label{canonical isomorphism of relative determinants for short exact sequence}
Let $S$ be a scheme over $E$ and let $W$ be a c-lattice of $T$. For $\mathscr{L}_1,\mathscr{L}_2\in\Grass_T(S)$ such that $\mathscr{L}_1\subset\mathscr{L}_2$, we have a canonical isomorphism
\[\det(\mathscr{L}_1,W)\otimes\det(\mathscr{L}_2/\mathscr{L}_1)\cong\det(\mathscr{L}_2,W)\]
which commutes with base change.\qed
\end{Lemma}

\begin{Lemma}\label{canonical isomorphism of relative determinants for Frobenius pullback}
Suppose the base field $E$ is $\mathbb{F}_q$, and $S$ is a scheme over $E$. For $\mathscr{L}\in\Grass_T(S)$ and a c-lattice $W$ of $T$, we have a canonical isomorphism
\[\det(\Fr_S^*\mathscr{L},W)\cong\det(\mathscr{L},W)^{\otimes q}\]
which commutes with base change.\qed
\end{Lemma}

\subsection{Definition of Tate toy shtukas}
Let $T$ be a nondiscrete noncompact Tate space over $\mathbb{F}_q$. Let $\Dim_T$ denote the dimension torsor of $T$.

\begin{Definition}
A \emph{Tate toy shtuka} for $T$ over an $\mathbb{F}_q$-scheme $S$ of dimension $n\in\Dim_T$ is an element $\mathscr{L}\in\Grass_T^n(S)$ such that the composition
\[\Fr_S^*\mathscr{L}\hookrightarrow\Fr_S^*(T\widehat{\otimes}\mathscr{O}_S)=T\widehat{\otimes}\mathscr{O}_S
\twoheadrightarrow(T\widehat{\otimes}\mathscr{O}_S)/\mathscr{L}\]
has rank at most 1. (In other words, the corresponding morphism $\bigwedge^2\Fr_S^*\mathscr{L}\to\bigwedge^2((T\widehat{\otimes}\mathscr{O}_S)/\mathscr{L})$ is zero.)
\end{Definition}

For $n\in \Dim_T$, let $\ToySht_T^n$ be the functor which associates to each $\mathbb{F}_q$-scheme $S$ the set of isomorphism classes of Tate toy shtukas for $T$ over $S$ of dimension $n$.

As in the finite dimensional cases, $\ToySht_T^n$ is representable by a closed subscheme of $\Grass_T^n$.

We denote $\oToySht_T^n:=\ToySht_T^n-\Grass_T^n(\mathbb{F}_q)$. As in \cref{description of the nontrivial locus} we know that $\oToySht_T^n$ is an open subscheme of $\ToySht_T^n$.

\section{Open charts of the scheme of Tate toy shtukas}
\subsection{Notation}\label{(Section)notation for open charts of the scheme of Tate toy shtukas}
For a finite dimensional vector space $V$ over $\mathbb{F}_q$ and two subspaces $V'\subset V''$ of $V$, denote $\Grass_V^{n,V',V''}$ to be the open subscheme of $\Grass_V^n$, such that for any $\mathbb{F}_q$-algebra $R$, $\Grass_V^{n,V',V''}(R)$ consists of $L\in\Grass_V^n(R)$ satisfying the following two conditions:

(i) the morphism $L\to (V/V')\otimes R$ is injective and its cokernel is projective;

(ii) the morphism $L\to (V/V'')\otimes R$ is surjective.

We have a morphism $\Grass_V^{n,V',V''}\to \Grass_{V''/V'}^{n-\dim V/V''}$ which maps $L\in \Grass_V^{n,V',V''}(R)$ to $\im(L''\to(V''/V')\otimes R)\in\Grass_{V''/V'}^{n-\dim V/V''}(R)$, where $L''=\ker(L\to(V/V'')\otimes R)$.

Let $\oGrass_V^{n,V',V''}$ be the fiber product
\begin{equation}\label{Cartesian diagram for oGrass}
\begin{tikzcd}
\oGrass_V^{n,V',V''}\arrow[r]\arrow[d,"f"]&\Grass_V^{n,V',V''}\arrow[d]\\
\Grass_{V''/V'}^{n-\dim V/V''}-\Grass_{V''/V'}^{n-\dim V/V''}(\mathbb{F}_q)\arrow[r]&\Grass_{V''/V'}^{n-\dim V/V''}\
\end{tikzcd}
\end{equation}

We denote $\oToySht_V^{n,V',V''}=\oToySht_V^n\cap\oGrass_V^{n,V',V''}$. We have a morphism $\oToySht_V^{n,V',V''}\to\oToySht_{V''/V'}^{n-\dim V/V''}$ induced by $f$.

\begin{Lemma}
For a finite dimensional vector space $V$ over $\mathbb{F}_q$ and two subspaces $V'\subset V''$ of $V$, the morphism $\Grass_V^{n,V',V''}\to\Grass_{V''/V'}^{n-\dim V/V''}$ is affine.\qed
\end{Lemma}

\begin{Lemma}\label{transition maps between open charts of Grassmannians are affine}
For finite dimensional $\mathbb{F}_q$-vectors spaces $V'_2\subset V'_1\subset V'_0\subset V''_0\subset V''_1\subset V''_2$, the morphism $\oGrass_{V''_2/V'_2}^{n,V'_0/V'_2,V''_0/V'_2}\to\oGrass_{V''_1/V'_1}^{n-\dim V''_2/V''_1,V'_0/V'_1,V''_0/V'_1}$ is affine. \qed
\end{Lemma}

For a nondiscrete noncompact Tate space $T$ over $\mathbb{F}_q$, two c-lattices $\Lambda'\subset\Lambda''$ of $T$ and $n\in\Dim_T$, we define the notation $\Grass_T^{n,\Lambda',\Lambda''}$, $\oGrass_T^{n,\Lambda',\Lambda''}$ and $\oToySht_T^{n,\Lambda',\Lambda''}$ similarly.

\subsection{Admissible pairs of c-lattices}\label{(Section)admissible pairs of c-lattices}
Let $T$ be a nondiscrete noncompact Tate space over $\mathbb{F}_q$.

For two pairs of c-lattices $(\widetilde{\Lambda'},\widetilde{\Lambda''})$ and $(\Lambda',\Lambda'')$ of $T$, we say that $(\widetilde{\Lambda'},\widetilde{\Lambda''})$ is \emph{greater than} $(\Lambda',\Lambda'')$ (denoted by $(\widetilde{\Lambda'},\widetilde{\Lambda''})\succ(\Lambda',\Lambda'')$), if $\widetilde{\Lambda'}\subset\Lambda'$ and $\Lambda''\subset\widetilde{\Lambda''}$. All pairs of c-lattices of $T$ form a directed set under this partial order.

If $(\widetilde{\Lambda'},\widetilde{\Lambda''})\succ(\Lambda',\Lambda'')$, then $\Grass_T^{n,\Lambda',\Lambda''}\subset\Grass_T^{n,\widetilde{\Lambda'},\widetilde{\Lambda''}}$, $\oGrass_T^{n,\Lambda',\Lambda''}\subset\oGrass_T^{n,\widetilde{\Lambda'},\widetilde{\Lambda''}}$, $\oToySht_T^{n,\Lambda',\Lambda''}\subset\oToySht_T^{n,\widetilde{\Lambda'},\widetilde{\Lambda''}}$.

A pair of c-lattices $(\Lambda',\Lambda'')$ of $T$ is said to be \emph{admissible with respect to $n\in\Dim_T$} if $n(\Lambda')\le -2, n(\Lambda'')\ge 2$, and $\Lambda'\subset\Lambda''$. Let $AP_n(T)$ denote the set of pairs of c-lattices of $T$ that are admissible with respect to $n$. It is a directed set with respect to the partial order above. We see that $\oToySht_T^n$ is covered by the union of $\oToySht_T^{n,\Lambda',\Lambda''}$ for all $(\Lambda',\Lambda'')\in AP_n(T)$.

\subsection{$\oToySht_T^{n,\Lambda',\Lambda''}$ as a projective limit}\label{(Section)open subscheme of Tate toy shtukas as a projective limit}
Let $T$ be a nondiscrete noncompact Tate space over $\mathbb{F}_q$. Let $n\in\Dim_T$ and $(\Lambda',\Lambda'')\in AP_n(T)$.

In this subsection we describe the open subscheme $\oToySht_T^{n,\Lambda',\Lambda''}$ of $\oToySht_T^n$ as a projective limit of schemes.

For any $(\widetilde{\Lambda'},\widetilde{\Lambda''})\succ(\Lambda',\Lambda'')$, we denote $U_{\widetilde{\Lambda'},\widetilde{\Lambda''}}^{n,\Lambda',\Lambda''}
=\oToySht_{\widetilde{\Lambda''}/\widetilde{\Lambda'}}^{n(\widetilde{\Lambda''}),\Lambda'/\widetilde{\Lambda'}, \Lambda''/\widetilde{\Lambda'}}$. In particular, $U_{\Lambda',\Lambda''}^{n,\Lambda',\Lambda''}=\oToySht_{\Lambda''/\Lambda'}^{n(\Lambda'')}$. We have a morphism $\oToySht_T^{n,\Lambda',\Lambda''}\to U_{\widetilde{\Lambda'},\widetilde{\Lambda''}}^{n,\Lambda',\Lambda''}$ induced by the morphism $\oGrass_T^{n,\Lambda',\Lambda''}\to\oGrass_{\widetilde{\Lambda''}/\widetilde{\Lambda'}}^{n(\widetilde{\Lambda''}), \Lambda'/\widetilde{\Lambda'},\Lambda''/\widetilde{\Lambda'}}$. For any $(\Lambda'_2,\Lambda''_2)\succ(\Lambda'_1,\Lambda''_1)\succ(\Lambda',\Lambda'')$, we have a transition map $U_{\Lambda'_2,\Lambda''_2}^{n,\Lambda',\Lambda''}\to U_{\Lambda'_1,\Lambda''_1}^{n,\Lambda',\Lambda''}$, which is affine by \cref{transition maps of open charts of toy shtukas are affine}.

\begin{Lemma}\label{open subscheme of Tate toy shtukas as a projective limit}
The morphisms $\oToySht_T^{n,\Lambda',\Lambda''}\to U_{\widetilde{\Lambda'},\widetilde{\Lambda''}}^{n,\Lambda',\Lambda''}$ induce an isomorphism
\[\oToySht_T^{n,\Lambda',\Lambda''}\xrightarrow{\sim}
\varprojlim_{(\widetilde{\Lambda'},\widetilde{\Lambda''})\succ(\Lambda',\Lambda'')} U_{\widetilde{\Lambda'},\widetilde{\Lambda''}}^{n,\Lambda',\Lambda''}.\]
\end{Lemma}

\subsection{Transition maps are affine}
\begin{Lemma}\label{source is affine over a base implies morphism is affine}
Let $f:Y_1\to Y_2$ be a morphism of schemes over a scheme $S$. If $Y_1$ is affine over $S$ and $Y_2$ is separated, then $f$ is affine. \qed
\end{Lemma}

\begin{Lemma}\label{transition maps of toy shtukas are affine}
Let $V'_2\subset V'_1\subset V'_0\subset V''_0\subset V''_1\subset V''_2$ be finite dimensional vector spaces over $\mathbb{F}_q$. Then the morphism $\oToySht_{V''_2/V'_2}^{n,V'_0/V'_2,V''_0/V'_2}\to \oToySht_{V''_1/V'_1}^{n-\dim V''_2/V''_1,V'_0/V'_1,V''_0/V'_1}$ is affine.
\end{Lemma}
\begin{proof}
The morphism $\oToySht_{V''_2/V'_2}^{n,V'_0/V'_2,V''_0/V'_2}\to\oGrass_{V''_2/V'_2}^{n,V'_0/V'_2,V''_0/V'_2}$ is affine since it is a closed immersion. The morphism $\oGrass_{V''_2/V'_2}^{n,V'_0/V'_2,V''_0/V'_2}\to\oGrass_{V''_1/V'_1}^{n-\dim V''_2/V''_1,V'_0/V'_1,V''_0/V'_1}$ is affine by \cref{transition maps between open charts of Grassmannians are affine}. Applying \cref{source is affine over a base implies morphism is affine} with $Y_1=\oToySht_{V''_2/V'_2}^{n,V'_0/V'_2,V''_0/V'_2}, Y_2=\oToySht_{V''_1/V'_1}^{n-\dim V''_2/V''_1,V'_0/V'_1,V''_0/V'_1}$ and $S=\oGrass_{V''_1/V'_1}^{n-\dim V''_2/V''_1,V'_0/V'_1,V''_0/V'_1}$, the statement follows.
\end{proof}

\begin{Lemma}\label{transition maps of open charts of toy shtukas are affine}
Let $T$ be a nondiscrete noncompact Tate space over $\mathbb{F}_q$. Let $n\in\Dim_T$ and $(\Lambda'_0,\Lambda''_0)\in AP_n(T)$. Let $(\Lambda'_1,\Lambda''_1),(\Lambda'_2,\Lambda''_2)$ be two c-lattices such that $(\Lambda'_2,\Lambda''_2)\succ(\Lambda'_1,\Lambda''_1)\succ(\Lambda'_0,\Lambda''_0)$. Then the morphism $U_{\Lambda'_2,\Lambda''_2}^{n,\Lambda'_0,\Lambda''_0}\to U_{\Lambda'_1,\Lambda''_1}^{n,\Lambda'_0,\Lambda''_0}$ is affine.
\end{Lemma}
\begin{proof}
The statement follows from \cref{transition maps of toy shtukas are affine}.
\end{proof}

\subsection{Transition maps are smooth}
\begin{Lemma}\label{transition maps of toy shtukas for subspaces are smooth}
Let $V$ be a finite dimensional vector space over $\mathbb{F}_q$. Let $W$ be a subspace of $V$. Then the morphism $\oToySht_V^{n,0,W}\to\oToySht_W^{n-\dim V/W}$ is smooth.
\end{Lemma}
\begin{proof}
$\Grass_W^{n-\dim V/W}$ is covered by $\Grass_W^{n-\dim V/W,W',W'}$ when $W'$ runs through all subspaces of $W$ of dimension $(\dim V-n)$. So in view of the Cartesian diagram
\[\begin{tikzcd}
  \Grass_V^{n,W',W'}\arrow[r]\arrow[d] &\Grass_V^{n,0,W}\arrow[d] \\
  \Grass_W^{n-\dim V/W,W',W'}\arrow[r] &\Grass_W^{n-\dim V}
\end{tikzcd}\]
it suffices to show that the morphism
\[f:\oToySht_V^{n,0,W}\cap\Grass_V^{n,W',W'}\to\oToySht_W^{n-\dim V/W}\cap\Grass_W^{n-\dim V/W,W',W'}\]
is smooth for each subspace $W'\subset W$ of dimension $(\dim V-n)$.

Fix one such $W'$. Choose splittings $W=W'\oplus W'', V=W\oplus V'$.

We define $\mathbb{F}_q$-schemes $M,M',M''$ where
\[M=\{a\in\underline{\Hom}(W''\oplus V',W')|\rank a=\rank(W''\xrightarrow{a|_{W''}}W')=1\}\]
\[M'=\underline{\Hom}(V',W'),\quad M''=\underline{\Hom}^{\rank=1}(W'',W')\]
So $M$ is a locally closed subscheme of $\underline{\Hom}(W''\oplus V',W')$ and $M''$ is a locally closed subscheme of $\underline{\Hom}(W'',W')$.

We denote Artin-Schreier maps
\begin{align*}
AS&=\Id-\Fr:\underline{\Hom}(W''\oplus V',W')\to\underline{\Hom}(W''\oplus V',W'),\\
AS'&=\Id-\Fr:\underline{\Hom}(V',W')\to\underline{\Hom}(V',W'),\\
AS''&=\Id-\Fr:\underline{\Hom}(W'',W')\to\underline{\Hom}(W'',W').
\end{align*}
From the explicit local description of $\ToySht$ in \cref{(Section)explicit local discription of the scheme of toy shtukas} and the Cartesian diagram (\ref{Cartesian diagram for oGrass}), we know that
\[\oToySht_V^{n,0,W}\cap\Grass_V^{n,W',W'}=AS^{-1}(M),\]
\[\oToySht_W^{n-\dim V/W}\cap\Grass_W^{n-\dim V/W,W',W'}=(AS'')^{-1}(M''),\]
and that the morphism $f:AS^{-1}(M)\to (AS'')^{-1}(M'')$ is induced by the projection $\underline{\Hom}(W''\oplus V',W')\to\underline{\Hom}(W'',W')$.

We denote
\[\mathbb{P}_{W'}\star M'=\{(L,A)\in\mathbb{P}_{W'}\times M'|L\supset\im A\}\]
to be the closed subscheme of $\mathbb{P}_{W'}\times M'$.

We have to prove that the morphism $f:AS^{-1}(M)\to (AS'')^{-1}(M'')$ is smooth. To this end, we will construct a Cartesian diagram
\[\begin{tikzcd}
AS^{-1}(M)\arrow[d,"f"]\arrow[r,"g"] &\mathbb{P}_{W'}\star M'\arrow[d,"h"]\\
(AS'')^{-1}(M'') \arrow[r,"u"] &\mathbb{P}_{W'}
\end{tikzcd}\]
and prove that $h$ is smooth.

The maps in the diagram are as follows. The morphism $h$ is the composition $\mathbb{P}_{W'}\star M'\xrightarrow{\Id\times AS'}
\mathbb{P}_{W'}\star M'\xrightarrow{\pr}\mathbb{P}_{W'}$. The composition $AS^{-1}(M)\hookrightarrow\underline{\Hom}(W''\oplus V', W')\xrightarrow{\operatorname{res}_{W''}}M''\xrightarrow{\im}\mathbb{P}_{W'}$ and the projection $AS^{-1}(M)\to M'$ induce a morphism $AS^{-1}(M)\to\mathbb{P}_{W'}\times M'$ which factors through $\mathbb{P}_{W'}\star M'$. This gives $g$. The morphism $u$ is the composition $(AS'')^{-1}(M'') \xrightarrow{AS''}M'' \xrightarrow{\im}\mathbb{P}_{W'}$.

One can check that the diagram is commutative and Cartesian.

The morphism $\mathbb{P}_{W'}\star M'\xrightarrow{\Id\times AS'}
\mathbb{P}_{W'}\star M'$ is smooth since $AS'$ is \'etale. The morphism $\mathbb{P}_{W'}\star M'\xrightarrow{\pr}\mathbb{P}_{W'}$ is smooth since it is the projection morphism for a vector bundle. So $h$ is smooth. Hence is $f$.
\end{proof}

The following statement is dual to \cref{transition maps of toy shtukas for subspaces are smooth}.

\begin{Lemma}\label{transition maps of toy shtukas for quotients are smooth}
Let $V$ be a finite dimensional vector space over $\mathbb{F}_q$. Let $W$ be a subspace of $V$. Then the morphism $\oToySht_V^{n,W,V}\to\oToySht_{V/W}^n$ is smooth. \qed
\end{Lemma}

\begin{Proposition}\label{transition maps of toy shtukas for subquotients are smooth}
Let $V'\subset V''\subset V$ be finite dimensional vector spaces over $\mathbb{F}_q$. Then the morphism $\oToySht_V^{n,V',V''}\to\oToySht_{V''/V'}^{n-\dim V/V''}$ is smooth. \qed
\end{Proposition}
\begin{proof}
The statement follows from \cref{transition maps of toy shtukas for subspaces are smooth} and \cref{transition maps of toy shtukas for quotients are smooth}.
\end{proof}

\begin{Corollary}\label{transition maps of toy shtukas are smooth}
With the same notation and assumptions of \cref{transition maps of open charts of toy shtukas are affine}, the morphism $U_{\Lambda'_2,\Lambda''_2}^{n,\Lambda'_0,\Lambda''_0}\to U_{\Lambda'_1,\Lambda''_1}^{n,\Lambda'_0,\Lambda''_0}$ is smooth.
\end{Corollary}
\begin{proof}
The statement follows from \cref{transition maps of toy shtukas for subquotients are smooth}.
\end{proof}

\section{Functorial properties of toy horospherical divisors}
In this section, we use the notation of \cref{(Section)notation for open charts of the scheme of Tate toy shtukas,(Section)admissible pairs of c-lattices}.

\subsection{Notation}\label{(Section)notation for functorial properties of toy horospherical divisors}
Suppose $M$ is a vector space over $\mathbb{F}_q$ such that $3\le\dim M<\infty$.

For $J\in\mathbf{P}_M$, we consider $\ToySht_{M/J}^{n-1}$ as a closed subscheme of $\ToySht_M^n$ as before. Denote $\Delta_{M,J}^n=\ToySht_{M/J}^{n-1}\cap\oToySht_M^n$. For two subspaces $M'\subset M''$ of $M$, denote $\Delta_{M,J}^{n,M',M''}=\Delta_{M,J}^n\cap\oGrass_M^{n,M',M''}$.

For $H\in\mathbf{P}_{M^*}$, we consider $\ToySht_H^n$ as a closed subscheme of $\ToySht_M^n$ as before. Denote $\Delta_{M,H}^n=\ToySht_H^n\cap\oToySht_M^n$. For two subspaces $M'\subset M''$ of $M$, denote $\Delta_{M,H}^{n,M',M''}=\Delta_{M,H}^n\cap\oGrass_M^{n,M',M''}$.

We denote $\Delta_M^n=(\bigcup_{H\in\mathbf{P}_{M^*}}\Delta_{M,H}^n)\cup(\bigcup_{J\in\mathbf{P}_M}\Delta_{M,J}^n)$ and $\Delta_M^{n,M',M''}=\Delta_M^n\cap\oGrass_M^{n,M',M''}$. They are the union of toy horospherical divisors of $\oToySht_M^n$ and $\oToySht_M^{n,M',M''}$ respectively.

By \cref{irreducibility of the scheme of toy shtukas}, each $\Delta_{M,J}^n, \Delta_{M,J}^{n,M',M''},\Delta_{M,H}^n, \Delta_{M,H}^{n,M',M''}$ is reduced and irreducible if it is nonempty.

\subsection{Pullbacks of toy horospherical divisors under transition maps}
In this subsection, we calculate the pullback of a toy horospherical divisor under transition maps. (The fact that this pullback is well-defined is clear from \cref{transition maps of toy shtukas are smooth}.)

\begin{Lemma}\label{set-theoretic containment for the pullback for quotient of toy horospherical divisor}
Let $W\subset V$ be finite dimensional vector spaces over $\mathbb{F}_q$. Let $J'\in\mathbf{P}_{V/W}$ and denote $\mathfrak{J}=\{J\in\mathbf{P}_V|\im(J\to(V/W))=J'\}$. Let $E$ be a field over $\mathbb{F}_q$. Suppose $L\in\oToySht_V^{n,W,V}(E)$ satisfies $\im(L\to((V/W)\otimes E))\supset J'\otimes E$. Then $L\supset J\otimes E$ for some $J\in\mathfrak{J}$.
\end{Lemma}
\begin{proof}
Let $J''$ be the unique subspace of $V$ containing $W$ such that $J''/W=J'$. Denote $L'=L\cap((W+J'')\otimes E), L''=\im(L\to(V/(W+J''))\otimes E),L'''=\im(L\to(V/W)\otimes E)$. Since $L'''$ is a nontrivial toy shtuka and $L'''\supset J'\otimes E$, $L''$ is also a nontrivial toy shtuka. Applying \cref{triviality of intersection or quotient of intersection} to $L$ and $W+J''$, we get $\Fr_E^*L'=L'$. Since $L\in\oToySht_V^{n,W,V}(E)$, we have $L\cap(W\otimes E)=0$. Hence $\dim_E L'=1$. Thus $L'=J\otimes E$ for some $J\in\mathfrak{J}$.
\end{proof}

\begin{Lemma}\label{pullback of toy horospherical divisors for quotient space}
Let $V$ be a finite dimensional vector space over $\mathbb{F}_q$. Let $W$ be a subspace of $V$ such that $\dim(V/W)\ge 3$. Let $u:\oToySht_V^{n,W,V}\to\oToySht_{V/W}^n$ be the transition map. Let $H'\in\mathbf{P}_{(V/W)^*}$ and let $H$ be the hyperplane of $V$ such that $H\supset W$ and $H/W=H'$. Let $J'\in\mathbf{P}_{V/W}$ and denote $\mathfrak{J}=\{J\in\mathbf{P}_V|\im(J\to(V/W))=J'\}$.

When $1\le n\le \dim(V/W)-2$, we have an equality of divisors $u^*(\Delta_{V/W,H'}^n)=\Delta_{V,H}^{n,W,V}$. When $2\le n\le \dim(V/W)-1$, we have an equality of divisors $u^*(\Delta_{V/W,J'}^n)=\sum_{J\in\mathfrak{J}}\Delta_{V,J}^{n,W,V}$
\end{Lemma}
\begin{proof}
Since $u$ is smooth by \cref{transition maps of toy shtukas are smooth}, and $\Delta_{V/W,H'}^n,\Delta_{V,H}^{n,W,V},\Delta_{V/W,J'}^n,\Delta_{V,J}^{n,W,V}(J\in\mathfrak{J})$ are reduced by \cref{basic properties of the scheme of toy shtukas}, it suffices to prove the statements set-theoretically.

When $1\le n\le \dim(V/W)-2$, it is obvious that $u^{-1}(\Delta_{V/W,H'}^n)=\Delta_{V,H}^{n,W,V}$ as subsets of $\oToySht_V^{n,W,V}$.

When $2\le n\le \dim(V/W)-1$, we have $u(\Delta_{V,J}^{n,W,V})\subset\Delta_{V/W,J'}^n$ for $J\in\mathfrak{J}$, and the inclusion $u^{-1}(\Delta_{V/W,J'}^n)\subset\bigcup_{J\in\mathfrak{J}}\Delta_{V,J}^{n,W,V}$ is set-theoretically true by \cref{set-theoretic containment for the pullback for quotient of toy horospherical divisor}.
\end{proof}

\begin{Lemma}\label{pullback of toy horospherical divisors for subspace}
Let $V$ be a finite dimensional vector space over $\mathbb{F}_q$. Let $W$ be a subspace of $V$ such that $\dim W\ge 3$. Let $u:\oToySht_V^{n,0,W}\to\oToySht_W^{n-\dim(V/W)}$ be the transition map. Let $J\in\mathbf{P}_W$, $H'\in\mathbf{P}_{W^*}$ and denote $\mathfrak{H}=\{H\in\mathbf{P}_{V^*}|H\cap W=H'\}$.

When $\dim V/W+2\le n\le \dim V-1$, we have an equality of divisors $u^*(\Delta_{W,J}^{n-\dim V/W})=\Delta_{V,J}^{n,0,W}$. When $\dim V/W+1\le n\le \dim V-2$, we have an equality of divisors $u^*(\Delta_{W,H'}^{n-\dim V/W})=\sum_{H\in\mathfrak{H}}\Delta_{V,H}^{n,0,W}$.
\end{Lemma}
\begin{proof}
The statement is dual to \cref{pullback of toy horospherical divisors for quotient space}.
\end{proof}

\begin{Proposition}\label{pullback of toy horospherical divisors under transition maps}
Let $T$ be a nondiscrete noncompact Tate space over $\mathbb{F}_q$. Let $n\in\Dim_T$ and let $(\Lambda'_0,\Lambda''_0)\in AP_n(T)$. Let $(\Lambda'_1,\Lambda''_1),(\Lambda'_2,\Lambda''_2)$ be two pairs of c-lattices of $T$ such that $(\Lambda'_2,\Lambda''_2)\succ(\Lambda'_1,\Lambda''_1)\succ(\Lambda'_0,\Lambda''_0)$. Let $u: U_{\Lambda'_2,\Lambda''_2}^{n,\Lambda'_0,\Lambda''_0}\to U_{\Lambda'_1,\Lambda''_1}^{n,\Lambda'_0,\Lambda''_0}$ be the transition map as in \cref{(Section)open subscheme of Tate toy shtukas as a projective limit}.

For $J_1\in\mathbf{P}_{\Lambda''_1/\Lambda'_1}$ such that $J_1\not\subset\Lambda'_0/\Lambda'_1$, put $\mathfrak{J}_2=\{J_2\in\mathbf{P}_{\Lambda''_2/\Lambda'_2}|\im(J_2\cap\Lambda''_1\to\Lambda''_1/\Lambda'_1)=J_1\}$. For $H_1\in\mathbf{P}_{(\Lambda''_1/\Lambda'_1)^*}$ such that $H_1\not\supset\Lambda''_0/\Lambda'_1$, put $\mathfrak{H}_2=\{H_2\in\mathbf{P}_{(\Lambda''_2/\Lambda'_2)^*}|\im(H_2\cap\Lambda''_1\to\Lambda''_1/\Lambda'_1)=H_1\}$.

Then we have equalities of divisors
\[u^*(\Delta_{\Lambda''_1/\Lambda'_1,J_1}^{n(\Lambda''_1)}\cap U_{\Lambda'_1,\Lambda''_1}^{n,\Lambda'_0,\Lambda''_0}) =\sum_{J_2\in\mathfrak{J}_2}\Delta_{\Lambda''_2/\Lambda'_2,J_2}^{n(\Lambda''_2)}\cap U_{\Lambda'_2,\Lambda''_2}^{n,\Lambda'_0,\Lambda''_0},  \]
\[u^*(\Delta_{\Lambda''_1/\Lambda'_1,H_1}^{n(\Lambda''_1)}\cap U_{\Lambda'_1,\Lambda''_1}^{n,\Lambda'_0,\Lambda''_0}) =\sum_{H_2\in\mathfrak{H}_2}\Delta_{\Lambda''_2/\Lambda'_2,H_2}^{n(\Lambda''_2)}\cap U_{\Lambda'_2,\Lambda''_2}^{n,\Lambda'_0,\Lambda''_0}.  \]
\end{Proposition}
\begin{proof}
The statement follows from \cref{pullback of toy horospherical divisors for quotient space,pullback of toy horospherical divisors for subspace,Cartesian diagram for transition maps of toy shtukas and open charts}.
\end{proof}

\begin{Lemma}\label{Cartesian diagram for transition maps of toy shtukas and open charts}
With the same notation as in \cref{pullback of toy horospherical divisors under transition maps}, the following commutative diagram is Cartesian. \qed
\end{Lemma}
\[\begin{tikzcd}
U_{\Lambda'_2,\Lambda''_2}^{n,\Lambda'_0,\Lambda''_0}\arrow[r]\arrow[d]
&\oToySht_{\Lambda''_2/\Lambda'_2}^{n(\Lambda''_2),\Lambda'_2/\Lambda'_1,\Lambda''_2/\Lambda'_1}\arrow[d]\\
 U_{\Lambda'_1,\Lambda''_1}^{n,\Lambda'_0,\Lambda''_0}\arrow[r] &\oToySht_{\Lambda''_1/\Lambda'_1}^{n(\Lambda''_1)}
\end{tikzcd}\]

\begin{Lemma}\label{transition maps between toy horospherical divisors are dominant}
  Use the same notation of \cref{pullback of toy horospherical divisors under transition maps}. Let $J_2\in\mathfrak{J}_2$. Put $Y_1=\Delta_{\Lambda''_1/\Lambda'_1,J_1}^{n(\Lambda''_1)}\cap U_{\Lambda'_1,\Lambda''_1}^{n,\Lambda'_0,\Lambda''_0}$ and $Y_2=\Delta_{\Lambda''_2/\Lambda'_2,J_2}^{n(\Lambda''_2)}\cap U_{\Lambda'_2,\Lambda''_2}^{n,\Lambda'_0,\Lambda''_0}$. Then we have $u(Y_2)\subset Y_1$, and the morphism $Y_2\to Y_1$ induced by $u$ is dominant.
\end{Lemma}
\begin{proof}
  It is clear that $u(Y_2)\subset Y_1$. Since $(\Lambda'_0,\Lambda''_0)\in AP_n(T)$, both $Y_1$ and $Y_2$ are nonempty. Hence they are irreducible by \cref{irreducibility of the scheme of toy shtukas}. Then \cref{pullback of toy horospherical divisors under transition maps} shows that $Y_2$ is an irreducible component of $u^{-1}(Y_1)$. Now the statement follows from \cref{transition maps of toy shtukas are smooth}.
\end{proof}

\subsection{Functoriality for non-horospherical toy shtukas}
\begin{Lemma}\label{transition maps of toy shtukas are surjective}
Let $V'\subset V''\subset V$ be finite dimensional vector spaces over $\mathbb{F}_q$. Then the morphism $\oToySht_V^{n,V',V''}\to\oToySht_{V''/V'}^{n-\dim V/V''}$ is surjective.
\end{Lemma}
\begin{proof}
Fix a splitting $V=V'\oplus V''/V'\oplus V/V''$. Let $L\in\oToySht_{V''/V'}^{n-\dim V/V''}(E)$, where $E$ is a field over $\mathbb{F}_q$. Then $\widetilde{L}=L\oplus((V/V'')\otimes E)\subset V\otimes E$ is an $E$-point of $\oToySht_V^{n,V',V''}$ which maps to $L$.
\end{proof}

\begin{Lemma}\label{smoothness of restriction of transition map to extra toy horospherical divisors}
Let $V'\subset V''\subset V$ be finite dimensional vector spaces over $\mathbb{F}_q$ such that $\dim V''/V'\ge 3$. Suppose $J\in\mathbf{P}_V$ satisfies $J\not\subset V''$. Then the composition $\Delta_{V,J}^{n,V',V''}\hookrightarrow\oToySht_V^{n,V',V''}\to\oToySht_{V''/V'}^{n-\dim V/V''}$ is smooth.
\end{Lemma}
\begin{proof}
We choose a subspace $V'''\subset V$ such that $V'''\supset V''$ and $V'''\oplus J=V$. Then we get an isomorphism $\oToySht_{V'''}^{n-1,V',V''}\xrightarrow{\sim}\Delta_{V,J}^{n,V',V''}$. The composition $\oToySht_{V'''}^{n-1,V',V''}\xrightarrow{\sim}\Delta_{V,J}^{n,V',V''}\hookrightarrow\oToySht_V^{n,V',V''}\to\oToySht_{V''/V'}^{n-\dim V/V''}$ is the natural one, and is smooth by \cref{transition maps of toy shtukas for subquotients are smooth}. The statement follows.
\end{proof}

Recall that the notation $\ooToySht_V^n$ is defined in \cref{(Section)short complexes related with principal toy horospherical divisors}.

\begin{Lemma}\label{complement of open chart in the scheme of toy shtukas is horospherical}
Let $V'\subset V''\subset V$ be finite dimensional vector spaces over $\mathbb{F}_q$ such that $\dim V''/V'\ge 3$. Assume $\dim V/V''+1\le n\le \dim V/V'-1$. Then $\ooToySht_V^n$ is contained in $\oToySht_V^{n,V',V''}$.
\end{Lemma}
\begin{proof}
Let $L\in\ooToySht_V^n(E)$, where $E$ is a field over $\mathbb{F}_q$. From the description of Schubert divisors in \cref{description of Schubert divisors of the scheme of toy shtukas}, we see that $L\cap(W\otimes E)=0$ for any subspace $W\subset V$ of codimension $n$. Thus $\dim_E(V\cap(U\otimes E))=\max\{0,n-\dim V/U\}$ for any subspace $U\subset V$. In particular, we have $L\in\Grass_V^{n,V',V''}(E)$.

Let $\overline{L}=\im(L\cap(V''\otimes E)\to(V''/V')\otimes E)$. Then
\begin{equation}\label{equation: minimal possible dimension of intersection}
\dim_E(\overline{L}\cap(M\otimes E))=\max\{0,n-\dim V/V'+\dim M\}
\end{equation}
for any subspace $M\subset V''/V'$. Suppose $\overline{L}=M\otimes E$ for some subspace $M\subset V''/V'$. Then $\dim_E(\overline{L}\cap(M\otimes E))=\dim_E\overline{L}=n-\dim V/V''$. On the other hand, $\dim_E(\overline{L}\cap(M\otimes E))=\dim M$. We get a contradiction to (\ref{equation: minimal possible dimension of intersection}) from the assumption $\dim V/V''+1\le n\le \dim V/V'-1$. Thus $\overline{L}$ is a nontrivial toy shtuka for $V''/V'$.
\end{proof}

\begin{Corollary}\label{complement of (Tate) open chart in the scheme of toy shtukas is horospherical}
For $(\Lambda',\Lambda''),(\widetilde{\Lambda'},\widetilde{\Lambda''})\in AP_n(T)$ such that $(\widetilde{\Lambda'},\widetilde{\Lambda''})\succ(\Lambda',\Lambda'')$, $\ooToySht_{\widetilde{\Lambda''}/\widetilde{\Lambda'}}^{n(\widetilde{\Lambda''})}$ is contained in $U_{\widetilde{\Lambda'},\widetilde{\Lambda''}}^{n,\Lambda',\Lambda''}$. \qed
\end{Corollary}

Recall that $\Delta_V^n$ denotes the union of all toy horospherical divisors of $\oToySht_V^n$.

\begin{Lemma}\label{set-theoretic containment for the pullback of toy horospherical divisor}
Let $V'\subset V''\subset V$ be finite dimensional vector spaces over $\mathbb{F}_q$ such that $\dim V''/V'\ge 3$. Assume $\dim V/V''+1\le n\le \dim V/V'-1$. Then the inverse image of $\Delta_{V''/V'}^{n-\dim V/V''}$ under the morphism $\oToySht_V^{n,V',V''}\to\oToySht_{V''/V'}^{n-\dim V/V''}$ is set-theoretically contained in $\Delta_V^n$.
\end{Lemma}
\begin{proof}
The statement follows from \cref{set-theoretic containment for the pullback for quotient of toy horospherical divisor} and the dual of it.
\end{proof}

\begin{Remark}
Let $V'\subset V''\subset V$ be finite dimensional vector spaces over $\mathbb{F}_q$ such that $\dim V''/V'\ge 3$. Assume $\dim V/V''+1\le n\le \dim V/V'-1$. Then \cref{complement of open chart in the scheme of toy shtukas is horospherical} and \cref{set-theoretic containment for the pullback of toy horospherical divisor} show that the morphism $\oToySht_V^{n,V',V''}\to\oToySht_{V''/V'}^{n-\dim V/V''}$ induces a morphism $\ooToySht_V^n\to\ooToySht_{V''/V'}^{n-\dim V/V''}$.
\end{Remark}

\begin{Lemma}\label{properties of transition map for non-horospherical toy shtukas}
Let $V'\subset V''\subset V$ be finite dimensional vector spaces over $\mathbb{F}_q$ such that $\dim V''/V'\ge 3$. Assume $\dim V/V''+1\le n\le \dim V/V'-1$. Then the morphism $\ooToySht_V^n\to\ooToySht_{V''/V'}^{n-\dim V/V''}$ is affine, smooth and surjective.
\end{Lemma}
\begin{proof}
The morphism is affine by \cref{transition maps of toy shtukas are affine} and is smooth by \cref{transition maps of toy shtukas for subquotients are smooth}.

Now we prove that the morphism is surjective.

Denote $u:\oToySht_V^{n,V',V''}\to\oToySht_{V''/V'}^{n-\dim V/V''}$. Let $x$ be a point of $\ooToySht_{V''/V'}^{n-\dim V/V''}$. We show that $u^{-1}(x)\cap\ooToySht_V^n$ is nonempty.

The morphism $u$ is surjective by \cref{transition maps of toy shtukas are surjective}. Hence $u^{-1}(x)$ is nonempty.

We have
\[\ooToySht_V^n=\oToySht_V^{n,V',V''}-
(\bigcup_{\substack{H\in\mathbf{P}_{V^*}\\H\not\supset V''}}\Delta_{V,H}^{n,V',V''})\cup
(\bigcup_{\substack{J\in\mathbf{P}_{V}\\J\not\subset V'}}\Delta_{V,J}^{n,V',V''}).\]

For $J\in\mathbf{P}_V$ such that $J\subset V''$ and $J\not\subset V'$, we have $u(\Delta_{V,J}^{n,V',V''})\subset\Delta_{V''/V',\bar{J}}^{n-\dim V/V''}$, where $\bar{J}=\im(J\to V''/V')$. So $u^{-1}(x)\cap\Delta_{V,J}^{n,V',V''}$ is empty. For $J\in\mathbf{P}_V$ such that $J\not\subset V''$, \cref{smoothness of restriction of transition map to extra toy horospherical divisors} implies that $u^{-1}(x)\cap\Delta_{V,J}^{n,V',V''}$ has codimension 1 in $u^{-1}(x)$. Thus $u^{-1}(x)\cap(\bigcup_{\substack{J\in\mathbf{P}_{V}\\J\not\subset V'}}\Delta_{V,J}^{n,V',V''})$ has codimension 1 in $u^{-1}(x)$.

Similarly, $u^{-1}(x)\cap(\bigcup_{\substack{H\in\mathbf{P}_{V^*}\\H\not\supset V''}}\Delta_{V,H}^{n,V',V''})$ has codimension 1 in $u^{-1}(x)$.

The statement follows.
\end{proof}

Let $T$ be a nondiscrete noncompact Tate space over $\mathbb{F}_q$ and $n\in\Dim_T$. We denote
\[\ooToySht_T^n=\varprojlim_{(\Lambda',\Lambda'')\in AP_n(T)}\ooToySht_{\Lambda''/\Lambda'}^{n(\Lambda'')}\]

\begin{Lemma}\label{description of ooToySht_T^n}
We have
\[\ooToySht_T^n=\oToySht_T^n-(\bigcup_{H\in\mathbf{P}_{T^*}}\Delta_{T,H}^n)\cup(\bigcup_{J\in\mathbf{P}_T}\Delta_{T,J}^n).\]
\end{Lemma}
\begin{proof}
Let $J\in\mathbf{P}_T$ and $x\in \Delta_{T,J}^n$. Choose $(\Lambda',\Lambda'')\in AP_n(T)$ such that $J\not\subset\Lambda'$, $J\subset\Lambda''$ and  $x\in\oToySht_T^{n,\Lambda',\Lambda''}$. Then the image of $x$ under the morphism $\oToySht_T^{n,\Lambda',\Lambda''}\to\oToySht_{\Lambda''/\Lambda'}^{n(\Lambda'')}$ is contained in $\Delta_{\Lambda''/\Lambda'}^{n(\Lambda'')}$. There is a similar statement for $y\in\Delta_{T,H}^n, H\in\mathbf{P}_{T^*}$.

On the other hand, \cref{open subscheme of Tate toy shtukas as a projective limit} and \cref{set-theoretic containment for the pullback of toy horospherical divisor} imply that for any $(\Lambda',\Lambda'')\in AP_n(T)$ the inverse image of $\Delta_{\Lambda''/\Lambda'}^{n(\Lambda'')}$ under the morphism $\oToySht_T^{n,\Lambda',\Lambda''}\to\oToySht_{\Lambda''/\Lambda'}^{n(\Lambda'')}$ is set-theoretically contained in the union of $\Delta_{T,H}^n$ for $H\in\mathbf{P}_{T^*}$ and $\Delta_{T,J}^n$ for $J\in\mathbf{P}_T$.

The statement follows.
\end{proof}

\begin{Lemma}\label{complement of open chart in the scheme of Tate toy shtukas is horospherical}
For $(\Lambda',\Lambda'')\in AP_n(T)$, $\ooToySht_T^n$ is contained in $\oToySht_T^{n,\Lambda',\Lambda''}$.
\end{Lemma}
\begin{proof}
The proof is similar to that of \cref{complement of open chart in the scheme of toy shtukas is horospherical}.
\end{proof}

\begin{Lemma}\label{surjectivity of transition map for non-horospherical (Tate) toy shtukas}
For $(\Lambda',\Lambda'')\in AP_n(T)$, the morphism $\ooToySht_T^n\to\ooToySht_{\Lambda''/\Lambda'}^{n(\Lambda'')}$ is surjective.
\end{Lemma}
\begin{proof}
The statement follows from \cref{properties of transition map for non-horospherical toy shtukas}.
\end{proof}
\section{Tate toy horospherical subschemes}
Fix a nondiscrete noncompact Tate space $T$ over $\mathbb{F}_q$ and fix $n\in\Dim_T$.

In this section, we use the notation of \cref{(Section)notation for open charts of the scheme of Tate toy shtukas,(Section)admissible pairs of c-lattices}. We also frequently use \cref{open subscheme of Tate toy shtukas as a projective limit} to describe open charts of $\oToySht_T^n$ as a projective limit.

\subsection{Basic properties of Tate toy horospherical subschemes}\label{(Section)basic properties of Tate toy horospherical subschemes}
\begin{Proposition}\label{irreducibility of the moduli scheme of Tate toy shtukas}
  $\oToySht_T^n$ is irreducible.
\end{Proposition}
\begin{proof}
For $(\Lambda',\Lambda''),(\widetilde{\Lambda'},\widetilde{\Lambda''})\in AP_n(T)$ satisfying $(\widetilde{\Lambda'},\widetilde{\Lambda''})\succ(\Lambda',\Lambda'')$, $U_{\widetilde{\Lambda'},\widetilde{\Lambda''}}^{n,\Lambda',\Lambda''}$ is a nonempty open subscheme of $\oToySht_{\widetilde{\Lambda'}/\widetilde{\Lambda''}}^{n(\widetilde{\Lambda''})}$, hence is irreducible by \cref{irreducibility of the scheme of toy shtukas}. So $\oToySht_T^{n,\Lambda',\Lambda''}$ is irreducible as the projective limit. Note that $\oToySht_T^n$ is the union of $\oToySht_T^{n,\Lambda',\Lambda''}$ for all $(\Lambda',\Lambda'')\in AP_n(T)$, and $(\widetilde{\Lambda'},\widetilde{\Lambda''})\succ(\Lambda',\Lambda'')$ implies $\oToySht_T^{n,\Lambda',\Lambda''}\subset\oToySht_T^{n,\widetilde{\Lambda'},\widetilde{\Lambda''}}$. Thus $\oToySht_T^n$ is irreducible.
\end{proof}

For $J\in \mathbf{P}_T$, denote $\Delta_{T,J}^n=\oToySht_T^n\cap\ToySht_{T/J}$. It is called a \emph{Tate toy horospherical subscheme} of $\oToySht_T^n$. Since $T/J$ is a nondiscrete noncompact Tate space over $\mathbb{F}_q$, \cref{irreducibility of the moduli scheme of Tate toy shtukas} implies that $\Delta_{T,J}^n$ is irreducible. Let $\eta_J^n$ be the generic point of $\Delta_{T,J}^n$. For two c-lattices $\Lambda'\subset\Lambda''$ of $T$, denote $\Delta_{T,J}^{n,\Lambda',\Lambda''}=\Delta_{T,J}^n\cap\oGrass_T^{n,\Lambda',\Lambda''}$.

For $H\in \mathbf{P}_{T^*}$, denote $\Delta_{T,H}^n=\oToySht_T^n\cap\ToySht_{H}$. It is called a \emph{Tate toy horospherical subscheme} of $\oToySht_T^n$. Since $H$ is a nondiscrete noncompact Tate space over $\mathbb{F}_q$, \cref{irreducibility of the moduli scheme of Tate toy shtukas} implies that $\Delta_{T,H}^n$ is irreducible. Let $\eta_H^n$ be the generic point of $\Delta_{T,H}^n$. For two c-lattices $\Lambda'\subset\Lambda''$ of $T$, denote $\Delta_{T,H}^{n,\Lambda',\Lambda''}=\Delta_{T,H}^n\cap\oGrass_T^{n,\Lambda',\Lambda''}$.

\begin{Lemma}\label{separation of a d-lattice and a finite dimensional subspace}
  Let $E$ be a field over $\mathbb{F}_q$. Let $P$ be a d-lattice of the Tate space $T\widehat{\otimes}E$ over $E$. Let $W$ be a finite dimensional subspace of $T$ such that $P\cap(W\otimes E)=0$. Then there exists a c-lattice $L$ of $T$ such that $L\supset W$ and $P\cap(L\widehat{\otimes}E)=0$. \qed
\end{Lemma}

\begin{Lemma}\label{Tate toy horospherical divisor as a projective limit}
  Let $J\in\mathbf{P}_T$ and $(\Lambda',\Lambda'')\in AP_n(T)$.
For any $(\widetilde{\Lambda'},\widetilde{\Lambda''})\succ(\Lambda',\Lambda'')$ such that $J\subset\widetilde{\Lambda''}$, we denote $J_{(\widetilde{\Lambda'},\widetilde{\Lambda''})} =\im(J\to\widetilde{\Lambda''}/\widetilde{\Lambda'})\in\mathbf{P}_{\widetilde{\Lambda''}/\widetilde{\Lambda'}}$ and $Z_{(\widetilde{\Lambda'},\widetilde{\Lambda''}),J}
=\Delta_{\widetilde{\Lambda''}/\widetilde{\Lambda'},J_{(\widetilde{\Lambda'},\widetilde{\Lambda''})}}^{n(\widetilde{\Lambda''})}\cap U_{\widetilde{\Lambda'},\widetilde{\Lambda''}}^{n,\Lambda',\Lambda''}
=\Delta_{\widetilde{\Lambda''}/\widetilde{\Lambda'}, J_{(\widetilde{\Lambda'},\widetilde{\Lambda''})}}^{n(\widetilde{\Lambda''}),\Lambda'/\widetilde{\Lambda'},\Lambda''/\widetilde{\Lambda'}}$.
  Then the isomorphism in \cref{open subscheme of Tate toy shtukas as a projective limit} induces an isomorphism
  \[\Delta_{T,J}^{n,\Lambda',\Lambda''}\xrightarrow{\sim}
  \varprojlim\limits_{(\widetilde{\Lambda'},\widetilde{\Lambda''})\succ(\Lambda',\Lambda'')} Z_{(\widetilde{\Lambda'},\widetilde{\Lambda''}),J}.\]
\end{Lemma}
\begin{proof}
  It is clear that the image of $\Delta_{T,J}^{n,\Lambda',\Lambda''}$ is contained in $Z_{(\widetilde{\Lambda'},\widetilde{\Lambda''}),J}$ under each morphism $\oToySht_T^{n,\Lambda',\Lambda''}\to U_{\widetilde{\Lambda'},\widetilde{\Lambda''}}^{n,\Lambda',\Lambda''}$.

  Suppose $x\in\oToySht_T^{n,\Lambda',\Lambda''}(E)$ but $x\notin\Delta_{T,J}^{n,\Lambda',\Lambda''}(E)$, where $E$ is a field over $\mathbb{F}_q$. So $x$ corresponds to a d-lattice $P$ of the Tate space $T\widehat{\otimes}E$ over $E$ such that $P\cap(J\otimes E)=0$. Then \cref{separation of a d-lattice and a finite dimensional subspace} shows that there exists a c-lattice $\Lambda'_0$ of $T$ such that $P\cap(\Lambda'_0\widehat{\otimes}E)=0$. Let $\widetilde{\Lambda'}=\Lambda'\cap\Lambda'_0$ and $\widetilde{\Lambda''}=\Lambda''$. Then the image of $x$ in $U_{\widetilde{\Lambda'},\widetilde{\Lambda''}}^{n,\Lambda',\Lambda''}$ is not contained in $Z_{(\widetilde{\Lambda'},\widetilde{\Lambda''}),J}$. Therefore, the intersection of the preimages of $Z_{(\widetilde{\Lambda'},\widetilde{\Lambda''}),J}$ in $\oToySht_T^{n,\Lambda',\Lambda''}$ for all $(\widetilde{\Lambda'},\widetilde{\Lambda''})\succ(\Lambda',\Lambda'')$ is set-theoretically equal to $\Delta_{T,J}^{n,\Lambda',\Lambda''}$.

  Since $\Delta_{T,J}^{n,\Lambda',\Lambda''}$ and all $Z_{(\widetilde{\Lambda'},\widetilde{\Lambda''}),J}$ are reduced, the statement follows.
\end{proof}

\begin{Lemma}\label{condition for a Tate toy horospherical subscheme to be generically in an open chart}
  For $J\in\mathbf{P}_T$ and $(\Lambda',\Lambda'')\in AP_n(T)$, $\eta_J^n\in\oToySht_T^{n,\Lambda',\Lambda''}$ if and only if $J\not\subset\Lambda'$. For $H\in\mathbf{P}_{T^*}$ and $(\Lambda',\Lambda'')\in AP_n(T)$, $\eta_H^n\in\oToySht_T^{n,\Lambda',\Lambda''}$ if and only if $H\not\supset\Lambda''$.
\end{Lemma}
\begin{proof}
If $J\subset\Lambda'$, then $\ToySht_{T/J}\cap\oToySht_T^{n,\Lambda',\Lambda''}$ is empty, so $\eta_J^n\notin\oToySht_T^{n,\Lambda',\Lambda''}$.

Suppose $J\not\subset\Lambda'$. For any $(\widetilde{\Lambda'},\widetilde{\Lambda''})\succ(\Lambda',\Lambda'')$ such that $J\subset\widetilde{\Lambda''}$, let $J_{(\widetilde{\Lambda'},\widetilde{\Lambda''})}$ and $Z_{(\widetilde{\Lambda'},\widetilde{\Lambda''}),J}$ be as in \cref{Tate toy horospherical divisor as a projective limit}. Since $(\Lambda',\Lambda'')\in AP_n(T)$, we have $n(\widetilde{\Lambda''})\ge 2$, so $Z_{(\widetilde{\Lambda'},\widetilde{\Lambda''}),J}$ is nonempty. \cref{transition maps between toy horospherical divisors are dominant,Tate toy horospherical divisor as a projective limit} imply that $\Delta_{T,J}^{n,\Lambda',\Lambda''}$ is nonempty. Since $\oToySht_T^{n,\Lambda',\Lambda''}$ is open in $\oToySht_T^n$, we deduce that $\eta_J^n\in\oToySht_T^{n,\Lambda',\Lambda''}$.

The proof about $\eta_H^n$ is similar.
\end{proof}

\begin{Lemma}\label{inductive limit of discrete valuation rings}
A filtered inductive limit of discrete valuation rings with common uniformizer is a discrete valuation ring with the same uniformizer.\qed
\end{Lemma}

\begin{Lemma}\label{local ring for a line is discrete}
For $J\in \mathbf{P}_T$, the local ring of $\eta_J^n$ viewed as a point of $\oToySht_T^n$ is a discrete valuation ring.
\end{Lemma}
\begin{proof}
Choose $(\Lambda',\Lambda'')\in AP_n(T)$ such that $J\not\subset\Lambda'$ and $J\subset\Lambda''$. By \cref{condition for a Tate toy horospherical subscheme to be generically in an open chart}, $\eta_J^n$ is in $\oToySht_T^{n,\Lambda',\Lambda''}$.  For any $(\widetilde{\Lambda'},\widetilde{\Lambda''})\succ(\Lambda',\Lambda'')$ such that $J\subset\widetilde{\Lambda''}$, let $J_{(\widetilde{\Lambda'},\widetilde{\Lambda''})}$ and $Z_{(\widetilde{\Lambda'},\widetilde{\Lambda''}),J}$ be as in \cref{Tate toy horospherical divisor as a projective limit}. Since $J\subset\Lambda''$, each $Z_{(\widetilde{\Lambda'},\widetilde{\Lambda''}),J}$ is nonempty, hence irreducible by \cref{irreducibility of the scheme of toy shtukas}. Let $A_{(\widetilde{\Lambda'},\widetilde{\Lambda''}),J}$ denote the local ring of $Z_{(\widetilde{\Lambda'},\widetilde{\Lambda''}),J}$ in $U_{\widetilde{\Lambda'},\widetilde{\Lambda''}}^{n,\Lambda',\Lambda''}$. It is a discrete valuation ring since each $Z_{(\widetilde{\Lambda'},\widetilde{\Lambda''}),J}$ has codimension 1 in $U_{\widetilde{\Lambda'},\widetilde{\Lambda''}}^{n,\Lambda',\Lambda''}$. Let $A_J$ be the local ring of $\eta_J^n$ viewed as a point of $\oToySht_T^n$. We see that $A_J$ is the local ring of $\Delta_{T,J}^{n,\Lambda',\Lambda''}$ in $\oToySht_T^{n,\Lambda',\Lambda''}$. Then the projective limit $\Delta_{T,J}^{n,\Lambda',\Lambda''}=\varprojlim_{(\widetilde{\Lambda'},\widetilde{\Lambda''})\succ(\Lambda',\Lambda'')} Z_{(\widetilde{\Lambda'},\widetilde{\Lambda''}),J}$ gives an isomorphism $A_J=\varinjlim_{(\widetilde{\Lambda'},\widetilde{\Lambda''})\succ(\Lambda',\Lambda'')}A_{(\widetilde{\Lambda'},\widetilde{\Lambda''}),J}$. \cref{pullback of toy horospherical divisors under transition maps} shows that for any $(\Lambda'_2,\Lambda''_2)\succ(\Lambda'_1,\Lambda''_1)\succ(\Lambda',\Lambda'')$, the morphism $A_{(\Lambda'_1,\Lambda''_1),J}\to A_{(\Lambda'_2,\Lambda''_2),J}$ maps the uniformizer to the uniformizer. So the statement follows from \cref{inductive limit of discrete valuation rings}.
\end{proof}

The following statement is dual to \cref{local ring for a line is discrete}.

\begin{Lemma}\label{local ring for a hyperplane is discrete}
For $H\in \mathbf{P}_{T^*}$, the local ring of $\eta_H^n$ viewed as a point of $\oToySht_T^n$ is a discrete valuation ring.
\end{Lemma}

\begin{Remark}
Although $\Delta_{T,J}^n (J\in\mathbf{P}_T)$ and $\Delta_{T,H}^n (J\in\mathbf{P}_{T^*})$ have codimension 1 in $\oToySht_T^n$, they are not Cartier divisors. (See \cref{space of Tate toy horospherical subschemes}.)
\end{Remark}

\subsection{The group of divisors supported on the union of Tate toy horospherical subschemes}
\subsubsection{Formulation of the main result}\label{(Section)formulation of the main result for the group of Tate toy horospherical divisors}
For an open subset $U$ of a Tate space $M$, we denote $C^\infty(U)$ to be the (partially) ordered abelian group of locally constant $\mathbb{Z}$-valued functions on $U$, and denote $C_c^\infty(U)$ to be the (partially) ordered abelian group of locally constant $\mathbb{Z}$-valued functions on $U$ with compact support. We define
\[C_+(U):=\{f\in C^\infty(U)| \text{ $\supp f$ is contained in some compact subset of $M$}\}\]

\begin{Lemma}
We have an isomorphism of ordered abelian groups $C_+(M-\{0\})\cong\varprojlim C_c^\infty(M-\Lambda)$, where the projective limit is taken with respect to the directed set of all c-lattices $\Lambda$ of $M$.
\end{Lemma}
\begin{proof}
\cref{Hausdorff property of a Tate space} and \cref{c-lattices form a basis of neighborhood of 0} imply that $(M-\{0\})$ is the union of $(M-\Lambda)$ for all c-lattices $\Lambda$ of $M$. The statement follows.
\end{proof}

\begin{Definition}
  Let $\mathfrak{O}_T^n$ denote the ordered abelian group of those Cartier divisors of $\oToySht_T^n$ whose restrictions to $\ooToySht_T^n$ are zero. Elements in $\mathfrak{O}_T^n$ are called \emph{(Tate) toy horospherical divisors} of $\oToySht_T^n$.
\end{Definition}

\begin{Remark}
\cref{description of ooToySht_T^n} shows that a Cartier divisor of $\oToySht_T^n$ is an element of $\mathfrak{O}_T^n$ if and only if its support\footnote{The \emph{support} of a divisor is the smallest closed subset such that the
 restriction of the divisor to its complement is zero.} is contained in the union of $\Delta_{T,H}^n$ and $\Delta_{T,J}^n$ for all $H\in\mathbf{P}_{T^*}$ and $J\in\mathbf{P}_{T}$.
\end{Remark}

From \cref{local ring for a line is discrete,local ring for a hyperplane is discrete}, we get an ordered homomorphism $\mathfrak{O}_T^n\to\Maps(\mathbf{P}_{T^*}\coprod\mathbf{P}_T,\mathbb{Z})$, where $\Maps(\mathbf{P}_{T^*}\coprod\mathbf{P}_T,\mathbb{Z})$ is the ordered abelian group of $\mathbb{Z}$-valued functions on $\mathbf{P}_{T^*}\coprod\mathbf{P}_T$.

The goal of this section is to prove the following explicit description of $\mathfrak{O}_T^n$.

\begin{Theorem}\label{space of Tate toy horospherical subschemes}
The homomorphism $\mathfrak{O}_T^n\to \Maps(\mathbf{P}_{T^*}\coprod\mathbf{P}_T,\mathbb{Z})$ induces an isomorphism of ordered abelian groups $\mathfrak{O}_T^n\cong C_+(T^*-\{0\})^{\mathbb{F}_q^\times}\oplus C_+(T-\{0\})^{\mathbb{F}_q^\times}$.
\end{Theorem}

\subsubsection{Reduction of \cref{space of Tate toy horospherical subschemes} to
\cref{space of Tate toy horospherical subschemes on an open chart}}\label{(Section)reduction of proposition about Tate toy horospherical divisors}
For two c-lattices $\Lambda'\subset\Lambda''$ of $T$, let $\mathfrak{O}_T^{n,\Lambda',\Lambda''}$ denote the (partially) ordered abelian group of those Cartier divisors of $\oToySht_T^{n,\Lambda',\Lambda''}$ whose restrictions to $\ooToySht_T^n$ are zero.

From \cref{local ring for a line is discrete,local ring for a hyperplane is discrete,condition for a Tate toy horospherical subscheme to be generically in an open chart}, we get an ordered homomorphism $\mathfrak{O}_T^{n,\Lambda',\Lambda''}\to\Maps((T^*-(\Lambda'')^\perp)\coprod(T-\Lambda'),\mathbb{Z})$.

The isomorphism in \cref{open subscheme of Tate toy shtukas as a projective limit} induces an isomorphism of ordered abelian groups
\[\mathfrak{O}_T^n\xrightarrow{\sim}\varprojlim_{(\Lambda',\Lambda'')\in AP_n(T)}\mathfrak{O}_T^{n,\Lambda',\Lambda''}.\]
So \cref{space of Tate toy horospherical subschemes} follows from \cref{space of Tate toy horospherical subschemes on an open chart}.

\begin{Lemma}\label{space of Tate toy horospherical subschemes on an open chart}
For $(\Lambda',\Lambda'')\in AP_n(T)$, the homomorphism $\mathfrak{O}_T^{n,\Lambda',\Lambda''}\to\Maps((T^*-(\Lambda'')^\perp) \coprod(T-\Lambda'),\mathbb{Z})$ induces an isomorphism of ordered abelian groups $\mathfrak{O}_T^{n,\Lambda',\Lambda''}\cong C_c^\infty(T^*-(\Lambda'')^\perp)^{\mathbb{F}_q^\times}\oplus C_c^\infty(T-\Lambda')^{\mathbb{F}_q^\times}$.
\end{Lemma}

\subsubsection{Reduction of \cref{space of Tate toy horospherical subschemes on an open chart} to
\cref{projective limit of charts induces inductive limit of spaces of toy horospherical divisors},
\cref{inductive limit of spaces of Tate toy horospherical subschemes} and
\cref{space of Tate toy horospherical subschemes from a subquotient}}\label{(Section)inductive limit of spaces of toy horospherical divisors}
Let $\Lambda',\Lambda'',\widetilde{\Lambda'},\widetilde{\Lambda''}$ be c-lattices of $T$ such that $\widetilde{\Lambda'}\subset\Lambda'\subset\Lambda''\subset\widetilde{\Lambda''}$. We have injective homomorphisms of ordered abelian groups
\begin{equation}\label{subgroup containing pullbacks of toy horospherical divisors of type 2}
C^\infty((\widetilde{\Lambda''}/\widetilde{\Lambda'})-(\Lambda'/\widetilde{\Lambda'}))^{\mathbb{F}_q^\times}\to
C_c^\infty(\widetilde{\Lambda''}-\Lambda')^{\mathbb{F}_q^\times}\to
C_c^\infty(T-\Lambda')^{\mathbb{F}_q^\times},
\end{equation}
where the first homomorphism is induced by the map $(\widetilde{\Lambda''}-\Lambda')\to((\widetilde{\Lambda''}/\widetilde{\Lambda'})-(\Lambda'/\widetilde{\Lambda'}))$, the second homomorphism is extension by zero from $(\widetilde{\Lambda''}-\Lambda')$ to $(T-\Lambda')$. We consider $C^\infty((\widetilde{\Lambda''}/\widetilde{\Lambda'})-(\Lambda'/\widetilde{\Lambda'}))^{\mathbb{F}_q^\times}$ as an ordered subgroup of $C_c^\infty(T-\Lambda')^{\mathbb{F}_q^\times}$.

Similarly, we have injective homomorphisms of ordered abelian groups
\begin{equation}\label{subgroup containing pullbacks of toy horospherical divisors of type 1}
C^\infty((\widetilde{\Lambda''}/\widetilde{\Lambda'})^*-(\Lambda''/\widetilde{\Lambda'})^\perp)^{\mathbb{F}_q^\times}\to
C_c^\infty((\widetilde{\Lambda'})^\perp-(\Lambda'')^\perp)^{\mathbb{F}_q^\times}\to
C_c^\infty(T^*-(\Lambda'')^\perp)^{\mathbb{F}_q^\times}.
\end{equation}
We consider $C^\infty((\widetilde{\Lambda''}/\widetilde{\Lambda'})^*-(\Lambda''/\widetilde{\Lambda'})^\perp)^{\mathbb{F}_q^\times}$ as an ordered subgroup of $C_c^\infty(T^*-(\Lambda'')^\perp)^{\mathbb{F}_q^\times}$.

For $(\Lambda',\Lambda'')\in AP_n(T)$ and $(\widetilde{\Lambda'},\widetilde{\Lambda''})\succ(\Lambda',\Lambda'')$, we define $\mathfrak{O}_{\widetilde{\Lambda'},\widetilde{\Lambda''}}^{n,\Lambda',\Lambda''}$ to be the ordered abelian group of those Cartier divisors of $U_{\widetilde{\Lambda'},\widetilde{\Lambda''}}^{n,\Lambda',\Lambda''}$ whose restrictions to $\ooToySht_{\widetilde{\Lambda''}/\widetilde{\Lambda'}}^{n(\widetilde{\Lambda''})}$ are zero.

Now \cref{space of Tate toy horospherical subschemes on an open chart} follows from Lemmas \ref{projective limit of charts induces inductive limit of spaces of toy horospherical divisors}, \ref{inductive limit of spaces of Tate toy horospherical subschemes}, and \ref{space of Tate toy horospherical subschemes from a subquotient}.

\begin{Lemma}\label{projective limit of charts induces inductive limit of spaces of toy horospherical divisors}
For $(\Lambda',\Lambda'')\in AP_n(T)$, the isomorphism
\[\oToySht_T^{n,\Lambda',\Lambda''}\xrightarrow{\sim}
\varprojlim_{(\widetilde{\Lambda'},\widetilde{\Lambda''})\succ(\Lambda',\Lambda'')} U_{\widetilde{\Lambda'},\widetilde{\Lambda''}}^{n,\Lambda',\Lambda''}\]
in \cref{open subscheme of Tate toy shtukas as a projective limit} induces an isomorphism of ordered abelian groups
\[\xi:\varinjlim_{(\widetilde{\Lambda'},\widetilde{\Lambda''})\succ(\Lambda',\Lambda'')}
\mathfrak{O}_{\widetilde{\Lambda'},\widetilde{\Lambda''}}^{n,\Lambda',\Lambda''}
\xrightarrow{\sim}\mathfrak{O}_T^{n,\Lambda',\Lambda''}.\]
\end{Lemma}

\begin{Lemma}\label{inductive limit of spaces of Tate toy horospherical subschemes}
Let $\Lambda',\Lambda''$ be two c-lattices of $T$ such that $\Lambda'\subset\Lambda''$. Then the map $\phi$ induced by (\ref{subgroup containing pullbacks of toy horospherical divisors of type 2}) and the map $\psi$ induced by (\ref{subgroup containing pullbacks of toy horospherical divisors of type 1}) are isomorphisms of ordered abelian groups.
\[\phi:\varinjlim_{(\widetilde{\Lambda'},\widetilde{\Lambda''})\succ(\Lambda',\Lambda'')}
C^\infty((\widetilde{\Lambda''}/\widetilde{\Lambda'})-(\Lambda'/\widetilde{\Lambda'}))^{\mathbb{F}_q^\times}\xrightarrow{\sim}
C_c^\infty(T-\Lambda')^{\mathbb{F}_q^\times}\]
\[\psi:\varinjlim_{(\widetilde{\Lambda'},\widetilde{\Lambda''})\succ(\Lambda',\Lambda'')}
C^\infty((\widetilde{\Lambda''}/\widetilde{\Lambda'})^*-(\Lambda''/\widetilde{\Lambda'})^\perp)^{\mathbb{F}_q^\times}\xrightarrow{\sim}
C_c^\infty(T^*-(\Lambda'')^\perp)^{\mathbb{F}_q^\times}\]
\end{Lemma}

\begin{Lemma}\label{space of Tate toy horospherical subschemes from a subquotient}
For $(\Lambda',\Lambda'')\in AP_n(T)$, the composition
$\mathfrak{O}_{\widetilde{\Lambda'},\widetilde{\Lambda''}}^{n,\Lambda',\Lambda''}\to\mathfrak{O}_T^{n,\Lambda',\Lambda''}
\to\Maps((T^*-(\Lambda'')^\perp)\coprod(T-\Lambda'),\mathbb{Z})$ induces an isomorphism of ordered abelian groups \[\mathfrak{O}_{\widetilde{\Lambda'},\widetilde{\Lambda''}}^{n,\Lambda',\Lambda''}\xrightarrow{\sim}
C^\infty((\widetilde{\Lambda''}/\widetilde{\Lambda'})^*-(\Lambda''/\widetilde{\Lambda'})^\perp)^{\mathbb{F}_q^\times} \oplus
C^\infty((\widetilde{\Lambda''}/\widetilde{\Lambda'})-(\Lambda'/\widetilde{\Lambda'}))^{\mathbb{F}_q^\times}.\]
\end{Lemma}

\subsubsection{Proofs of
\cref{inductive limit of spaces of Tate toy horospherical subschemes}
and \cref{space of Tate toy horospherical subschemes from a subquotient}}
\begin{proof}[Proof of \cref{inductive limit of spaces of Tate toy horospherical subschemes}]
It is obvious that each homomorphism $\phi_{(\widetilde{\Lambda'},\widetilde{\Lambda''})}:
C^\infty((\widetilde{\Lambda''}/\widetilde{\Lambda'})-(\Lambda'/\widetilde{\Lambda'}))^{\mathbb{F}_q^\times}\to
C_c^\infty(T-\Lambda')^{\mathbb{F}_q^\times}$ is injective, and $\phi_{(\widetilde{\Lambda'},\widetilde{\Lambda''})}(f)\le\phi_{(\widetilde{\Lambda'},\widetilde{\Lambda''})}(g)$ if and only if $f\le g$.

Now to prove that $\phi$ is an isomorphism of ordered abelian groups, it suffices to show that $\phi$ is surjective.

Let $f\in C_c^\infty(T-\Lambda')^{\mathbb{F}_q^\times}$.

Since $f$ is compactly supported, its support is contained in some c-lattice $\widetilde{\Lambda''}$.

Since $f$ is locally constant, from \cref{c-lattices form a basis of neighborhood of 0} we see that for every $x\in (T-\Lambda')$ there is a c-lattice $\Lambda_x$ of $T$ such that $f$ is constant on $x+\Lambda_x$. Since $f$ is compactly supported, we can find finitely many $x_1,\dots,x_r\in T$ and c-lattices $\Lambda_1,\dots, \Lambda_r$ of $T$ such that $\supp f\subset\bigcup_{i=1}^r (x_i+\Lambda_i)$ and $f$ is constant on each $x_i+\Lambda_i(1\le i\le r)$. We get a c-lattice $\widetilde{\Lambda'}=\bigcap_{i=1}^r \Lambda_i$.

Thus $f$ is in the image of $C^\infty((\widetilde{\Lambda''}/\widetilde{\Lambda'})-(\Lambda'/\widetilde{\Lambda'}))^{\mathbb{F}_q^\times}$. This shows that $\phi$ is surjective.

The proof for $\psi$ is similar.
\end{proof}

\begin{proof}[Proof of \cref{space of Tate toy horospherical subschemes from a subquotient}]
The statement follows from \cref{equality of multiplicities at toy horospherical divisors for a projective limit}.
\end{proof}

\begin{Lemma}\label{equality of multiplicities at toy horospherical divisors for a projective limit}
Let $(\Lambda',\Lambda'')\in AP_n(T)$ and $Z\in\mathfrak{O}_T^{n,\Lambda',\Lambda''}$. Assume that $Z$ equals the pullback of a Cartier divisor $\widetilde{Z}$ on $U_{\widetilde{\Lambda'},\widetilde{\Lambda''}}^{n,\Lambda',\Lambda''}$ for some $(\widetilde{\Lambda'},\widetilde{\Lambda''})\succ(\Lambda',\Lambda'')$.

Suppose $J\in\mathbf{P}_T$ satisfies $J\not\subset\Lambda'$. If $J\not\subset\widetilde{\Lambda''}$, the multiplicity of $Z$ at $\Delta_{T,J}^{n,\Lambda',\Lambda''}$ is zero. If $J\in\widetilde{\Lambda''}$, the multiplicity of $Z$ at $\Delta_{T,J}^{n,\Lambda',\Lambda''}$ equals the multiplicity of $\widetilde{Z}$ at $\Delta_{\widetilde{\Lambda''}/\widetilde{\Lambda'},\widetilde{J}}^{n(\widetilde{\Lambda''})}\cap U_{\widetilde{\Lambda'},\widetilde{\Lambda''}}^{n,\Lambda',\Lambda''}$, where $\widetilde{J}=\im(J\to \widetilde{\Lambda''}/\widetilde{\Lambda'})$.

Suppose $H\in\mathbf{P}_{T^*}$ satisfies $H\not\supset\Lambda''$. If $H\not\supset\widetilde{\Lambda'}$, then the multiplicity of $Z$ at $\Delta_{T,H}^{n,\Lambda',\Lambda''}$ is zero. If $H\supset\widetilde{\Lambda'}$, then the multiplicity of $Z$ at $\Lambda_{T,H}^{n,\Lambda',\Lambda''}$ equals the multiplicity of $\widetilde{Z}$ at $\Delta_{\widetilde{\Lambda''}/\widetilde{\Lambda'},\widetilde{H}}^{n(\widetilde{\Lambda''})}\cap U_{\widetilde{\Lambda'},\widetilde{\Lambda''}}^{n,\Lambda',\Lambda''}$, where $\widetilde{H}=\im(H\cap\widetilde{\Lambda''}\to \widetilde{\Lambda''}/\widetilde{\Lambda'})$.
\end{Lemma}
\begin{proof}
The statement follows from \cref{pullback of toy horospherical divisors under transition maps}, \cref{Tate toy horospherical divisor as a projective limit} and \cref{inductive limit of discrete valuation rings}.
\end{proof}

\subsubsection{Proof of \cref{projective limit of charts induces inductive limit of spaces of toy horospherical divisors}}\label{(Section)proof of Lemma "projective limit of charts induces inductive limit of spaces of toy horospherical divisors"}
\begin{Lemma}\label{non-mixing property of toy horospherical divisors}
Let $(\widetilde{\Lambda'},\widetilde{\Lambda''})\succ(\Lambda',\Lambda'')$. Suppose $D$ is a Cartier divisor of $U_{\widetilde{\Lambda'},\widetilde{\Lambda''}}^{n,\Lambda',\Lambda''}$ whose pullback to $\oToySht_T^{n,\Lambda',\Lambda''}$ is an element of $\mathfrak{O}_T^{n,\Lambda',\Lambda''}$. Then $D$ is an element of $\mathfrak{O}_{\widetilde{\Lambda'},\widetilde{\Lambda''}}^{n,\Lambda',\Lambda''}$.
\end{Lemma}
\begin{proof}
The statement follows from \cref{description of ooToySht_T^n,surjectivity of transition map for non-horospherical (Tate) toy shtukas}.
\end{proof}

\begin{proof}[Proof of \cref{projective limit of charts induces inductive limit of spaces of toy horospherical divisors}]
By \cref{equality of multiplicities at toy horospherical divisors for a projective limit}, each homomorphism $\xi_{(\widetilde{\Lambda'},\widetilde{\Lambda''})}:
\mathfrak{O}_{\widetilde{\Lambda'},\widetilde{\Lambda''}}^{n,\Lambda',\Lambda''}\to\mathfrak{O}_T^{n,\Lambda',\Lambda''}$ is injective, and $\xi_{(\widetilde{\Lambda'},\widetilde{\Lambda''})}(f)\le\xi_{(\widetilde{\Lambda'},\widetilde{\Lambda''})}(g)$ if and only if $f\le g$.

Now it suffices to show that $\xi$ is surjective. Let $D_0\in\mathfrak{O}_T^{n,\Lambda',\Lambda''}$. By \cref{open subscheme of Tate toy shtukas as a projective limit}, there exists $(\widetilde{\Lambda'},\widetilde{\Lambda''})\succ(\Lambda',\Lambda'')$ and a Cartier divisor $D$ of $U_{\widetilde{\Lambda'},\widetilde{\Lambda''}}^{n,\Lambda',\Lambda''}$ such that $D_0$ equals the pullback of $D$. Now \cref{non-mixing property of toy horospherical divisors} implies that $D\in\mathfrak{O}_{\widetilde{\Lambda'},\widetilde{\Lambda''}}^{n,\Lambda',\Lambda''}$.
\end{proof}

\section{Schubert divisors of $\oToySht_T^n$}
Fix a nondiscrete noncompact Tate space $T$ over $\mathbb{F}_q$ and fix $n\in\Dim_T$.

\subsection{Schubert divisors of Sato Grassmannians}
Let $W$ be a c-lattice of $T$ such that $n(W)=0$. We have a perfect complex
\[\mathscr{S}_W^{\bullet}=(\mathscr{S}_W^{-1}\to\mathscr{S}_W^{0})\]
on $\Grass_T^n$, where $\mathscr{S}_W^{-1}\subset T\widehat{\otimes}\mathscr{O}_{\Grass_T^n}$ is the universal locally free sheaf on $\Grass_T^n$, $\mathscr{S}_W^{0}=(T/W)\otimes\mathscr{O}_{\oGrass_T^n}$, and the map $\mathscr{S}_W^{-1}\to\mathscr{S}_W^{0}$ is the natural one.

Since $n(W)=0$, $\mathscr{S}_W^{-1}\to\mathscr{S}_W^{0}$ is an isomorphism on the big cell $\Grass_T^{n,W,W}$, which is an open dense subscheme of $\Grass_T^n$. So the complex $\mathscr{S}_W^{\bullet}$ is good in the sense of Knudsen-Mumford.

We define the \emph{Schubert divisor} of $\Grass_T^n$ for $W$ to be
\[\Schub_T^W:=\Div(\mathscr{S}_W^{\bullet}).\]

\begin{Lemma}
  $\Schub_T^W\cap\oToySht_T^n$ is a Cartier divisor of $\oToySht_T^n$.
\end{Lemma}
\begin{proof}
  We know that $\oToySht_T^n$ is irreducible and reduced, and $\Schub_T^W$ is a Cartier divisor of $\Grass_T^n$. To prove the statement, it suffices to show that $\oToySht_T^n$ is not contained in $\Schub_T^W$. We have $\oToySht_T^n-\Schub_T^W=\oToySht_T^{n,W,W}$. Choose $(\Lambda',\Lambda'')\in AP_n(T)$ such that $\Lambda'\subset W\subset\Lambda''$. Then $\oToySht_T^{n,\Lambda',\Lambda''}\subset\oToySht_T^{n,W,W}$, and $\oToySht_T^{n,\Lambda',\Lambda''}$ is nonempty by \cref{open subscheme of Tate toy shtukas as a projective limit} and the nonemptiness of each $U_{\widetilde{\Lambda'},\widetilde{\Lambda''}}^{n,\Lambda',\Lambda''}$ in the projective limit there.
\end{proof}

We call $\Schub_T^W\cap\oToySht_T^n$ the \emph{Schubert divisor} of $\oToySht_T^n$ for $W$.

\subsection{Schubert divisors of the scheme of Tate toy shtukas}
\begin{Theorem}\label{description of Schubert divisors of the scheme of Tate toy shtukas}
For a c-lattice $W$ of $T$, the Schubert divisor $\Schub_T^W\cap\oToySht_T^n$ corresponds to the function $(\mathds{1}_{W^\perp-\{0\}},\mathds{1}_{W-\{0\}})$ via \cref{space of Tate toy horospherical subschemes}.
\end{Theorem}
\begin{proof}
Let $\mathscr{L}$ be the universal Tate toy shtuka on $\oToySht_T^n$. Choose $(\Lambda',\Lambda'')\in AP_n(T)$ such that $\Lambda'\subset W\subset \Lambda''$.

On $\oToySht_T^{n,\Lambda',\Lambda''}$, the complex $\mathscr{S}_W^{\bullet}$ is quasi-isomorphic to its subcomplex
\begin{equation}\label{subcomplex of the perfect complex for the Schubert divisor}
\mathscr{L}\cap(\Lambda''\widehat{\otimes}\mathscr{O}_{\oToySht_T^n})\to(\Lambda''/W)\otimes\mathscr{O}_{\oToySht_T^n}.
\end{equation}
From the definition of the projection morphism $u:\oToySht_T^{n,\Lambda',\Lambda''}\to\oToySht_{\Lambda''/\Lambda'}^{n(\Lambda'')}$ in \cref{(Section)notation for open charts of the scheme of Tate toy shtukas}, we see that the complex (\ref{subcomplex of the perfect complex for the Schubert divisor}) on $\oToySht_T^{n,\Lambda',\Lambda''}$ is isomorphic to the pullback of the complex $\mathscr{S}_{\Lambda''/\Lambda',W/\Lambda'}^\bullet$ on $\oToySht_{\Lambda''/\Lambda'}^{n(\Lambda'')}$ by $u$. Then \cref{equivalence of two definitions of the Schubert divisor} implies that
\[\Schub_T^W\cap\oToySht_T^{n,\Lambda',\Lambda''}= u^*(\Schub_{\Lambda''/\Lambda'}^{W/\Lambda'}\cap\oToySht_{\Lambda''/\Lambda'}^{n(\Lambda'')}).\]

Since $\oToySht_T^n$ is covered by the union of $\oToySht_T^{\Lambda',\Lambda''}$ for $(\Lambda',\Lambda'')\in AP_n(T)$ such that $\Lambda'\subset W\subset\Lambda''$, the statement now follows from \cref{equality of multiplicities at toy horospherical divisors for a projective limit} and \cref{description of Schubert divisors of the scheme of toy shtukas}.
\end{proof}

\begin{Proposition}\label{linear equivalence of Schubert divisors for Tate spaces}
For two c-lattices $W_1,W_2$ of $T$ satisfying $n(W_1)=n(W_2)=0$, the two divisors $\Schub_T^{W_1}$ and $\Schub_T^{W_2}$ of $\Grass_T^n$ are linearly equivalent.
\end{Proposition}
\begin{proof}
The statement follows from Theorem 3.3 of~\cite{PM}.
\end{proof}

\begin{Corollary}\label{linear equivalence of Schubert divisors for Tate toy shtukas}
For two c-lattices $W_1,W_2$ of $T$ satisfying $n(W_1)=n(W_2)=0$, the two divisors $\Schub_T^{W_1}\cap\oToySht_T^n$ and $\Schub_T^{W_2}\cap\oToySht_T^n$ of $\oToySht_T^n$ are linearly equivalent.
\end{Corollary}
\begin{proof}
The statement follows from \cref{description of Schubert divisors of the scheme of Tate toy shtukas,linear equivalence of Schubert divisors for Tate spaces}.
\end{proof}

\section{Principal toy horospherical $\mathbb{Z}[\frac{1}{p}]$-divisors of $\oToySht_T^n$}
Fix a nondiscrete noncompact Tate space $T$ over $\mathbb{F}_q$ and fix $n\in\Dim_T$.

We normalize the Haar measure on $T$ by the condition that the measure of any c-lattice $\Lambda$ equals $q^{n(\Lambda)}$.

Recall that the notation $\mathfrak{O}_{\widetilde{\Lambda'},\widetilde{\Lambda''}}^{n,\Lambda',\Lambda''}$ is defined in \cref{(Section)inductive limit of spaces of toy horospherical divisors}. In particular, for $(\Lambda',\Lambda'')\in AP_n(T)$, $\mathfrak{O}_{\Lambda',\Lambda''}^{n,\Lambda',\Lambda''}$ is the ordered abelian group of horospherical divisors of $\oToySht_{\Lambda''/\Lambda'}^{n(\Lambda'')}$.
\subsection{Formulation of the main result}
Recall that $\mathfrak{O}_T^n$ is defined in \cref{(Section)formulation of the main result for the group of Tate toy horospherical divisors}.

Let $\mathfrak{R}_T^n$ be the subgroup of $\mathfrak{O}_T^n$ generated by principal divisors.

Fix a nontrivial additive character $\psi:\mathbb{F}_q\to\mathbb{C}^\times$. We define the Fourier transform
\[\Four_\psi:C_c^\infty(T;\mathbb{C})\to C_c^\infty(T^*;\mathbb{C})\]
such that for any $f\in C_c^\infty(T;\mathbb{C}), \omega\in T^*$, we have
\begin{equation}\label{Fourier transform on a Tate space}
\Four_\psi(f)(\omega)=\int_T f(v)\psi(\omega(v))dv.
\end{equation}
Recall the Haar measure on $T$ is normalized by the condition that the measure of any c-lattice $\Lambda$ equals $q^{n(\Lambda)}$.

When $f\in C_c^\infty(T;\mathbb{Z}[\frac{1}{p}])^{\mathbb{F}_q^\times}$, we have $\Four_\psi(f)\in C_c^\infty(T^*;\mathbb{Z}[\frac{1}{p}])^{\mathbb{F}_q^\times}$, and $\Four_\psi(f)$ does not depend on the choice of $\psi$.

For an open subset $U\subset T$, we define $C_0^\infty(U):=\{f\in C_c^\infty(U)|\int_U f(v)dv=0\}$. Note that the definition does not depend on the choice of the Haar measure on $T$.

\begin{Lemma}
$\im(C_0^\infty(T-\{0\};\mathbb{Z}[\frac{1}{p}])^{\mathbb{F}_q^\times} \xrightarrow{\Four_\psi}C_c^\infty(T^*;\mathbb{Z}[\frac{1}{p}])^{\mathbb{F}_q^\times})\subset C_0^\infty(T^*-\{0\};\mathbb{Z}[\frac{1}{p}])^{\mathbb{F}_q^\times}$. \qed
\end{Lemma}

We make identifications of ordered abelian groups via \cref{space of Tate toy horospherical subschemes}
\[\mathfrak{O}_T^n\otimes\mathbb{Z}[\tfrac{1}{p}]\cong C_+(T^*-\{0\};\mathbb{Z}[\tfrac{1}{p}])^{\mathbb{F}_q^\times}\oplus C_+(T-\{0\};\mathbb{Z}[\tfrac{1}{p}])^{\mathbb{F}_q^\times}.\]

The goal of this section is to prove the following statement.

\begin{Theorem}\label{space of principal rational Tate toy horospherical divisors}
An element $(f_1,f_2)$ of $C_+(T^*-\{0\};\mathbb{Z}[\frac{1}{p}])^{\mathbb{F}_q^\times}\oplus C_+(T-\{0\};\mathbb{Z}[\frac{1}{p}])^{\mathbb{F}_q^\times}\cong\mathfrak{O}_T^n\otimes\mathbb{Z}[\frac{1}{p}]$ is contained in $\mathfrak{R}_T^n\otimes\mathbb{Z}[\frac{1}{p}]$ if and only if $f_2\in C_0^\infty(T-\{0\};\mathbb{Z}[\frac{1}{p}])^{\mathbb{F}_q^\times}$ and $f_1=\Four_\psi(f_2)$.
\end{Theorem}

\subsection{Descent of principal divisors}
Recall that the notation $U_{\widetilde{\Lambda'},\widetilde{\Lambda''}}^{n,\Lambda',\Lambda''}$ is defined in \cref{(Section)open subscheme of Tate toy shtukas as a projective limit}.

\begin{Lemma}\label{descent of principal divisor on some open chart}
Let $h$ be a nonzero rational function on $\oToySht_T^n$ such that $\Div(h)\in\mathfrak{O}_T^n$. Then there exists $(\Lambda',\Lambda'')\in AP_n(T)$ and a rational function $f$ on $\oToySht_{\Lambda''/\Lambda'}^{n(\Lambda'')}$, such that $\Div(f)\in\mathfrak{O}_{\Lambda',\Lambda''}^{n,\Lambda',\Lambda''}$ and the pullback of $f$ to $\oToySht_T^{n,\Lambda',\Lambda''}$ equals $h$.
\end{Lemma}
\begin{proof}
Choose $(\Lambda'_0,\Lambda''_0)\in AP_n(T)$. By \cref{open subscheme of Tate toy shtukas as a projective limit}, there exists $(\Lambda',\Lambda'')\succ(\Lambda'_0,\Lambda''_0)$ and a rational function $f_0$ on $U_{\Lambda',\Lambda''}^{n,\Lambda'_0,\Lambda''_0}$ such that the pullback of $f_0$ to $\oToySht_T^{n,\Lambda'_0,\Lambda''_0}$ equals $h$. Let $f$ be the (unique) rational function on $\oToySht_{\Lambda''/\Lambda'}^{n(\Lambda'')}$ extending $f_0$. The pullback of $f$ to $\oToySht_T^{n,\Lambda',\Lambda''}$ equals $h$ since $\oToySht_T^{n,\Lambda'_0,\Lambda''_0}$ is dense in $\oToySht_T^n$. Now \cref{non-mixing property of toy horospherical divisors} shows that $\Div(f)\in\mathfrak{O}_{\Lambda',\Lambda''}^{n,\Lambda',\Lambda''}$.
\end{proof}

\subsection{Support of extension of pullback of principal toy horospherical divisors}
Recall that in \cref{(Section)notation for functorial properties of toy horospherical divisors}, we defined the notation $\Delta_{V,H}^m,\Delta_{V,J}^m$ for $H\in\mathbf{P}_{V^*},J\in\mathbf{P}_V$. Also we denoted $\Delta_{V,H}^{m,W,V}:=\ToySht_H^m\cap\oGrass_V^{m,W,V}$.

\begin{Lemma}\label{nonemptiness of toy horospherical divisor on an open chart}
Let $V$ be a finite dimensional vector space over $\mathbb{F}_q$ such that $\dim V\ge 3$. Let $W$ be a subspace of $V$. When $1\le m\le \dim V/W-1$, $\Delta_{V,H}^{m,W,V}$ is nonempty for any $H\in\mathbf{P}_V$.
\end{Lemma}
\begin{proof}
Note that $\Delta_{V,H}^{m,W,V}=\oToySht_H^{m,W\cap H,H}$. From the assumptions we know that $\dim H\ge 2$ and $\dim(W\cap H)+m\le\dim H$. The statement follows.
\end{proof}

\begin{Lemma}\label{support of extension of pullback for quotient of principal toy horospherical divisors}
Let $V'\subset V$ be finite dimensional vector spaces over $\mathbb{F}_q$. Assume $2\le m\le\dim V/V'-2$. Let $f$ be a nonzero rational function on $\oToySht_{V/V'}^m$ such that the restriction of $\Div(f)$ to $\ooToySht_{V/V'}^n$ is zero. Let $g_0$ be the pullback of $f$ under the morphism $\oToySht_V^{m,V',V}\to\oToySht_{V/V'}^m$. Let $g$ be the (unique) rational function on $\oToySht_V^m$ extending $g_0$. Then $\Div(g)$ is supported on the union of $\Delta_{V,H}^m$ for $H\in\mathbf{P}_{V^*},H\supset V'$ and $\Delta_{V,J}^m$ for $J\in\mathbf{P}_V,J\not\subset V'$.
\end{Lemma}
\begin{proof}
\cref{complement of open chart in the scheme of toy shtukas is horospherical} implies that $\Div(g)$ is supported on $\Delta_V^m$. So it suffices to show that the multiplicities of $\Div(g)$ are zero at $\Delta_{V,H}^m$ for $H\in\mathbf{P}_{V^*},H\not\supset V'$ and $\Delta_{V,J}^m$ for $J\in\mathbf{P}_V,J\subset V'$.

Let $H\in\mathbf{P}_{V^*}$. Recall that $\Delta_{V,H}^m$ is irreducible by \cref{irreducibility of the scheme of toy shtukas}. We denote $\lambda_H$ to be the multiplicity of $\Div(g)$ at $\Delta_{V,H}^m$, and we denote $\lambda'_{H/V'}$ to be the multiplicity of $\Div(f)$ at $\Delta_{V/V',H/V'}^m$. Then \cref{nonemptiness of toy horospherical divisor on an open chart} and \cref{pullback of toy horospherical divisors for quotient space} show that $\lambda_H=\lambda'_{H/V'}$  if $H\supset V'$, and $\lambda_H=0$ if $H\not\supset V'$.

Let $J\in\mathbf{P}_V$. Denote $\mu_J$ to be the multiplicity of $\Div(g)$ at $\Delta_{V,J}^m$. Applying \cref{space of principal rational toy horospherical divisors} to $V$, we get
\[\mu_J=q^{1-m}\sum_{\substack{H\in\mathbf{P}_{V^*}\\H\supset J}}\lambda_H.\]
Suppose $J\subset V'$. Then the above result about $\lambda_H$ shows that
\[\mu_J=q^{1-m}\sum_{H'\in\mathbf{P}_{(V/V')^*}}\lambda'_{H'}.\]
Applying \cref{space of principal rational toy horospherical divisors} to $V/V'$, we have
\[\sum_{H'\in\mathbf{P}_{(V/V')^*}}\lambda'_{H'}=0.\]
The statement follows.
\end{proof}

The following statement is dual to \cref{support of extension of pullback for quotient of principal toy horospherical divisors}.

\begin{Lemma}\label{support of extension of pullback for subspace of principal toy horospherical divisors}
Let $V''\subset V$ be finite dimensional vector spaces over $\mathbb{F}_q$. Assume $2+\dim V/V''\le m\le \dim V-2$. Let $f$ be a nonzero rational function on $\oToySht_{V''}^{m-\dim V/V''}$ such that the restriction of $\Div(f)$ to $\ooToySht_{V''/V'}^{m-\dim V/V''}$ is zero. Let $g_0$ be the pullback of $f$ under the morphism $\oToySht_V^{m,0,V''}\to\oToySht_{V''}^{m-\dim V/V''}$. Let $g$ be the (unique) rational function on $\oToySht_V^m$ extending $g_0$. Then $\Div(g)$ is supported on the union of $\Delta_{V,H}^m$ for $H\in\mathbf{P}_{V^*},H\not\supset V''$ and $\Delta_{V,J}^m$ for $J\in\mathbf{P}_V,J\subset V''$. \qed
\end{Lemma}

\begin{Lemma}\label{support of extension of pullback for subquotient of principal toy horospherical divisors}
Let $V'\subset V''\subset V$ be finite dimensional vector spaces over $\mathbb{F}_q$. Assume $\dim V/V''+2\le m\le \dim V/V'-2$. Let $f$ be a nonzero rational function on $\oToySht_{V''/V'}^{m-\dim V/V''}$ such that the restriction of $\Div(f)$ to $\ooToySht_{V''/V'}^{m-\dim V/V''}$ is zero. Let $g_0$ be the pullback of $f$ under the morphism $\oToySht_V^{m,V',V''}\to\oToySht_{V''/V'}^{m-\dim V/V''}$. Let $g$ be the (unique) rational function on $\oToySht_V^m$ extending $g_0$. Then $\Div(g)$ is supported on the union of $\Delta_{V,H}^m$ for $H\in\mathbf{P}_{V^*},H\supset V',H\not\supset V''$ and $\Delta_{V,J}^m$ for $J\in\mathbf{P}_V,J\not\subset V',J\subset V''$.
\end{Lemma}
\begin{proof}
The statement follows from \cref{support of extension of pullback for subspace of principal toy horospherical divisors,support of extension of pullback for quotient of principal toy horospherical divisors}.
\end{proof}

\subsection{Extension of pullback of principal toy horospherical divisors}
For $(\Lambda',\Lambda'')\in AP_n(T)$, we define $\varepsilon_{\Lambda',\Lambda''}$ to be the composition of homomorphism of ordered abelian groups
\[\varepsilon_{\Lambda',\Lambda''}:C(\mathbf{P}_{\Lambda''/\Lambda'})\to C^\infty(\Lambda''-\Lambda')\to C_c^\infty(T-\{0\}),\]
where the first homomorphism is induced by the map $(\Lambda''-\Lambda')\to\mathbf{P}_{\Lambda''/\Lambda'}$ and the second homomorphism is extension by zero.

Similarly, we define
\[\varepsilon_{\Lambda',\Lambda''}^*:C(\mathbf{P}_{(\Lambda''/\Lambda')^*})\to C^\infty((\Lambda')^\perp-(\Lambda'')^\perp)\to C_c^\infty(T^*-\{0\})\]

For $(\Lambda',\Lambda'')\in AP_n(T)$, we identify of ordered abelian groups
\[\mathfrak{O}_{\Lambda',\Lambda''}^{n,\Lambda',\Lambda''}\cong C(\mathbf{P}_{(\Lambda''/\Lambda')^*})\oplus C(\mathbf{P}_{\Lambda''/\Lambda'}).\]
As before, we make an identification via \cref{space of principal rational Tate toy horospherical divisors}
\[\mathfrak{O}_T^n\cong C_+(T^*-\{0\})^{\mathbb{F}_q^\times}\oplus C_+(T-\{0\})^{\mathbb{F}_q^\times}.\]
We view $\varepsilon_{\Lambda',\Lambda''}^*\oplus\varepsilon_{\Lambda',\Lambda''}$ as a homomorphism $\mathfrak{O}_{\Lambda',\Lambda''}^{n,\Lambda',\Lambda''}\to\mathfrak{O}_T^n$.

\begin{Lemma}\label{extension of pullback of principal toy horospherical divisors}
Let $(\Lambda',\Lambda'')\in AP_n(T)$ and let $f$ be a rational function on $\oToySht_{\Lambda''/\Lambda'}^{n(\Lambda'')}$ such that $\Div(f)\in\mathfrak{O}_{\Lambda',\Lambda''}^{n,\Lambda',\Lambda''}$. Let $h_0$ be the pullback of $f$ to $\oToySht_T^{n,\Lambda',\Lambda''}$ and let $h$ be the (unique) extension of $h_0$ to $\oToySht_T^n$. Then $\Div(h)\in\mathfrak{O}_T^n$ and $\Div(h)=(\varepsilon_{\Lambda',\Lambda''}^*\oplus\varepsilon_{\Lambda',\Lambda''})(\Div(f))$.
\end{Lemma}
\begin{proof}
\cref{complement of open chart in the scheme of Tate toy shtukas is horospherical} shows that $\Div(h)\in\mathfrak{O}_T^n$.

From  \cref{support of extension of pullback for subquotient of principal toy horospherical divisors} we know that $\Div(h)$ is supported on the union of $\Delta_{T,H}^n$ for $H\in\mathbf{P}_{T^*},H\supset\Lambda',H\not\supset\Lambda''$ and $\Delta_{T,J}^n$ for $J\in\mathbf{P}_T,J\not\subset\Lambda',J\subset\Lambda''$. So $\Div(h)$ is supported on $((\Lambda')^\perp-(\Lambda'')^\perp)\coprod(\Lambda''-\Lambda')$. For $v\in(\Lambda''-\Lambda')$, we have $\Div(h)(v)=\Div(f)(\bar{v})$ by \cref{equality of multiplicities at toy horospherical divisors for a projective limit}, where $\bar{v}$ is the image of $v$ in $\mathbf{P}_{\Lambda''/\Lambda'}$. There is a similar statement for $\omega\in((\Lambda')^\perp-(\Lambda'')^\perp)$. The statement follows.
\end{proof}

\subsection{Proof of \cref{space of principal rational Tate toy horospherical divisors}}
For a finite set $S$, let $C_0(S)$ denote the (partially) ordered abelian group of $\mathbb{Z}$-valued functions on $S$ whose sum over $S$ equals 0.

For $(\Lambda',\Lambda'')\in AP_n(T)$, we define a homomorphism (Radon transform)
\[R_{\Lambda',\Lambda''}^n:C_0(\mathbf{P}_{\Lambda''/\Lambda'};\mathbb{Z}[\tfrac{1}{p}])\to C_0(\mathbf{P}_{(\Lambda''/\Lambda')^*};\mathbb{Z}[\tfrac{1}{p}])\]
such that for any $f\in C_0(\mathbf{P}_{\Lambda''/\Lambda'};\mathbb{Z}[\frac{1}{p}]),H\in\mathbf{P}_{(\Lambda''/\Lambda')^*}$, we have
\[R_{\Lambda',\Lambda''}^n(f)(H)=q^{n(\Lambda')+1}\sum_{\substack{J\in\mathbf{P}_{\Lambda''/\Lambda'}\\J\subset H}}f(J).\]

\begin{Remark}\label{space of principal rational toy horospherical divisors described by Radon transform}
Since we have $n(\Lambda'')-(\dim(\Lambda''/\Lambda')-1)=n(\Lambda')+1$, \cref{space of principal rational toy horospherical divisors} shows that $(f_1,f_2)\in C_0(\mathbf{P}_{(\Lambda''/\Lambda')^*};\mathbb{Z}[\frac{1}{p}])\oplus C_0(\mathbf{P}_{\Lambda''/\Lambda'};\mathbb{Z}[\frac{1}{p}])$ corresponds to a principal toy horospherical $\mathbb{Z}[\frac{1}{p}]$-divisor of $\oToySht_{\Lambda''/\Lambda'}^{n(\Lambda'')}$ if and only if $f_1=R_{\Lambda',\Lambda''}^n(f_2)$.
\end{Remark}

\begin{Lemma}
$\im\left(C_0(\mathbf{P}_{\Lambda''/\Lambda'}) \xrightarrow{\varepsilon_{\Lambda',\Lambda''}}C_c^\infty(T-\{0\})^{\mathbb{F}_q^\times}\right)
\subset C_0^\infty(T-\{0\})^{\mathbb{F}_q^\times}$. \qed
\end{Lemma}

\begin{Lemma}\label{commutativity for Radon transform and Fourier transform}
The following diagram is commutative.
\[\begin{tikzcd}
  C_0(\mathbf{P}_{\Lambda''/\Lambda'};\mathbb{Z}[\frac{1}{p}]) \arrow[r,"R_{\Lambda',\Lambda''}^n"]\arrow[d,"\varepsilon_{\Lambda',\Lambda''}\otimes\mathbb{Z}{[\frac{1}{p}]}"]
  &C_0(\mathbf{P}_{(\Lambda''/\Lambda')^*};\mathbb{Z}[\frac{1}{p}]) \arrow[d,"\varepsilon_{\Lambda',\Lambda''}^*\otimes\mathbb{Z}{[\frac{1}{p}]}"]\\
  C_0^\infty(T-\{0\};\mathbb{Z}[\frac{1}{p}])^{\mathbb{F}_q^\times} \arrow[r,"\Four_\psi"]&C_0^\infty(T^*-\{0\};\mathbb{Z}[\frac{1}{p}])^{\mathbb{F}_q^\times}
\end{tikzcd}\]
\end{Lemma}
\begin{proof}
Fix $f\in C_0(\mathbf{P}_{\Lambda''/\Lambda'};\mathbb{Z}[\frac{1}{p}])$. Let $g=(\varepsilon_{\Lambda',\Lambda''}\otimes\mathbb{Z}[\frac{1}{p}])(f)$, $h_1=(\Four_\psi\bcirc(\varepsilon_{\Lambda',\Lambda''}\otimes\mathbb{Z}[\frac{1}{p}]))(f)$, and $h_2=((\varepsilon_{\Lambda',\Lambda''}^*\otimes\mathbb{Z}[\frac{1}{p}])\bcirc R_{\Lambda',\Lambda''}^n)(f)$.

Let $\omega\in(T^*-\{0\})$.

Since $\psi$ is nontrivial and $g(v)=g(v+x)$ for any $v\in(T-\{0\})$ and $x\in \Lambda'$, we see that $h_1(\omega)=0$ when $\omega\notin(\Lambda')^\perp$. Since $\supp g\subset \Lambda''$ and $\sum_{J\in\mathbf{P}_{\Lambda''/\Lambda'}} f(J)=0$, we see that $h_1(\omega)=0$ when $\omega\in(\Lambda'')^\perp$.

From the definition of $\varepsilon_{\Lambda',\Lambda''}^*$ we see that $h_2(\omega)=0$ when $\omega\notin(\Lambda')^\perp$ or $\omega\in(\Lambda'')^\perp$.

Now we assume $\omega\in ((\Lambda')^\perp-(\Lambda'')^\perp)$.

Denote $H_\omega=\{v\in T|\omega(v)=0\}$ and $\mathbf{P}_{\Lambda''/\Lambda'}^{\omega}= \{J\in\mathbf{P}_{\Lambda''/\Lambda'}|J\subset(H_\omega\cap\Lambda'')/\Lambda'\}$.

From the definitions of $R_{\Lambda',\Lambda''}^n$ and $\varepsilon_{\Lambda',\Lambda''}^*$ we see that
\[h_2(\omega)=q^{n(\Lambda')+1}\sum_{J\in\mathbf{P}_{\Lambda''/\Lambda'}^{\omega}}f(J).\]

On the other hand, we have
\[h_1(\omega)=\int_T g(v)\psi(\omega(v))dv=\mu(\Lambda')\sum_{z\in((\Lambda''/\Lambda')-\{0\})}f(z)\psi(\omega(z)).\]
Since $\psi$ is nontrivial, we get
\[h_1(\omega)=\mu(\Lambda')(\sum_{J\in(\mathbf{P}_{\Lambda''/\Lambda'}-\mathbf{P}_{\Lambda''/\Lambda'}^{\omega})}-f(J)
+\sum_{J\in\mathbf{P}_{\Lambda''/\Lambda'}^{\omega}}(q-1)f(J)).\]
Since $\sum_{J\in\mathbf{P}_{\Lambda''/\Lambda'}} f(J)=0$, we get
\[h_1(\omega)=q^{n(\Lambda')}\cdot q\cdot\sum_{J\in\mathbf{P}_{\Lambda''/\Lambda'}^{\omega}}f(J).\]
(Recall the Haar measure on $T$ is normalized by the condition that the measure of any c-lattice $\Lambda$ equals $q^{n(\Lambda)}$.) Therefore, we have $h_1=h_2$.
\end{proof}

\begin{Lemma}\label{C_0 as the union of functions on finite projective spaces}
We have
\[C_0^\infty(T-\{0\};\mathbb{Z}[\tfrac{1}{p}])^{\mathbb{F}_q^\times}=
\bigcup_{(\Lambda',\Lambda'')\in AP_n(T)} \im(C_0(\mathbf{P}_{\Lambda''/\Lambda'};\mathbb{Z}[\tfrac{1}{p}])
\xrightarrow{\varepsilon_{\Lambda',\Lambda''}\otimes\mathbb{Z}[\tfrac{1}{p}]} C_c^\infty(T-\{0\};\mathbb{Z}[\tfrac{1}{p}])^{\mathbb{F}_q^\times}). \qed\]
\end{Lemma}

\begin{proof}[Proof of \cref{space of principal rational Tate toy horospherical divisors}]
The statement follows from \cref{descent of principal divisor on some open chart}, \cref{extension of pullback of principal toy horospherical divisors}, \cref{space of principal rational toy horospherical divisors described by Radon transform}, \cref{commutativity for Radon transform and Fourier transform}, and \cref{C_0 as the union of functions on finite projective spaces}.
\end{proof}

\section{Partial Frobeniuses for Tate toy shtukas}
Fix a nondiscrete noncompact Tate space $T$ over $\mathbb{F}_q$. Let $\Dim_T$ denote the dimension torsor of $T$.
\subsection{Definition of left/right Tate toy shtukas}
\begin{Definition}
A \emph{right Tate toy shtuka} for $T$ over an $\mathbb{F}_q$-scheme $S$ of dimension $n\in\Dim_T$ is a pair $\mathscr{L}\in\Grass_T^n(S), \mathscr{L}'\in\Grass_T^{n+1}(S)$, such that $\mathscr{L}\subset\mathscr{L}', \Fr_S^*\mathscr{L}\subset\mathscr{L}'$ with $\mathscr{L}'/\mathscr{L}$ and $\mathscr{L}'/\Fr_S^*\mathscr{L}$ being invertible.
\end{Definition}

\begin{Definition}
A \emph{left Tate toy shtuka} for $T$ over an $\mathbb{F}_q$-scheme $S$ of dimension $n\in\Dim_T$ is a pair $\mathscr{L}\in\Grass_T^n(S), \mathscr{L}'\in\Grass_T^{n-1}(S)$, such that $\mathscr{L}'\subset\mathscr{L}, \mathscr{L}'\subset\Fr_S^*\mathscr{L}$ with $\mathscr{L}/\mathscr{L}'$ and $\Fr_S^*\mathscr{L}/\mathscr{L}'$ being invertible.
\end{Definition}

For $n\in\Dim_T$, let $\LToySht_T^n$ (resp. $\RToySht_T^n$) be the functor which associates to each $\mathbb{F}_q$-scheme $S$ the set of isomorphism classes of left (resp. right) Tate toy shtukas over $S$ of dimension $n$.

As in the finite dimensional case, $\LToySht_T^n$ (resp. $\RToySht_T^n$) is representable by a closed subscheme of $\Grass_T^n\times\Grass_T^{n-1}$ (resp. $\Grass_T^n\times\Grass_T^{n+1}$).

\subsection{Partial Frobeniuses}
We have the following constructions for left/right Tate toy shtukas:

(i) For a left Tate toy shtuka $\mathscr{L}'\subset\mathscr{L}$ over an $\mathbb{F}_q$-scheme $S$, the pair $\mathscr{L}'\subset\Fr_S^*\mathscr{L}$ forms a right Tate toy shtuka over $S$.

(ii) For a right Tate toy shtuka $\mathscr{L}\subset\mathscr{L}'$ over an $\mathbb{F}_q$-scheme $S$, the pair $\Fr_S^*\mathscr{L}\subset\mathscr{L}'$ forms a left Tate toy shtuka over $S$.

For $n\in\Dim_T$, we define partial Frobeniuses $F_{T,n}^-:\LToySht_T^n\to\RToySht_T^{n-1}, F_{T,n}^+:\RToySht_T^n\to\LToySht_T^{n+1}$ induced by the above constructions.

We have $F_{T,n-1}^+\bcirc F_{T,n}^-=\Fr_{\LToySht_T^n}, F_{T,n+1}^-\bcirc F_{T,n}^+=\Fr_{\RToySht_T^n}$.

\subsection{Identification of $L\mathfrak{O}_T^n$ and $R\mathfrak{O}_T^n$ with $\mathfrak{O}_T^n$}
As in the finite dimensional case, $\oToySht_T^n$ is an open subscheme of $\LToySht_T^n$ and $\RToySht_T^n$.

Let $L\mathfrak{O}_T^n$ (resp. $R\mathfrak{O}_T^n$) be the (partially) ordered abelian group of Cartier divisors of $\LToySht_T^n$ (resp. $\RToySht_T^n$) whose restrictions to $\ooToySht_T^n$ are zero.

The goal of this subsection is to prove the following statement.

\begin{Lemma}\label{identification of left/right Tate toy horospherical divisors}
The open immersion $\oToySht_T^n\hookrightarrow\LToySht_T^n$ (resp. $\oToySht_T^n\hookrightarrow\RToySht_T^n$) induces an isomorphism of ordered abelian groups $L\mathfrak{O}_T^n\xrightarrow{\sim}\mathfrak{O}_T^n$ (resp. $R\mathfrak{O}_T^n\xrightarrow{\sim}\mathfrak{O}_T^n$).
\end{Lemma}

In the following we only prove the statement for $\LToySht_T^n$, since the statements for $\LToySht_T^n$ and $\RToySht_T^n$ are dual to each other.

\begin{Lemma}\label{open chart for oToySht contained in open chart for LToySht}
Let $E$ be a field over $\mathbb{F}_q$ and let $L\in\oToySht_T^{n,\Lambda',\Lambda''}(E)$ for two c-lattices $\Lambda'\subset\Lambda''$ of $T$. Put $L'=L\cap\Fr_E^*L$. Then $L'\in\Grass_T^{n-1,\Lambda',\Lambda''}(E)$.
\end{Lemma}
\begin{proof}
For any subspace $M$ of $T$, we denote $M_E:=M\widehat{\otimes}E$.

Since $L\in\oGrass_T^{n,\Lambda',\Lambda''}(E)$, we have $L\cap\Lambda'_E=0$. Hence $L'\cap\Lambda'_E=0$. So it suffices to show that $L'+\Lambda''_E=T_E$.

We have $L+\Lambda''_E=T_E$. Hence $\dim(L\cap\Lambda''_E)=\dim(L)(\Lambda''_E)=n(\Lambda'')$. (See \cref{definition of the dimension of a d-lattice} for the definition of the dimension of a d-lattice.) From the definition of $\oGrass$ in diagram (\ref{Cartesian diagram for oGrass}) we know that $L\cap\Lambda''_E$ is a nontrivial toy shtuka over $\Spec E$. Thus $\dim(L'\cap\Lambda''_E)=\dim(L\cap\Lambda''_E)-1=n(\Lambda'')-1$. Therefore, in the Tate space $T_E$ over $E$, the d-lattice $L'$ and the c-lattice $\Lambda''_E$ satisfy $\dim(L')(\Lambda''_E)=\dim(L'\cap\Lambda''_E)$. This implies that $L'+\Lambda''_E=T_E$.
\end{proof}

For $(\Lambda',\Lambda'')\in AP_n(T)$, denote $\LToySht_T^{n,\Lambda',\Lambda''}=\LToySht_T^n\cap(\Grass_T^{n,\Lambda',\Lambda''}\times\Grass_T^{n-1,\Lambda',\Lambda''})$. \cref{open chart for oToySht contained in open chart for LToySht} shows that $\oToySht_T^{n,\Lambda',\Lambda''}\subset\LToySht_T^{n,\Lambda',\Lambda''}$.

For $(\Lambda',\Lambda'')\in AP_n(T)$, let $L\mathfrak{O}_T^{n,\Lambda',\Lambda''}$ be the (partially) ordered abelian group of Cartier divisors of $\LToySht_T^{n,\Lambda',\Lambda''}$ whose restrictions to $\ooToySht_T^n$ are zero.

The isomorphisms
\[\oToySht_T^n=\varinjlim\limits_{(\Lambda',\Lambda'')\in AP_n(T)}\oToySht_T^{n,\Lambda',\Lambda''},\]
\[\LToySht_T^n=\varinjlim\limits_{(\Lambda',\Lambda'')\in AP_n(T)}\LToySht_T^{n,\Lambda',\Lambda''}\]
give isomorphisms of ordered abelian groups
\[\mathfrak{O}_T^n\xrightarrow{\sim}\varprojlim_{(\Lambda',\Lambda'')\in AP_n(T)}\mathfrak{O}_T^{n,\Lambda',\Lambda''},\]
\[L\mathfrak{O}_T^n\xrightarrow{\sim}\varprojlim_{(\Lambda',\Lambda'')\in AP_n(T)}L\mathfrak{O}_T^{n,\Lambda',\Lambda''}.\]

For different $(\Lambda',\Lambda'')\in AP_n(T)$, the open immersions $\oToySht_T^{n,\Lambda',\Lambda''}\hookrightarrow\LToySht_T^{n,\Lambda',\Lambda''}$ induce ordered homomorphisms $L\mathfrak{O}_T^{n,\Lambda',\Lambda''}\to\mathfrak{O}_T^{n,\Lambda',\Lambda''}$ which are compatible with transition homomorphisms. Thus \cref{identification of left/right Tate toy horospherical divisors} follows from \cref{identification of left Tate toy horospherical divisors on open charts}.

\begin{Lemma}\label{identification of left Tate toy horospherical divisors on open charts}
For $(\Lambda',\Lambda'')\in AP_n(T)$, the open immersion $\oToySht_T^{n,\Lambda',\Lambda''}\hookrightarrow\LToySht_T^{n,\Lambda',\Lambda''}$ induces an isomorphism of ordered abelian groups $L\mathfrak{O}_T^{n,\Lambda',\Lambda''}\xrightarrow{\sim}\mathfrak{O}_T^{n,\Lambda',\Lambda''}$.
\end{Lemma}
\begin{proof}
For $(\widetilde{\Lambda'},\widetilde{\Lambda''})\succ(\Lambda',\Lambda'')$, denote $LU_{\widetilde{\Lambda'},\widetilde{\Lambda''}}^{n,\Lambda',\Lambda''}
=\LToySht_{\widetilde{\Lambda''}/\widetilde{\Lambda'}}^{n(\widetilde{\Lambda''}), \Lambda'/\widetilde{\Lambda'},\Lambda''/\widetilde{\Lambda'}}$, and denote $L\mathfrak{O}_{\widetilde{\Lambda'},\widetilde{\Lambda''}}^{n,\Lambda',\Lambda''}$ to be the ordered abelian group of Cartier divisors of $LU_{\widetilde{\Lambda'},\widetilde{\Lambda''}}^{n,\Lambda',\Lambda''}$ whose restrictions to $\ooToySht_{\widetilde{\Lambda''}/\widetilde{\Lambda'}}^{n(\widetilde{\Lambda''}), \Lambda'/\widetilde{\Lambda'},\Lambda''/\widetilde{\Lambda'}}$ are zero.

The complement of $\oToySht_{\widetilde{\Lambda''}/\widetilde{\Lambda'}}^{n(\widetilde{\Lambda''})}$ in $\LToySht_{\widetilde{\Lambda''}/\widetilde{\Lambda'}}^{n(\widetilde{\Lambda''})}$ has codimension $-n(\widetilde{\Lambda'})\ge 2$ since $(\widetilde{\Lambda'},\widetilde{\Lambda''})\in AP_n(T)$. Hence the complement of $U_{\widetilde{\Lambda'},\widetilde{\Lambda''}}^{n,\Lambda',\Lambda''}$ in $LU_{\widetilde{\Lambda'},\widetilde{\Lambda''}}^{n,\Lambda',\Lambda''}$ has codimension at least 2. Therefore, the open immersion of regular schemes $U_{\widetilde{\Lambda'},\widetilde{\Lambda''}}^{n,\Lambda',\Lambda''}\hookrightarrow LU_{\widetilde{\Lambda'},\widetilde{\Lambda''}}^{n,\Lambda',\Lambda''}$ induces an isomorphism of ordered abelian groups $L\mathfrak{O}_{\widetilde{\Lambda'},\widetilde{\Lambda''}}^{n,\Lambda',\Lambda''}\xrightarrow{\sim} \mathfrak{O}_{\widetilde{\Lambda'},\widetilde{\Lambda''}}^{n,\Lambda',\Lambda''}$. Such isomorphisms are compatible with transition homomorphisms for different pairs $(\widetilde{\Lambda'},\widetilde{\Lambda''})$, and thus give rise to an isomorphism of ordered abelian groups
\[\lambda:
\varinjlim\limits_{(\widetilde{\Lambda'},\widetilde{\Lambda''})\succ(\Lambda',\Lambda'')}
L\mathfrak{O}_{\widetilde{\Lambda'},\widetilde{\Lambda''}}^{n,\Lambda',\Lambda''}
\xrightarrow{\sim}
\varinjlim\limits_{(\widetilde{\Lambda'},\widetilde{\Lambda''})\succ(\Lambda',\Lambda'')}
\mathfrak{O}_{\widetilde{\Lambda'},\widetilde{\Lambda''}}^{n,\Lambda',\Lambda''}.\]

Similarly to \cref{open subscheme of Tate toy shtukas as a projective limit}, we have an isomorphism
\[\LToySht_T^{n,\Lambda',\Lambda''}\xrightarrow{\sim}
\varprojlim_{(\widetilde{\Lambda'},\widetilde{\Lambda''})\succ(\Lambda',\Lambda'')} LU_{\widetilde{\Lambda'},\widetilde{\Lambda''}}^{n,\Lambda',\Lambda''}.\]
It induces an ordered homomorphism
\[L\xi:
\varinjlim\limits_{(\widetilde{\Lambda'},\widetilde{\Lambda''})\succ(\Lambda',\Lambda'')}
L\mathfrak{O}_{\widetilde{\Lambda'},\widetilde{\Lambda''}}^{n,\Lambda',\Lambda''} \to
L\mathfrak{O}_T^{n,\Lambda',\Lambda''}.\]
Similarly to the proof of \cref{projective limit of charts induces inductive limit of spaces of toy horospherical divisors} in \cref{(Section)proof of Lemma "projective limit of charts induces inductive limit of spaces of toy horospherical divisors"}, one can show that $L\xi$ is surjective using \cref{surjectivity of transition map for non-horospherical (Tate) toy shtukas}.

We get a commutative diagram of ordered abelian groups.
\[\begin{tikzcd}
\varinjlim\limits_{(\widetilde{\Lambda'},\widetilde{\Lambda''})\succ(\Lambda',\Lambda'')}
L\mathfrak{O}_{\widetilde{\Lambda'},\widetilde{\Lambda''}}^{n,\Lambda',\Lambda''} \arrow[r,"L\xi"]\arrow[d,"\lambda"]
&L\mathfrak{O}_T^{n,\Lambda',\Lambda''} \arrow[d,"\lambda'"]\\
\varinjlim\limits_{(\widetilde{\Lambda'},\widetilde{\Lambda''})\succ(\Lambda',\Lambda'')}
\mathfrak{O}_{\widetilde{\Lambda'},\widetilde{\Lambda''}}^{n,\Lambda',\Lambda''} \arrow[r,"\xi"]
&\mathfrak{O}_T^{n,\Lambda',\Lambda''}
\end{tikzcd}\]
We proved that $\lambda$ is an isomorphism of ordered abelian groups. So is $\xi$ by \cref{projective limit of charts induces inductive limit of spaces of toy horospherical divisors}. Both $L\xi$ and $\lambda'$ are ordered homomorphisms, and we proved that $L\xi$ is surjective. Hence $\lambda'$ is an isomorphism of ordered abelian groups.
\end{proof}

\subsection{Pullback of Tate toy horospherical subschemes under partial Frobeniuses}
We identify $L\mathfrak{O}_T^n$ and $R\mathfrak{O}_T^n$ with $\mathfrak{O}_T^n$ by \cref{identification of left/right Tate toy horospherical divisors}.

Partial Frobeniuses induce homomorphisms between ordered abelian groups
\[\begin{tikzcd}
\dots\arrow[r,shift left=0.5ex]
&\mathfrak{O}_T^{n-1}\arrow[r,shift left=0.5ex,"(F_{T,n}^-)^*"]\arrow[l,shift left=0.5ex]
&\mathfrak{O}_T^n\arrow[r,shift left=0.5ex,"(F_{T,n+1}^-)^*"]\arrow[l,shift left=0.5ex,"(F_{T,n-1}^+)^*"]
&\mathfrak{O}_T^{n+1}\arrow[l,shift left=0.5ex,"(F_{T,n}^+)^*"]\arrow[r,shift left=0.5ex]
&\dots\arrow[l,shift left=0.5ex]
\end{tikzcd}\]

Similarly to \cref{multiplicity 1 of pullback of toy horospherical divisors by partial Frobeniuses,multiplicity q of pullback of toy horospherical divisors by partial Frobeniuses}, we have the following statement.
\begin{Lemma}\label{pullback of Tate toy horospherical divisors under partial Frobeniuses}
Identifying $\mathfrak{O}_T^{n-1},\mathfrak{O}_T^n,\mathfrak{O}_T^{n-1}$ with $C_+(T^*-\{0\})^{\mathbb{F}_q^\times}\oplus C_+(T-\{0\})^{\mathbb{F}_q^\times}$ via \cref{space of Tate toy horospherical subschemes}, the two homomorphisms
\begin{align*}
(F_{T,n}^-)^*:\mathfrak{O}_T^{n-1}&\to \mathfrak{O}_T^n,\\
(F_{T,n}^+)^*:\mathfrak{O}_T^{n+1}&\to \mathfrak{O}_T^n
\end{align*}
are given by
\begin{align*}
(F_{T,n}^-)^*: C_+(T^*-\{0\})^{\mathbb{F}_q^\times}\oplus C_+(T-\{0\})^{\mathbb{F}_q^\times}&\to C_+(T^*-\{0\})^{\mathbb{F}_q^\times}\oplus C_+(T-\{0\})^{\mathbb{F}_q^\times}\\
(\lambda_1,\lambda_2)&\mapsto(q\lambda_1,\lambda_2)\\
(F_{T,n}^+)^*: C_+(T^*-\{0\})^{\mathbb{F}_q^\times}\oplus C_+(T-\{0\})^{\mathbb{F}_q^\times}&\to C_+(T^*-\{0\})^{\mathbb{F}_q^\times}\oplus C_+(T-\{0\})^{\mathbb{F}_q^\times}\\
(\lambda_1,\lambda_2)&\mapsto(\lambda_1,q\lambda_2)
\end{align*}
\end{Lemma}

\section{A canonical subgroup of $\Pic(\oToySht_T^n)$}
Fix a nondiscrete noncompact Tate space $T$ over $\mathbb{F}_q$ and fix $n\in\Dim_T$.

\subsection{Definition}
Let $\mathscr{L}_{T,n}$ be the universal Tate toy shtuka over $\oToySht_T^n$. Let $\widetilde{\mathscr{L}}_{T,n}^L\subset\mathscr{L}_{T,n}^L$ (resp. $\mathscr{L}_{T,n}^R\subset\widetilde{\mathscr{L}}_{T,n}^R$) be the universal left (resp. right) Tate toy shtuka over $\LToySht_T^n$ (resp. $\RToySht_T^n$).

Denote $\mathscr{L}'_{T,n}=\mathscr{L}_{T,n}\cap\Fr_{\oToySht_T^n}^*\mathscr{L}_{T,n}, \mathscr{L}''_{T,n}=\mathscr{L}_{T,n}+\Fr_{\oToySht_T^n}^*\mathscr{L}_{T,n}$.

\begin{Remark}\label{intersection and sum with Frobenius pullback as pullback under partial Frobeniuses}
We see that $\mathscr{L}'_{T,n}$ is the pullback of $\mathscr{L}_{T,n-1}^R$ under the composition $\oToySht_T^n\hookrightarrow\LToySht_T^n\xrightarrow{F_{T,n}^-}\RToySht_T^{n-1}$, and $\mathscr{L}''_{T,n}$ is the pullback of $\mathscr{L}_{T,n+1}^L$ under the composition $\oToySht_T^n\hookrightarrow\RToySht_T^n\xrightarrow{F_{T,n}^+}\LToySht_T^{n+1}$.
\end{Remark}

\begin{Definition}
We define two invertible sheaves $\ell_{T,n,a}:=\mathscr{L}_{T,n}/\mathscr{L}'_{T,n}, \ell_{T,n,b}:=\mathscr{L}''_{T,n}/\mathscr{L}_{T,n}$ on $\oToySht_T^n$.
\end{Definition}

\begin{Definition}
For a c-lattice $W$ of $T$, we define an invertible sheaf $\ell_{T,n,\det}^W:=\det(\mathscr{L}_{T,n},W)$. (See \cref{definition of the determinant of a family of d-lattices relative to a c-lattice} for the definition of $\det(\mathscr{L}_{T,n},W)$.)
\end{Definition}

\begin{Remark}\label{isomorphism between determinant line bundles for different c-lattices}
For two c-lattices $W_1,W_2$ of $T$, \cref{change of c-lattices for relative determinant of a family of d-lattices} shows that there is a canonical isomorphism $\ell_{T,n,\det}^{W_1}\otimes\det_{W_1}^{W_2}\cong\ell_{T,n,\det}^{W_2}$, and the two invertible sheaves $\ell_{T,n,\det}^{W_1}$ and $\ell_{T,n,\det}^{W_2}$ are isomorphic. Also, $(\ell_{T,n,\det}^{W_1})^{\otimes(q-1)}$ and $(\ell_{T,n,\det}^{W_2})^{\otimes(q-1)}$ are canonically isomorphic since $(\det_{W_1}^{W_2})^{\otimes(q-1)}$ is canonically isomorphic to $\mathbb{F}_q$.
\end{Remark}

\begin{Lemma}\label{relation between 3 line bundles on the scheme of Tate toy shtukas}
Let $W$ be a c-lattice of $T$. We have a canonical isomorphism $\ell_{T,n,b}\otimes\ell_{T,n,a}^{-1}\cong(\ell_{T,n,\det}^W)^{\otimes(q-1)}$.
\end{Lemma}
\begin{proof}
For any $\mathscr{O}_{\oToySht_T^n}$-module $\mathscr{E}$, denote $\tensor[^\tau]{\mathscr{E}}{}=\Fr_{\oToySht_T^n}^*\mathscr{E}$.

The canonical isomorphism $\mathscr{L}_{T,n}/\mathscr{L}'_{T,n}\cong\mathscr{L}''_{T,n}/\tensor[^\tau]{\mathscr{L}}{_{T,n}}$ induces a canonical isomorphism $\det(\mathscr{L}_{T,n},W)\otimes\det(\mathscr{L}'_{T,n},W)^{-1}\cong \det(\mathscr{L}''_{T,n},W)\otimes\det(\tensor[^\tau]{\mathscr{L}}{_{T,n}},W)^{-1}$ by \cref{canonical isomorphism of relative determinants for short exact sequence}.

\cref{canonical isomorphism of relative determinants for Frobenius pullback} gives a canonical isomorphism $\det(\tensor[^\tau]{\mathscr{L}}{_{T,n}},W)\cong \det(\mathscr{L}_{T,n},W)^{\otimes q}$.

\cref{canonical isomorphism of relative determinants for short exact sequence} gives canonical isomorphisms
\[\ell_{T,n,a}\cong \det(\mathscr{L}_{T,n},W)\otimes\det(\mathscr{L}'_{T,n},W)^{-1},\]
\[\ell_{T,n,b}\cong \det(\mathscr{L}''_{T,n},W)\otimes\det(\mathscr{L}_{T,n},W)^{-1}.\]

The statement follows.
\end{proof}

\subsection{Some divisors in the classes of $\ell_{T,n,a}$ and $\ell_{T,n,b}$}
We have an inclusion of ordered abelian groups
\[C_c^\infty(T^*)^{\mathbb{F}_q^\times}\oplus C_c^\infty(T)^{\mathbb{F}_q^\times}\subset
C_+(T^*-\{0\})^{\mathbb{F}_q^\times}\oplus C_+(T-\{0\})^{\mathbb{F}_q^\times}.\]

For an invertible sheaf $\ell$ on $\oToySht_T^n$ and an element $(f_1,f_2)\in C_c^\infty(T^*)^{\mathbb{F}_q^\times}\oplus C_c^\infty(T)^{\mathbb{F}_q^\times}$, we write $\ell\sim(f_1,f_2)$ if $[\ell]=[\mathscr{O}_{\oToySht_T^n}(D_{(f_1,f_2)})]$ in $\Pic(\oToySht_T^n)$, where $D_{(f_1,f_2)}$ is the Cartier divisor of $\oToySht_T^n$ corresponding to $(f_1,f_2)$ via \cref{space of Tate toy horospherical subschemes}.

\begin{Lemma}\label{divisor of 2 line bundles on the scheme of Tate toy shtukas}
Let $W_{-1},W_0,W_1$ be c-lattices of $T$ such that $W_{-1}\subset W_0\subset W_1$ and $n(W_i)=i, (i=-1,0,1)$. Then we have
\begin{align*}
\ell_{T,n,a}&\sim(q\cdot\mathds{1}_{W_1^\perp}-\mathds{1}_{W_0^\perp},\mathds{1}_{W_1- W_0}),\\
\ell_{T,n,b}&\sim(-\mathds{1}_{W_{-1}^\perp- W_0^\perp},-q\cdot\mathds{1}_{W_{-1}}+\mathds{1}_{W_0}).
\end{align*}
\end{Lemma}
\begin{proof}
Recall that we denoted $\mathscr{L}'_{T,n}=\mathscr{L}_{T,n}\cap\Fr_{\oToySht_T^n}^*\mathscr{L}_{T,n}, \mathscr{L}''_{T,n}=\mathscr{L}_{T,n}+\Fr_{\oToySht_T^n}^*\mathscr{L}_{T,n}$.

Recall that $\ell_{T,n,a}=\mathscr{L}_{T,n}/\mathscr{L}'_{T,n}$. So \cref{canonical isomorphism of relative determinants for short exact sequence} gives an isomorphism $\ell_{T,n,a}\cong\det(\mathscr{L}_{T,n},W_0)\otimes\det(\mathscr{L}'_{T,n},W_0)^{-1}$. Hence $\ell_{T,n,a}\cong\det(\mathscr{L}_{T,n},W_0)\otimes\det(\mathscr{L}'_{T,n},W_1)^{-1}$ by \cref{change of c-lattices for relative determinant of a family of d-lattices}.

\cref{description of Schubert divisors of the scheme of Tate toy shtukas} implies that $\det(\mathscr{L}_{T,n},W_0)\sim(-\mathds{1}_{W_0^\perp},-\mathds{1}_{W_0})$.

\cref{description of Schubert divisors of the scheme of Tate toy shtukas} and \cref{identification of left/right Tate toy horospherical divisors} imply that $\det(\mathscr{L}_{T,n-1}^R,W_1)\sim(-\mathds{1}_{W_1^\perp},-\mathds{1}_{W_1})$. \cref{intersection and sum with Frobenius pullback as pullback under partial Frobeniuses} and \cref{base change for relative determinant of a family of d-lattices} show that $\det(\mathscr{L}'_{T,n},W_1)$ is the pullback of $\det(\mathscr{L}_{T,n-1}^R,W_1)$ under the composition $\oToySht_T^n\hookrightarrow\LToySht_T^n\xrightarrow{F_{T,n}^-}\RToySht_T^{n-1}$. Then \cref{pullback of Tate toy horospherical divisors under partial Frobeniuses} implies that $\det(\mathscr{L}'_{T,n},W_1)\sim(-q\cdot\mathds{1}_{W_1^\perp},-\mathds{1}_{W_1})$.

Therefore, we get
\[\ell_{T,n,a}\sim (-\mathds{1}_{W_0^\perp},-\mathds{1}_{W_0})-(-q\cdot\mathds{1}_{W_1^\perp},-\mathds{1}_{W_1}) =(q\cdot\mathds{1}_{W_1^\perp}-\mathds{1}_{W_0^\perp},\mathds{1}_{W_1- W_0}).\]

The proof of the statement about $\ell_{T,n,b}$ is similar.
\end{proof}

\subsection{The preimage of $\Pic_{\can}(\oToySht_T^n)\otimes\mathbb{Z}[\frac{1}{p}]$ in $\mathfrak{O}_T^n\otimes\mathbb{Z}[\frac{1}{p}]$}
We normalize the Haar measure on $T$ by the condition that the measure of any c-lattice $\Lambda$ equals $q^{n(\Lambda)}$.

Fix a nontrivial additive character $\psi:\mathbb{F}_q\to\mathbb{C}^\times$.

The Fourier transform $\Four_\psi$ is defined by the equation (\ref{Fourier transform on a Tate space}).

We have a  homomorphism
\[C_c^\infty(T;\mathbb{Z}[\tfrac{1}{p}])^{\mathbb{F}_q^\times}\xrightarrow{\Four_\psi\oplus 1} C_c^\infty(T^*;\mathbb{Z}[\tfrac{1}{p}])^{\mathbb{F}_q^\times}\oplus C_c^\infty(T;\mathbb{Z}[\tfrac{1}{p}])^{\mathbb{F}_q^\times}\]
and an inclusion of abelian groups
\[C_c^\infty(T^*;\mathbb{Z}[\tfrac{1}{p}])^{\mathbb{F}_q^\times}\oplus C_c^\infty(T;\mathbb{Z}[\tfrac{1}{p}])^{\mathbb{F}_q^\times}\subset
C_+(T^*-\{0\};\mathbb{Z}[\tfrac{1}{p}])^{\mathbb{F}_q^\times}\oplus C_+(T-\{0\};\mathbb{Z}[\tfrac{1}{p}])^{\mathbb{F}_q^\times}.\]
\cref{space of Tate toy horospherical subschemes} gives a homomorphism
\[C_+(T^*-\{0\};\mathbb{Z}[\tfrac{1}{p}])^{\mathbb{F}_q^\times}\oplus C_+(T-\{0\};\mathbb{Z}[\tfrac{1}{p}])^{\mathbb{F}_q^\times}\to
\Pic(\oToySht_T^n)\otimes\mathbb{Z}[\tfrac{1}{p}].\]
The composition of above homomorphisms gives rise to a homomorphism
\begin{equation*}
\gamma: C_c^\infty(T;\mathbb{Z}[\tfrac{1}{p}])^{\mathbb{F}_q^\times}\to\Pic(\oToySht_T^n)\otimes\mathbb{Z}[\tfrac{1}{p}].
\end{equation*}
\begin{Remark}
\cref{space of principal rational Tate toy horospherical divisors} implies that $\ker\gamma=C_0^\infty(T-\{0\};\mathbb{Z}[\tfrac{1}{p}])^{\mathbb{F}_q^\times}$. Also note that $C_0^\infty(T-\{0\};\mathbb{Z}[\tfrac{1}{p}])^{\mathbb{F}_q^\times}=\{f\in C_c^\infty(T;\mathbb{Z}[\tfrac{1}{p}])^{\mathbb{F}_q^\times}|f(0)=0,\int_Tf(v)dv=0\}$.

\end{Remark}

\begin{Proposition}\label{two-dimensional canonical subspace of the Picard group of the scheme of Tate toy shtukas}
For $f\in C_c^\infty(T;\mathbb{Z}[\tfrac{1}{p}])^{\mathbb{F}_q^\times}$, let $a(f)=\int_T f(v)dv$ and $b(f)=f(0)$. Then
\[(q-1)\gamma(f)=a(f)[\ell_{T,n,a}]-b(f)[\ell_{T,n,b}].\]
\end{Proposition}
\begin{Remark}
In the definition of $a(f)$, we normalize the Haar measure on $T$ by the condition that the measure of any c-lattice $\Lambda$ equals $q^{n(\Lambda)}$.
\end{Remark}
\begin{proof}[Proof of \cref{two-dimensional canonical subspace of the Picard group of the scheme of Tate toy shtukas}]
The statement follows from \cref{divisor of 2 line bundles on the scheme of Tate toy shtukas} and the facts that
\begin{align*}
&a(\mathds{1}_{W_1- W_0})=q-1,& &b(\mathds{1}_{W_1- W_0})=0,\\
&a(-q\cdot\mathds{1}_{W_{-1}}+\mathds{1}_{W_0})=0,& &b(-q\cdot\mathds{1}_{W_{-1}}+\mathds{1}_{W_0})=-(q-1),\\
&q\cdot\mathds{1}_{W_1^\perp}-\mathds{1}_{W_0^\perp}=\Four_\psi(\mathds{1}_{W_1- W_0}), &
&-\mathds{1}_{W_{-1}^\perp- W_0^\perp}=\Four_\psi(-q\cdot\mathds{1}_{W_{-1}}+\mathds{1}_{W_0}).
\end{align*}
\end{proof}

\begin{Definition}
Let $\Pic_{\can}(\oToySht_T^n)$ be the subgroup of $\Pic(\oToySht_T^n)$ generated by $[\ell_{T,n,a}], [\ell_{T,n,b}]$ and $[\ell_{T,n,\det}^W]$, where $W$ is a c-lattice of $T$.
\end{Definition}

\begin{Remark}
\cref{isomorphism between determinant line bundles for different c-lattices} implies that the class of $\ell_{T,n,\det}^W$ is independent of the choice of the c-lattice $W$. Hence $\Pic_{\can}(\oToySht_T^n)$ is well-defined.
\end{Remark}
\begin{Remark}
\cref{divisor of 2 line bundles on the scheme of Tate toy shtukas}, \cref{description of Schubert divisors of the scheme of Tate toy shtukas} and \cref{space of Tate toy horospherical subschemes} show that the image of $\mathfrak{O}_T^n$ in $\Pic(\oToySht_T^n)$ contains $\Pic_{\can}(\oToySht_T^n)$.
\end{Remark}

As before, we identify $\mathfrak{O}_T^n\otimes\mathbb{Z}[\frac{1}{p}]$ with $C_+(T^*-\{0\};\mathbb{Z}[\frac{1}{p}])^{\mathbb{F}_q^\times}\oplus C_+(T-\{0\};\mathbb{Z}[\frac{1}{p}])^{\mathbb{F}_q^\times}$ via \cref{space of Tate toy horospherical subschemes}, and we consider $C_c^\infty(T^*;\mathbb{Z}[\frac{1}{p}])^{\mathbb{F}_q^\times}\oplus C_c^\infty(T;\mathbb{Z}[\frac{1}{p}])^{\mathbb{F}_q^\times}$ as an ordered subgroup of $\mathfrak{O}_T^n\otimes\mathbb{Z}[\frac{1}{p}]$.

\begin{Theorem}\label{preimage of two-dimensional canonical subspace of the Picard group of the scheme of Tate toy shtukas}
The preimage of $\Pic_{\can}(\oToySht_T^n)\otimes\mathbb{Z}[\frac{1}{p}]$ in $\mathfrak{O}_T^n\otimes\mathbb{Z}[\frac{1}{p}]$ is the ordered $\mathbb{Z}[\frac{1}{p}]$-submodule
\[\{(f_1,f_2)\in C_c^\infty(T^*;\mathbb{Z}[\tfrac{1}{p}])^{\mathbb{F}_q^\times}\oplus C_c^\infty(T;\mathbb{Z}[\tfrac{1}{p}])^{\mathbb{F}_q^\times}| f_1=\Four_\psi(f_2)\}.\]
\end{Theorem}
\begin{proof}
Let $W_0, W_1, W_{-1}$ be as in \cref{divisor of 2 line bundles on the scheme of Tate toy shtukas}. By \cref{divisor of 2 line bundles on the scheme of Tate toy shtukas,description of Schubert divisors of the scheme of Tate toy shtukas}, the preimage is the $\mathbb{Z}[\frac{1}{p}]$-span of the three elements $(\mathds{1}_{W_0^\perp},\mathds{1}_{W_0})$, $(q\cdot\mathds{1}_{W_1^\perp}-\mathds{1}_{W_0^\perp},\mathds{1}_{W_1\backslash W_0})$, $(-\mathds{1}_{W_{-1}^\perp\backslash W_0^\perp},-q\cdot\mathds{1}_{W_{-1}}+\mathds{1}_{W_0})$ and the $\mathbb{Z}[\frac{1}{p}]$-module $\mathfrak{R}_T^n\otimes\mathbb{Z}[\frac{1}{p}]$. Note that the three elements are contained in the above $\mathbb{Z}[\frac{1}{p}]$-submodule of $\mathfrak{O}_T^n\otimes\mathbb{Z}[\frac{1}{p}]$. So the statement follows from the description of $\mathfrak{R}_T^n\otimes\mathbb{Z}[\frac{1}{p}]$ in \cref{space of principal rational Tate toy horospherical divisors}.
\end{proof}

\section{Review of Drinfeld shtukas}
\subsection{Notation and conventions}\label{(Section)notation and conventions for shtukas}
The following notation and conventions will be used in the rest of the article.

Let $X$ be a smooth projective geometrically connected curve over $\mathbb{F}_q$. Let $k$ be the field of rational functions on $X$. Let $\mathbb{A}$ be the ring of adeles of $k$. Let $O$ be the subring of integral adeles in $\mathbb{A}$.

For any scheme $S$ over $\mathbb{F}_q$, denote $\Phi_S=\Id_X\times \Fr_S:X\times S\to X\times S$, and let $\pi_S:X\times S\to S$ be the projection. We sometimes write $\Phi$ and $\pi$ instead of $\Phi_S$ and $\pi_S$ when there is no ambiguity about $S$.

For a scheme $S$ over $\mathbb{F}_q$ and two morphisms $\alpha,\beta:S\to X$, we say that they satisfy condition (\ref{Frobenius shifts of zero and pole mutually disjoint}) if
\begin{equation}\label{Frobenius shifts of zero and pole mutually disjoint}
\text{$\Gamma_{\Fr_X^i\bcirc\alpha},\Gamma_{\Fr_X^j\bcirc\beta},(i,j\in \mathbb{Z}_{\ge 0})$ are mutually disjoint subsets of $X\times S$.} \tag{$\ast$}
\end{equation}
We say that they satisfy condition (\ref{zero and pole not related by Frobenius}) if
\begin{equation}\label{zero and pole not related by Frobenius}
\text{$\Gamma_{\Fr_X^i\bcirc\alpha}\cap\Gamma_{\Fr_X^j\bcirc\beta}=\emptyset$ for all $i,j\in \mathbb{Z}_{\ge 0}$.} \tag{$+$}
\end{equation}
Condition (\ref{Frobenius shifts of zero and pole mutually disjoint}) is equivalent to the combination of (\ref{zero and pole not related by Frobenius}) and the condition that $\alpha$ and $\beta$ map $S$ to the generic point of $X$.

For a morphism $\alpha:S\to X$, we sometimes write $\alpha$ instead of $\Gamma_\alpha$.

\subsection{Definition of Drinfeld shtukas}
Let $S$ be a scheme over $\mathbb{F}_q$. Denote $\Phi=\Id_X\times\Fr_S:X\times S\to X\times S$.
\begin{Definition}[Drinfeld]
  A \emph{left shtuka} of rank $d$ over $S$ is a diagram
  \[\Phi^*\mathscr{F}\xhookleftarrow{j}\mathscr{F}'\xhookrightarrow{i}\mathscr{F}\]
  of locally free sheaves of rank $d$ on $X\times S$, where $i$ and $j$ are  injective morphisms, the cokernel of $i$ is an invertible sheaf on the graph $\Gamma_\alpha$ of some morphism $\alpha:S\to X$, the cokernel of $j$ is an invertible sheaf on the graph $\Gamma_\beta$ of some morphism $\beta:S\to X$.

  A \emph{right shtuka} of rank $d$ over $S$ is a diagram
  \[\Phi^*\mathscr{F}\xhookrightarrow{f}\mathscr{F}'\xhookleftarrow{g}\mathscr{F}\]
  of locally free sheaves of rank $d$ on $X\times S$, where $f$ and $g$ are  injective morphisms, the cokernel of $f$ is an invertible sheaf on the graph $\Gamma_\alpha$ of some morphism $\alpha:S\to X$, the cokernel of $g$ is an invertible sheaf on the graph $\Gamma_\beta$ of some morphism $\beta:S\to X$.

  We say that $\alpha$ is the zero of the shtuka and $\beta$ is the pole of the shtuka.
\end{Definition}

\begin{Definition}[Drinfeld]
For two morphisms $\alpha,\beta:S\to X$, a \emph{shtuka} of rank $d$ over $S$ with zero $\alpha$ and pole $\beta$ is a locally free sheaf $\mathscr{F}$ of rank $d$ on $X\times S$ equipped with an injective morphism $\Phi^*\mathscr{F}\to\mathscr{F}(\Gamma_\beta)$ inducing an isomorphism $\Phi^*\det\mathscr{F}\xrightarrow{\sim}(\det\mathscr{F})(\Gamma_\beta-\Gamma_\alpha)$, such that the image of the composition $\Phi^*\mathscr{F}\to\mathscr{F}(\Gamma_\beta)\to\mathscr{F}(\Gamma_\beta)/\mathscr{F}$ has rank at most 1 at $\Gamma_\beta$.
\end{Definition}

\begin{Remark}
  When the zero and the pole have disjoint graph, Construction A in \cref{(Section)general constructions for shtukas} shows that there is no difference between a shtuka, a left shtuka and a right shtuka. This remains true after applying any powers of partial Frobeniuses (see \cref{(Section)partial Frobenius}) if we impose condition (\ref{zero and pole not related by Frobenius}) on the zero and the pole.

  From now on, when the zero and the pole satisfy condition (\ref{zero and pole not related by Frobenius}), we do not distinguish left and right shtukas, and simply call them shtukas.
\end{Remark}

\begin{Lemma}
  Let $\mathscr{M}$ be a quasi-coherent sheaf on $X$ and let $\widetilde{\mathscr{M}}$ be the pullback of $\mathscr{M}$ under the projection $X\times S\to X$. Then we have a canonical isomorphism $\Phi^*\widetilde{\mathscr{M}}\cong\widetilde{\mathscr{M}}$. \qed
\end{Lemma}

\begin{Definition}[Drinfeld]
Let $\mathscr{F}$ be a shtuka (resp. a left shtuka, resp. a right shtuka) of rank $d$ over $S$ with zero $\alpha$ and pole $\beta$. Let $D$ be a finite subscheme of $X$ such that $\alpha$ and $\beta$ map to $X-D$. A structure of level $D$ on $\mathscr{F}$ is an isomorphism $\iota:\mathscr{F}\otimes\mathscr{O}_{D\times S}\xrightarrow{\sim}\mathscr{O}_{D\times S}^d$ such that the following diagram commutes:
\[\begin{tikzcd}
  \mathscr{F}\otimes\mathscr{O}_{D\times S}\arrow[r,"\iota","\sim"']&\mathscr{O}_{D\times S}^d\\
  \Phi^*\mathscr{F}\otimes\mathscr{O}_{D\times S}\arrow[r,"\Phi^*\iota","\sim"']\arrow[u,"\sim"]&\Phi^*\mathscr{O}_{D\times S}^d\arrow[u,"\sim"]
\end{tikzcd}\]
\end{Definition}

\subsection{General constructions for shtukas}\label{(Section)general constructions for shtukas}
We have the following constructions for shtukas, which induce morphisms between moduli stacks of shtukas.

\textbf{Construction A:}

(i)  Let $\Phi^*\mathscr{G}\xhookleftarrow{j}\mathscr{F}\xhookrightarrow{i}\mathscr{G}$ be a left shtuka with zero $\alpha$ and pole $\beta$ such that $\Gamma_\alpha\cap\Gamma_\beta=\emptyset$. We form the pushout diagram
\[\begin{tikzcd}
    \mathscr{G}\arrow{r}{g}&\mathscr{H}\\
    \mathscr{F}\arrow{r}{j}\arrow{u}{i}&\Phi^*\mathscr{G}\arrow{u}{f}
\end{tikzcd}\]
Then $\Phi^*\mathscr{G}\xhookrightarrow{f}\mathscr{H}\xhookleftarrow{g}\mathscr{G}$ is a right shtuka of the same rank, with the same zero and pole.

(ii)  Let $\Phi^*\mathscr{G}\xhookrightarrow{f}\mathscr{F}\xhookleftarrow{g}\mathscr{G}$ be a right shtuka with zero $\alpha$ and pole $\beta$ such that $\Gamma_\alpha\cap\Gamma_\beta=\emptyset$. We form the pullback diagram
\[\begin{tikzcd}
    \mathscr{G}\arrow{r}{g}&\mathscr{F}\\
    \mathscr{H}\arrow{r}{j}\arrow{u}{i}&\Phi^*\mathscr{G}\arrow{u}{f}
\end{tikzcd}\]
Then $\Phi^*\mathscr{G}\xhookleftarrow{j}\mathscr{H}\xhookrightarrow{i}\mathscr{G}$ is a left shtuka of the same rank, with the same zero and pole.

\textbf{Construction B:}

(i) For a left shtuka $\Phi^*\mathscr{G}\xhookleftarrow{j}\mathscr{F}\xhookrightarrow{i}\mathscr{G}$ with zero $\alpha$ and pole $\beta$, we can construct a right shtuka $\Phi^*\mathscr{F}\xhookrightarrow{\Phi^*i}\Phi^*\mathscr{G}\xhookleftarrow{j}\mathscr{F}$ of the same rank, with zero $\Fr_X\bcirc\alpha$ and pole $\beta$.

(ii) For a right shtuka $\Phi^*\mathscr{F}\xhookrightarrow{f}\mathscr{G}\xhookleftarrow{g}\mathscr{F}$ with zero $\alpha$ and pole $\beta$, we can construct a left shtuka $\Phi^*\mathscr{G}\xhookleftarrow{\Phi^*g}\Phi^*\mathscr{F}\xhookrightarrow{f}\mathscr{G}$ of the same rank, with zero $\alpha$ and pole $\Fr_X\bcirc\beta$.

\textbf{Construction C:}
For a shtuka $\mathscr{F}$ with zero $\alpha$ and pole $\beta$ satisfying $\Gamma_\alpha\cap\Gamma_\beta=\emptyset$, its dual $\mathscr{F}^\vee$ is a shtuka of the same rank, with zero $\beta$ and pole $\alpha$.

\textbf{Construction D:}
Let $\mathscr{F}$ be a shtuka over $S$. Let $\mathscr{L}$ be an invertible sheaf on $X$, and let $\widetilde{\mathscr{L}}$ be the pullback of $\mathscr{L}$ under the projection $X\times S\to X$. then $\mathscr{F}\otimes\widetilde{\mathscr{L}}$ is a shtuka of the same rank, with the same zero and pole.

\textbf{Construction D':}
Suppose in Construction D the shtuka $\mathscr{F}$ is equipped with a structure of level $D$ and the invertible sheaf $\mathscr{L}$ is trivialized at $D$. Then the shtuka $\mathscr{F}\otimes\widetilde{\mathscr{L}}$ is naturally equipped with a structure of level $D$.

\textbf{Construction E:}
Let $\mathscr{F}$ be a shtuka of rank $d$ equipped with a structure of level $D$. Suppose we are given an $\mathscr{O}_D$-submodule $\mathscr{R}\subset \mathscr{O}_D^d$. Let $\mathscr{F}'$ be the kernel of the composition $\mathscr{F}\to\mathscr{F}\otimes_{\mathscr{O}_{X\times S}}\mathscr{O}_{D\times S}\xrightarrow{\sim}\mathscr{O}_{D\times S}^d\to\mathscr{O}_{D\times S}^d/(\mathscr{R}\boxtimes\mathscr{O}_S)$. Then $\mathscr{F}'$ is a shtuka of the same rank, with the same zero and pole.

\textbf{Construction E':}
Suppose in Construction E we are given a surjective morphism of $\mathscr{O}_D$-modules $\mathscr{R}\to\mathscr{O}_{D'}^d$, where $D'$ is a subscheme of $D$. Then the composition $\mathscr{F}'\to\mathscr{R}\boxtimes\mathscr{O}_S\to\mathscr{O}_{D'\times S}^d$ defines a structure of level $D'$ on $\mathscr{F}'$.

\subsection{Group action}
Let $\Sht_{\all}^d$ be the moduli scheme of shtukas equipped with structures of all levels compatible with each other. We have a left action of $GL_d(\mathbb{A})$ on $\Sht_{\all}^d$ as shown in Section 3 of~\cite{Drinfeld1}. In particular, we have the following statement.

\begin{Lemma}\label{Euler characteristic of shtukas under group action}
   The action of $g\in GL_d(\mathbb{A})$ increases the Euler characteristic of a shtuka by $\deg g$. \qed
\end{Lemma}

\subsection{Partial Frobenius}\label{Section: partial Frobenius}
\begin{Definition}
Let $F_1$ be the construction of first applying A(ii) and then applying B(i). Let $F_2$ be the construction of first applying B(ii) and then applying A(i). They are called partial Frobenius.
\end{Definition}

\begin{Proposition}
For a shtuka $\mathscr{F}$ over $S$ with zero and pole satisfying condition (\ref{zero and pole not related by Frobenius}), we have natural isomorphisms $F_1F_2\mathscr{F}\cong F_2F_1\mathscr{F}\cong\Fr_S^*\mathscr{F}$. \qed
\end{Proposition}

\begin{Remark}
  If a shtuka $\mathscr{F}$ over $S$ has zero $\alpha$ and pole $\beta$, then $F_1\mathscr{F}$ has zero $\alpha\bcirc\Fr_S=\Fr_X\bcirc\alpha$ and pole $\beta$, and $F_2\mathscr{F}$ has zero $\alpha$ and pole $\beta\bcirc\Fr_S=\Fr_X\bcirc\beta$.
\end{Remark}

\section{Reducible shtukas over a field}
In this section, we recollect the results in Section 2 of~\cite{Ding17}.

We use the notation and conventions of \cref{(Section)notation and conventions for shtukas}.

Fix a field $E$ over $\mathbb{F}_q$. Denote $\Phi=\Id_X\otimes\Fr_E:X\otimes E\to X\otimes E$.

\subsection{Definitions}
\begin{Definition}\label{definition of reducible shtukas}
A shtuka $\mathscr{F}$ over $\Spec E$ of rank $d$ with zero $\alpha$ and pole $\beta$ is said to be \emph{reducible} if $\mathscr{F}$ contains a nonzero subsheaf $\mathscr{E}$ of rank $<d$ such that the image of $\Phi^*\mathscr{E}$ in $\mathscr{F}(\beta)$ is contained in $\mathscr{E}(\beta)$.
\end{Definition}

\begin{Remark}\label{two types of exact sequences for reducible shtukas}
  Let the notation be the same as in \cref{definition of reducible shtukas}. Let $\mathscr{G}$ be the saturation of $\mathscr{E}$ in $\mathscr{F}$. We have an exact sequence of locally free sheaves on $X\otimes E$
\[\begin{tikzcd}0\arrow[r]&\mathscr{G}\arrow[r]&\mathscr{F}\arrow[r]&\mathscr{H}\arrow[r]&0.\end{tikzcd}\]
Assume $\alpha\ne \beta$. Then one (and only one) of the following two possibilities holds.

(1) $\mathscr{G}$ is a shtuka with zero $\alpha$ and pole $\beta$, and the morphism $\Phi^*\mathscr{F}\hookrightarrow\mathscr{F}(\beta)$ induces an isomorphism $\Phi^*\mathscr{H}\xrightarrow{\sim}\mathscr{H}$.

(2) $\mathscr{H}$ is a shtuka with zero $\alpha$ and pole $\beta$, and the morphism $\Phi^*\mathscr{F}\hookrightarrow\mathscr{F}(\beta)$ induces an isomorphism $\Phi^*\mathscr{G}\xrightarrow{\sim}\mathscr{G}$.
\end{Remark}

\subsection{Maximal trivial sub and maximal trivial quotient}\label{section: maximal trivial sub and maximal trivial quotient}
The following three statements are proved in Section 2.2 of~\cite{Ding17}.

\begin{Lemma}\label{difference between divisor and its Frobenius}
Let $\alpha,\beta:\Spec E\to X$ be two morphisms. If an effective divisor $D$ of $X\otimes E$ satisfies $\Phi^*D+\alpha=D+\beta$, then $\beta=\Fr_X^n\bcirc\alpha$ for some $n\ge 0$.
\end{Lemma}

\begin{Proposition}\label{Frobenius stable subsheaf (equal rank)}
Let $\mathscr{G}$ be a shtuka over $\Spec E$ with zero $\alpha $ and pole $\beta$. Assume that there exists a subsheaf $\mathscr{E}\subset\mathscr{G}$ satisfying

(i) $\rank\mathscr{E}=\rank\mathscr{G}$;

(ii) the image of $\Phi^*\mathscr{E}$ in $\mathscr{G}(\beta)$ is contained in $\mathscr{E}$.

\noindent Then $\beta=\Fr_X^n\bcirc\alpha$ for some $n\ge 0$.
\end{Proposition}

\begin{Proposition}\label{maximal trivial sub and maximal trivial quotient}
Let $\mathscr{F}$ be a shtuka over $\Spec E$ with zero $\alpha$ and pole $\beta$ satisfying condition (\ref{zero and pole not related by Frobenius}). Let $S_1$ (resp. $S_2$) be the poset of all subsheaves $\mathscr{E}\subset\mathscr{F}$ satisfying the following condition (1) (resp. (2)).

(1) The image of $\Phi^*\mathscr{E}$ in $\mathscr{F}(\beta)$ is contained in $\mathscr{E}(\beta)$, the sheaf $\mathscr{F}/\mathscr{E}$ is locally free, and the morphism $\Phi^*(\mathscr{F}/\mathscr{E})\to(\mathscr{F}/\mathscr{E})(\beta)$ induces an isomorphism $\Phi^*(\mathscr{F}/\mathscr{E})\xrightarrow{\sim}(\mathscr{F}/\mathscr{E})$.

(2) The image of $\Phi^*\mathscr{E}$ in $\mathscr{F}(\beta)$ is $\mathscr{E}$.

Then the poset $S_1$ has a least element, denoted by $\mathscr{F}^{\I}$. The poset $S_2$ has a greatest element, denoted by $\mathscr{F}^{\II}$, and $\mathscr{F}/\mathscr{F}^{\II}$ is locally free.
\end{Proposition}

\begin{Remark}\label{criterion for irreducibility of shtuka using maximal trivial sub and quotient}
Suppose $\mathscr{F}$ is a right shtuka over $\Spec E$ with zero and pole satisfying condition (\ref{zero and pole not related by Frobenius}). Then $\mathscr{F}$ is irreducible if and only if $\mathscr{F}^{\I}=\mathscr{F}$ and $\mathscr{F}^{\II}=0$.
\end{Remark}

\begin{Corollary}\label{twist of maximal trivial sub and maximal trivial quotient}
  Let $\mathscr{F},\mathscr{F}^{\I},\mathscr{F}^{\II}$ be as in \cref{maximal trivial sub and maximal trivial quotient}. Let $\mathscr{M}$ be an invertible sheaf on $X$ and let $\mathscr{G}$ be the shtuka $\mathscr{F}\otimes(\mathscr{M}\otimes E)$ obtained by Construction E in \cref{(Section)general constructions for shtukas}. Then $\mathscr{G}^{\I}=\mathscr{F}^{\I}\otimes(\mathscr{M}\otimes E)$ and $\mathscr{G}^{\II}=\mathscr{F}^{\II}\otimes(\mathscr{M}\otimes E)$. \qed
\end{Corollary}

\section{Relation between shtukas and toy shtukas}\label{(Section)relation between shtukas and toy shtukas}
We use the notation and conventions of \cref{(Section)notation and conventions for shtukas}.

Let $S$ be a scheme over $\mathbb{F}_q$. For an $\mathscr{O}_{X\times S}$-module $\mathscr{F}$ and a point $s\in S$, we denote $\mathscr{F}_s$ to be the pullback of $\mathscr{F}$ to $X\times s$. Put $\Phi=\Id_X\times\Fr_S:X\times S\to X\times S$. Let $\pi:X\times S\to S$ be the projection.

\subsection{Right shtukas}
Let $\Phi^*\mathscr{F}\hookrightarrow\mathscr{F}'\hookleftarrow\mathscr{F}$ be a right shtuka of rank $d$ over $S$ equipped with a structure of level $D$. Suppose that $D$ viewed as an effective divisor is represented as $D''-D'$ so that for every point $s\in S$ one has
\begin{equation}\label{vanishing assumption for H^0 of right shtuka}
H^0(X\times s,\mathscr{F}'_s(D'\times s))=0,
\end{equation}
\begin{equation}\label{vanishing assumption for H^1 of right shtuka}
H^1(X\times s,\mathscr{F}_s(D''\times s))=0.
\end{equation}

Let $V=H^0(X,(\mathscr{O}_X(D'')/\mathscr{O}_X(D'))^d)$, $\mathscr{L}=\pi_*\mathscr{F}(D''\times S)$, $\mathscr{L}'=\pi_*\mathscr{F}'(D''\times S)$.

\begin{Proposition}\label{from a right shtuka to a right toy shtuka}
The pair $\mathscr{L},\mathscr{L}'$ forms a right toy shtuka for $V$ over $S$.
\end{Proposition}
\begin{proof}
Consider the composition
\[\begin{tikzcd}
\mathscr{F}\arrow[r]&\mathscr{F}\otimes\mathscr{O}_{D\times S}\arrow[r,"\sim"] &\mathscr{O}_{D\times S}^d
\end{tikzcd}\]
where the second morphism is the isomorphism from the structure of level $D$. Tensoring with $\mathscr{O}_{X\times S}(D''\times S)$, we get a composition
\[\begin{tikzcd}
\mathscr{F}(D''\times S)\arrow[r] &\mathscr{F}(D''\times S)/\mathscr{F}(D'\times S)\arrow[r,"\sim"]&
(\mathscr{O}_{X\times S}(D''\times S)/\mathscr{O}_{X\times S}(D'\times S))^d .
\end{tikzcd}\]
This induces a morphism $\mathscr{L}\to V\otimes\mathscr{O}_S$.

Similarly one gets a morphism $\mathscr{L}'\to V\otimes\mathscr{O}_S$.

Consider the exact sequence
\[\begin{tikzcd}
  0\arrow[r]& \mathscr{F}'(D'\times S)\arrow[r] & \mathscr{F}'(D''\times S)\arrow[r] & \mathscr{F}'(D''\times S)/\mathscr{F}'(D'\times S)\arrow[r]&0.
\end{tikzcd}\]
Assumption (\ref{vanishing assumption for H^0 of right shtuka}) and base change for cohomology imply that $\pi_*\mathscr{F}'(D'\times S)=0$. Hence the morphism $\mathscr{L}'\to V\otimes\mathscr{O}_S$ is injective.

Since $\mathscr{F}'/\mathscr{F}$ is torsion, assumption (\ref{vanishing assumption for H^1 of right shtuka}) implies $H^1(X\times s,\mathscr{F}'_s(D''\times s))=0$ for every point $s\in S$. From base change for cohomology we get $R^1\pi_*\mathscr{F}'(D''\times S)=0$. So $(V\otimes\mathscr{O}_S)/\mathscr{L}'=R^1\pi_*\mathscr{F}'(D'\times S)$. For any point $s\in S$, we have $H^0(X\times s,\mathscr{F}'_s(D'\times s))=0$ by assumption (\ref{vanishing assumption for H^0 of right shtuka}), so in particular the base change morphism $k(s)\otimes \pi_*\mathscr{F}'(D'\times s)\to H^0(X\times s,\mathscr{F}'_s(D'\times s))$ is surjective. By Theorem 12.11(b) of Chapter 3 of~\cite{Hart}, $(V\otimes\mathscr{O}_S)/\mathscr{L}'=R^1\pi_*\mathscr{F}'(D'\times S)$ is locally free.

A similar argument shows that $(V\otimes\mathscr{O}_S)/\mathscr{L}$ is locally free.

Since $\mathscr{F}\subset \mathscr{F}'$, we have $\mathscr{L}\subset\mathscr{L}'$. Since $\Phi^*\mathscr{F}\subset\mathscr{F}'$, we have $\Fr_S^*\mathscr{L}\subset\mathscr{L}'$.

Assumption (\ref{vanishing assumption for H^1 of right shtuka}) and base change for cohomology imply that $R^1\pi_*\mathscr{F}(D''\times S)$. Since $\Gamma_\beta\cap(D''\times S)=\emptyset$, we have an isomorphism $\mathscr{F}'(D''\times S)/\mathscr{F}(D''\times S)\cong \mathscr{F}'/\mathscr{F}$. So we have an exact sequence
\[\begin{tikzcd}
0\arrow[r]& \mathscr{L}\arrow[r]& \mathscr{L}'\arrow[r]& \pi_*(\mathscr{F}'/\mathscr{F})\arrow[r]& 0.
\end{tikzcd}\]
Since $\mathscr{F}'/\mathscr{F}$ is an invertible sheaf on $\Gamma_\beta$, and $\pi:X\times S\to S$ induces an isomorphism $\Gamma_\beta\to S$, we see that $\mathscr{L}'/\mathscr{L}$ is invertible.

Similarly, $\mathscr{L}'/\Fr_S^*\mathscr{L}$ is also invertible.
\end{proof}
\subsection{Left shtukas}
Let $\Phi^*\mathscr{F}\hookleftarrow\mathscr{F}'\hookrightarrow\mathscr{F}$ be a left shtuka of rank $d$ over $S$ equipped with a structure of level $D$. Suppose that $D$ viewed as an effective divisor is represented as $D''-D'$ so that for every point $s\in S$ one has
\begin{equation}\label{vanishing assumption for H^0 of left shtuka}
H^0(X\times s,\mathscr{F}_s(D'\times s))=0,
\end{equation}
\begin{equation}\label{vanishing assumption for H^1 of left shtuka}
H^1(X\times s,\mathscr{F}'_s(D''\times s))=0.
\end{equation}

Let $V=H^0(X,(\mathscr{O}_X(D'')/\mathscr{O}_X(D'))^d)$, $\mathscr{L}=\pi_*\mathscr{F}(D''\times S)$, $\mathscr{L}'=\pi_*\mathscr{F}'(D''\times S)$.

Similarly to the case of a right shtuka, one can prove the following statement.

\begin{Proposition}\label{from a left shtuka to a left toy shtuka}
The pair $\mathscr{L},\mathscr{L}'$ forms a left toy shtuka for $V$ over $S$. \qed
\end{Proposition}
\subsection{Shtukas}
Let $\mathscr{F}$ be a shtuka of rank $d$ over $S$ equipped with a structure of level $D$. Suppose that $D$ viewed as an effective divisor is represented as $D''-D'$ so that for every point $s\in S$ one has
\begin{equation}\label{vanishing assumption for H^0 of shtuka}
H^0(X\times s,\mathscr{F}_s(D'\times s))=0,
\end{equation}
\begin{equation}\label{vanishing assumption for H^1 of shtuka}
H^1(X\times s,\mathscr{F}_s(D''\times s))=0.
\end{equation}

Let $V=H^0(X,(\mathscr{O}_X(D'')/\mathscr{O}_X(D'))^d)$, $\mathscr{L}=\pi_*\mathscr{F}(D''\times S)$.

Similarly to the case of a right shtuka, one can prove the following statement.

\begin{Proposition}\label{from a shtuka to a toy shtuka}
$\mathscr{L}$ forms a toy shtuka for $V$ over $S$. \qed
\end{Proposition}

\section{The morphism from the moduli scheme of shtukas with structures of all levels to the moduli scheme of Tate toy shtukas}\label{(Section)from a shtuka with all level structures to a Tate toy shtuka}
We use the notation and conventions of \cref{(Section)notation and conventions for shtukas}.

Fix a field $E$ over $\mathbb{F}_q$ and two morphisms $\alpha,\beta:\Spec E\to X$ satisfying condition (\ref{zero and pole not related by Frobenius}). Let $d$ be a positive integer.

Let $\Sht_{E,\all}^d$ denote the moduli scheme of shtukas over $\Spec E$ with zero $\alpha$ and pole $\beta$ equipped with structures of all levels compatible with each other. Let $\Sht_{E,\all}^{d,\chi}$ be the components of $\Sht_{E,\all}^d$ on which the shtuka $\mathscr{F}$ has Euler characteristic $\chi$.

For any integer $\chi$, let $n_\chi\in\Dim_{\mathbb{A}^d}$ be such that $n_\chi(O^d)=\chi$.

For a divisor $D$ of $X$, we denote $O_D$ to be the c-lattice
\[\varprojlim_{D'}H^0(X,\mathscr{O}_X(D)/\mathscr{O}_X(D'))\]
of $\mathbb{A}$, where $D'$ runs through all divisors of $X$ such that $D'\le D$. In other words, $O_D$ consists of those adeles with poles bounded by $D$. For two divisors $D', D''$ of $X$ such that $D'\le D''$, we have
\[O_{D''}/O_{D'}=H^0(X,\mathscr{O}_X(D'')/\mathscr{O}_X(D')).\]
We have an isomorphism
\[\varinjlim_{D''}\varprojlim_{D'\le D''} O_{D''}/O_{D'}\xrightarrow{\sim}\mathbb{A}.\]

\subsection{Construction of the morphism $\theta$}
We first construct the morphism from $\Sht_{E,\all}^{d,\chi}$ to $\ToySht_{\mathbb{A}^d}^{n_\chi}$.

Let $S$ be a scheme over $\Spec E$ and let $\mathscr{F}\in\Sht_{E,\all}^{d,\chi}(S)$. Let $\pi:X\times S\to S$ be the projection.

Let $S^{D',D''}$ be the open subscheme of $S$ such that all $s\in S^{D',D''}$ satisfy conditions (\ref{vanishing assumption for H^0 of shtuka}) and (\ref{vanishing assumption for H^1 of shtuka}). \cref{from a shtuka to a toy shtuka} shows that $\pi_*\mathscr{F}(D'')$ is a toy shtuka over $S$ for $O_{D''}^d/O_{D'}^d$. Moreover, conditions (\ref{vanishing assumption for H^0 of shtuka}) and (\ref{vanishing assumption for H^1 of shtuka}) imply that $\pi_*\mathscr{F}(D'')$ has rank $\chi+d\cdot\deg D''$.

For divisors $\widetilde{D'},\widetilde{D''}$ such that $\widetilde{D'}\le D'\le D''\le \widetilde{D''}$, we have $S^{D',D''}\subset S^{\widetilde{D'},\widetilde{D''}}$, and the composition
\[S^{\widetilde{D'},\widetilde{D''}}\to\ToySht_{O_{\widetilde{D'}}^d/O_{\widetilde{D''}}^d}^{\chi+d\cdot\deg\widetilde{D''}} \to\ToySht_{O_{D''}^d/O_{D'}^d}^{\chi+d\cdot\deg D''}\]
when restricted to $S^{D',D''}$ coincides with the morphism $S^{D',D''}\to\ToySht_{O_{D''}^d/O_{D'}^d}^{\chi+d\cdot\deg D''}$.

For each $s\in S$, there exists a pair of divisors $D'\le D''$ such that $s\in S^{D',D''}$. Passing to the double limit, we see that
\[\mathscr{L}=\varinjlim_{D''}\pi_*\mathscr{F}(D'')\]
is a Tate toy shtuka over $S$ of dimension $n_\chi$ for $\mathbb{A}^d$.

\begin{Proposition}\label{construction of the morphism theta}
  For each $\chi\in\mathbb{Z}$, the above construction induces a morphism
  \[\theta_E^{d,\chi}: \Sht_{E,\all}^{d,\chi}\to\oToySht_{\mathbb{A}^d}^{n_\chi}.\]
\end{Proposition}
\begin{proof}
The above construction induces a morphism $\Sht_{E,\all}^{d,\chi}\to\ToySht_{\mathbb{A}^d}^{n_\chi}$.

Put $\mathcal{M}=\Sht_{E,\all}^{d}$. Let $\mathscr{F}$ be the universal right shtuka over $\mathcal{M}$. Denote $\Phi=\Id_X\times\Fr_\mathcal{M}:X\times \mathcal{M}\to X\times \mathcal{M}$.

For any point $s\in \mathcal{M}$, there exists a divisor $D''\subset X$ such that $\mathscr{F}_s(D''\times s)$ and $\Phi_s^*\mathscr{F}_s(D''\times s)$ are generated by their global sections. Since $\alpha\ne\beta$, we have $\mathscr{F}_s\ne\Phi^*\mathscr{F}_s$. Hence ${\pi_s}_*\mathscr{F}_s(D''\times s)\ne{\pi_s}_*\Phi^*\mathscr{F}_s(D''\times s)=\Fr_s^*{\pi_s}_*\mathscr{F}_s(D''\times s)$. This shows that the image of $s$ in $\ToySht_{\mathbb{A}^d}^{n_\chi}$ is contained in $\oToySht_{\mathbb{A}^d}^{n_\chi}$.
\end{proof}

\subsection{$GL_d(\mathbb{A})$-equivariance of the morphism $\theta$}
For $g\in GL_d(\mathbb{A}),L\in\Grass_{\mathbb{A}^d}(R)$, where $R$ is an $\mathbb{F}_q$-algebra, we define $g\cdot L$ to be the image of $L\hookrightarrow \mathbb{A}^d\widehat{\otimes}R\xrightarrow{g\otimes 1} \mathbb{A}^d\widehat{\otimes}R$. In this way we get a (left) action of $GL_d(\mathbb{A})$ on $\Grass_{\mathbb{A}^d}$. We see that this action preserves $\oToySht_{\mathbb{A}^d}$.

\begin{Proposition}
  The morphism $\theta_E^d:\Sht_{E,\all}^d\to\oToySht_{\mathbb{A}^d}$ is $GL_d(\mathbb{A})$-equivariant.
\end{Proposition}
\begin{proof}
Let $S$ be a scheme over $E$. Let $\mathscr{F}\in\Sht_{E,\all}^d(S)$.

For any $g\in GL_d(\mathbb{A})$, the definition of the action of $g$ in Section 3 of~\cite{Drinfeld1} implies that the following diagram commutes.
\[\begin{tikzcd}[column sep=tiny]
  \pi_*(\varinjlim\limits_{D''}\varprojlim\limits_{D'\le D''}(g^*\mathscr{F})(D''\times S)/(g^*\mathscr{F})(D'\times S)) \arrow[r]\arrow[d,"\sim"]&
  \pi_*(\varinjlim\limits_{D''}\varprojlim\limits_{D'\le D''}\mathscr{F}(D''\times S)/\mathscr{F}(D'\times S)) \arrow[d,"\sim"]\\
  \pi_*(\varinjlim\limits_{D''}\varprojlim\limits_{D'\le D''}(\mathscr{O}_{X\times S}(D''\times S)/\mathscr{O}_{X\times S}(D'\times S))^d) \arrow[d,"="]&
  \pi_*(\varinjlim\limits_{D''}\varprojlim\limits_{D'\le D''}(\mathscr{O}_{X\times S}(D''\times S)/\mathscr{O}_{X\times S}(D'\times S))^d) \arrow[d,"="] \\
  \mathbb{A}^d\widehat{\otimes}R \arrow[r,"g^{-1}\widehat{\otimes}1"] & \mathbb{A}^d\widehat{\otimes}R
\end{tikzcd}\]

The natural morphism
\[\varinjlim_{D''}(g^*\mathscr{F})(D''\times S)\to\varinjlim_{D''}\mathscr{F}(D''\times S)\]
induces an isomorphism
\[\varinjlim_{D''}\pi_*((g^*\mathscr{F})(D''\times S))\xrightarrow{\sim}\varinjlim_{D''}\pi_*(\mathscr{F}(D''\times S))\]
The statement follows.
\end{proof}

\section{Review of horospherical divisors}\label{(Section)review of horospherical divisors}
The goal of \cref{(Section)review of horospherical divisors,(Section)pullback of toy horospherical divisors under theta} is to prove \cref{formula for pullback of toy horospherical divisors under theta}, which relates horospherical divisors on the moduli scheme of shtukas with Tate toy horospherical divisors on the moduli scheme of Tate toy shtukas, and reduces algebraic geometry to representation theory.

We use the notation and conventions of \cref{(Section)notation and conventions for shtukas}.

Let $\eta$ denote the generic point of $X\times X$. Let $\alpha,\beta:\eta\to X$ be the first and second projection. In this section, all shtukas will have zero $\alpha$ and pole $\beta$.

For a finite subscheme $D\subset X$, let $\Sht_{\eta,D}^d$ denote the moduli stack which to each scheme $S$ over $\eta$ associates the groupoid of shtukas over $S$ with zero $\alpha$ and pole $\beta$ equipped with a structure of level $D$.

Let $\Sht_{\eta,\all}^d$ denote the moduli scheme which to each scheme $S$ over $\eta$ associates the set of isomorphism classes of shtukas over $S$ with zero $\alpha$ and pole $\beta$ equipped with structures of all levels compatible with each other.

For $\chi\in\mathbb{Z}$, we denote $n_\chi$ to be the element of $\Dim_{\mathbb{A}^d}$ such that $n_\chi(O^d)=\chi$.

Let $\mathfrak{Vect}_{X,D}^d$ (resp. $\mathfrak{Vect}_{X,\all}^d$) denote the set of isomorphism classes of locally free sheaves of rank $d$ on $X$ equipped with a structure of level $D$ (resp. equipped with structures of all levels compatible with each other).

\subsection{Trivial shtukas}
\begin{Definition}
For $d\ge 1$ and an effective divisor $D\subset X$, let $\TrSht_D^d$ denote the moduli stack which to each scheme $S$ over $\mathbb{F}_q$ associates the groupoid of locally free sheaves $\mathscr{M}$ on $X\times S$ equipped with the following data:

(i) a structure of level $D$, i.e., an isomorphism
$\gamma:\mathscr{M}\otimes\mathscr{O}_{D\times S}\xrightarrow{\sim}\mathscr{O}_{D\times S}^d$;

(ii) an isomorphisms $\Phi_S^*\mathscr{M}\xrightarrow{\sim}\mathscr{M}$ such that the diagram
\[\begin{tikzcd}
  \mathscr{M}\otimes\mathscr{O}_{D\times S} \arrow[r,"\gamma","\sim"']& \mathscr{O}_{D\times S}^d \\
  \Phi_S^*\mathscr{M}\otimes\mathscr{O}_{D\times S} \arrow[r,"\Phi_S^*\gamma","\sim"']\arrow[u,"\sim"]& \Phi_S^*\mathscr{O}_{D\times S}^d \arrow[u,"\sim"]
\end{tikzcd}\]
commutes.

An element of $\TrSht_D(S)$ is called a \emph{trivial shtuka} over $S$ equipped with a structure of level $D$.
\end{Definition}

For $\mathscr{E}\in\mathfrak{Vect}_{X,D}^d$, let $\TrSht_{\mathscr{E},D}$ denote the quotient stack $[\Spec\mathbb{F}_q/\Aut\mathscr{E}]$.

The following statement is Theorem 2 of Section 3 of Chapter I of~\cite{L.Lafforgue97}.

\begin{Proposition}
  The stack $\TrSht_D^d$ is a disjoint union
  \[\TrSht_{D}^d=\coprod_{\mathscr{E}\in\mathfrak{Vect}_{D}^d}\TrSht_{\mathscr{E},D}.\]
\end{Proposition}

\begin{Proposition}\label{descent of coherent sheaves for projective scheme over finite field}
  Let $S$ be a projective scheme over $\mathbb{F}_q$. Let $E$ be an separably closed field over $\mathbb{F}_q$. Then the functor $\mathscr{F}\mapsto\mathscr{F}\otimes E$ is an equivalence between the category of coherent sheaves $\mathscr{F}$ on $S$ and the category of coherent sheaves $\mathscr{M}$ on $S\otimes E$ equipped with an isomorphism $(\Id_S\otimes\Fr_E)^*\mathscr{M}\simto\mathscr{M}$. \qed
\end{Proposition}
\begin{Remark}
  When $E$ is algebraically closed, the above statement is Proposition 1.1 of~\cite{Drinfeld1}. The same proof applies when $E$ is separably closed.
\end{Remark}

\begin{Corollary}\label{descent of locally free sheaves with level structure}
If $E$ is separably closed and $\mathscr{E}\in\mathfrak{Vect}_{X,D}^d$, for any $\mathscr{M}\in\TrSht_{\mathscr{E},D}(E)$, we can find an isomorphism $\mathscr{M}\simto\mathscr{E}\otimes E$ compatible with the structure of level $D$.
\end{Corollary}

\begin{Lemma}\label{non-containment for subshtuka of the same rank}
  Let $E$ be a field over $\mathbb{F}_q$. Let $\mathscr{F}$ be a shtuka over $\Spec E$ with zero $\alpha$ and pole $\beta$ satisfying condition (\ref{Frobenius shifts of zero and pole mutually disjoint}). Let $\mathscr{G}$ be a subshtuka of $\mathscr{F}$ of the same rank with the same zero and pole. Then $\Phi_E^*\mathscr{G}\not\subset\mathscr{F}$ and $\mathscr{G}\not\subset\Phi_E^*\mathscr{F}$.
\end{Lemma}
\begin{proof}
   We apply $d$-th exterior power to all sheaves involved to reduce the problem to the case $d=1$.

   Suppose $\Phi_E^*\mathscr{G}\subset \mathscr{F}$. Then we have $\mathscr{F}=\mathscr{G}(\beta+W)$ for some effective divisor $W$ of $X\otimes E$. From the isomorphisms $\Phi_E^*\mathscr{G}\cong\mathscr{G}(\beta-\alpha)$ and $\Phi_E^*\mathscr{F}\cong\mathscr{F}(\beta-\alpha)$ we deduce that $\beta-\alpha+\Fr_X\bcirc\beta+\Phi^*W=\beta+W+\beta-\alpha$. Hence $\beta+W=\Fr_X\bcirc\beta+\Phi^*W$. Applying \cref{difference between divisor and its Frobenius} to the two morphisms $\beta,\Fr_X\bcirc\beta:\Spec E\to X$, we see that $\beta=\Fr_X^i\bcirc\beta$ for some $i\ge 1$, a contradiction to condition (\ref{Frobenius shifts of zero and pole mutually disjoint}).

   The proof of the second statement is similar.
\end{proof}

\begin{Lemma}\label{subshtuka of the same rank}
  Let $E$ be an algebraically closed field. Let $\mathscr{F}$ be a shtuka of rank $d$ over $\Spec E$ with zero and pole satisfying condition (\ref{Frobenius shifts of zero and pole mutually disjoint}). Let $\mathscr{G}$ be a subshtuka of $\mathscr{F}$ of the same rank with the same zero and pole. Then $\mathscr{F}/\mathscr{G}$ is supported on $D\otimes E$ for some finite subscheme $D\subset X$. Moreover, for any structure of level $D$ on $\mathscr{F}$, $\mathscr{G}$ is obtained from $\mathscr{F}$ by applying Construction E in \cref{(Section)general constructions for shtukas} with respect to that level structure and an $\mathscr{O}_D$-submodule $\mathscr{R}\subset\mathscr{O}_D^d$.
\end{Lemma}
\begin{proof}
  Let $\mathscr{F}'=\Phi_E^*\mathscr{F}+\mathscr{F},\mathscr{G}'=\Phi_E^*\mathscr{G}+\mathscr{G}$. \cref{non-containment for subshtuka of the same rank} shows that $\mathscr{F}\cap\mathscr{G}'=\mathscr{G}$ and $\Phi_E^*\mathscr{F}\cap\mathscr{G}'=\Phi_E^*\mathscr{G}$. Thus the morphisms $\mathscr{F}/\mathscr{G}\to\mathscr{F}'/\mathscr{G}'$ and $\Phi_E^*(\mathscr{F}/\mathscr{G})\to\mathscr{F}'/\mathscr{G}'$ are injective. The sheaves $\mathscr{F}/\mathscr{G}, \mathscr{F}'/\mathscr{G}', \Phi_E^*(\mathscr{F}/\mathscr{G})$ are torsion sheaves on $X\otimes E$, and we have $h^0(\mathscr{F}/\mathscr{G})=h^0(\mathscr{F}'/\mathscr{G}')=h^0(\Phi_E^*(\mathscr{F}/\mathscr{G}))$. Hence the morphisms $\mathscr{F}/\mathscr{G}\to\mathscr{F}'/\mathscr{G}'$ and $\Phi_E^*(\mathscr{F}/\mathscr{G})\to\mathscr{F}'/\mathscr{G}'$ are isomorphisms. So we have $\Phi_E^*(\mathscr{F}/\mathscr{G})\cong\mathscr{F}/\mathscr{G}$.
  \cref{descent of coherent sheaves for projective scheme over finite field} gives an isomorphism $\mathscr{F}/\mathscr{G}\cong\mathscr{M}\otimes E$ for some coherent sheaf $\mathscr{M}$ on $X$. Since $\mathscr{F}$ and $\mathscr{G}$ have the same rank, $\mathscr{M}$ is supported on a finite subscheme $D\subset X$.

  Equip $\mathscr{F}$ with a structure of level $D$. Let $\mathscr{P}=\mathscr{G}/\mathscr{F}(-D\otimes E)\subset\mathscr{F}/\mathscr{F}(-D\otimes E)\cong\mathscr{O}_{D\otimes E}^d$. Since $\mathscr{F}(-D\otimes E)$ is a subshtuka of $\mathscr{G}$, we get an isomorphism $\Phi_E^*\mathscr{P}\simto\mathscr{P}$ which is compatible with the natural isomorphism $\Phi_E^*\mathscr{O}_{D\otimes E}^d\simto\mathscr{O}_{D\otimes E}^d$. \cref{descent of coherent sheaves for projective scheme over finite field} gives an isomorphism $\mathscr{P}\cong\mathscr{R}\otimes E$ for some $\mathscr{O}_D$-submodule $\mathscr{R}\subset\mathscr{O}_D^d$. We see that $\mathscr{G}$ is obtained from $\mathscr{F}$ by applying Construction E with respect to $\mathscr{R}$.
\end{proof}

\subsection{Definition of horospherical cycles}
For $\mathscr{E}\in\mathfrak{Vect}_{X,D}^i$, denote $\TrSht_{\eta,\mathscr{E},D}$ to be the base change of $\TrSht_{\mathscr{E},D}$ from $\Spec \mathbb{F}_q$ to $\eta$.

\begin{Definition}\label{definition of the stack of reducible shtukas of type I}
Given $d\ge 2$, $1\le i\le d-1$ and $\mathscr{E}\in\mathfrak{Vect}_{X,D}^i$, we denote  $\RedSht_{\eta,\mathscr{E},D}^{d,i,\I}$ to be the moduli stack which to each scheme $S$ over $\eta$ associates the groupoid of exact sequences
\[\begin{tikzcd}
0 \arrow[r]& \mathscr{A} \arrow[r]& \mathscr{F} \arrow[r]& \mathscr{B} \arrow[r]& 0,
\end{tikzcd}\]
where

(i) $\mathscr{A}\in\Sht_{\eta,D}^{d-i}(S)$, $\mathscr{F}\in\Sht_{\eta,D}^{d,\chi=0}(S)$, $\mathscr{B}\in\TrSht_{\eta,\mathscr{E},D}(S)$;

(ii) the morphisms $\mathscr{A}\to\mathscr{F}$ and $\mathscr{F}\to\mathscr{B}$ are morphisms of shtukas with structures of level $D$;

(iii) the structures of level $D$ give the following commutative diagram
\[\begin{tikzcd}
0 \arrow[r]& \mathscr{A}\otimes\mathscr{O}_{D\times S} \arrow[r]\arrow[d,"\sim"]&\mathscr{F}\otimes\mathscr{O}_{D\times S} \arrow[r]\arrow[d,"\sim"]& \mathscr{B}\otimes\mathscr{O}_{D\times S} \arrow[r] \arrow[d,"\sim"] & 0 \\
0 \arrow[r]& \mathscr{O}_{D\times S}^{d-i} \arrow[r]& \mathscr{O}_{D\times S}^d=\mathscr{O}_{D\times S}^{d-i}\oplus\mathscr{O}_{D\times S}^i \arrow[r]& \mathscr{O}_{D\times S}^i \arrow[r]& 0
\end{tikzcd}\]
where the lower exact sequence is the standard one.
\end{Definition}

\begin{Definition}\label{definition of the stack of reducible shtukas of type II}
Given $d\ge 2$, $1\le i\le d-1$ and $\mathscr{E}\in\mathfrak{Vect}_{X,D}^i$, we denote  $\RedSht_{\eta,\mathscr{E},D}^{d,i,\II}$ to be the moduli stack which to each scheme $S$ over $\eta$ associates the groupoid of exact sequences
\[\begin{tikzcd}
0 \arrow[r]& \mathscr{A} \arrow[r]& \mathscr{F} \arrow[r]& \mathscr{B} \arrow[r]& 0,
\end{tikzcd}\]
where

(i) $\mathscr{A}\in\TrSht_{\eta,\mathscr{E},D}(S)$, $\mathscr{F}\in\Sht_{\eta,D}^{d,\chi=0}(S)$, $\mathscr{B}\in\Sht_{\eta,D}^{d-i}(S)$;

(ii) the morphisms $\mathscr{A}\to\mathscr{F}$ and $\mathscr{F}\to\mathscr{B}$ are morphisms of shtukas with structures of level $D$;

(iii) the structures of level $D$ give the following commutative diagram
\[\begin{tikzcd}
0 \arrow[r]& \mathscr{A}\otimes\mathscr{O}_{D\times S} \arrow[r]\arrow[d,"\sim"]&\mathscr{F}\otimes\mathscr{O}_{D\times S} \arrow[r]\arrow[d,"\sim"]& \mathscr{B}\otimes\mathscr{O}_{D\times S} \arrow[r] \arrow[d,"\sim"] & 0 \\
0 \arrow[r]& \mathscr{O}_{D\times S}^i \arrow[r]& \mathscr{O}_{D\times S}^d=\mathscr{O}_{D\times S}^i\oplus\mathscr{O}_{D\times S}^{d-i} \arrow[r]& \mathscr{O}_{D\times S}^{d-i} \arrow[r]& 0
\end{tikzcd}\]
where the lower exact sequence is the standard one.
\end{Definition}

\begin{Definition}
  Given $d\ge 2$, $1\le i\le d-1$, and $\mathscr{E}\in\mathfrak{Vect}_{X,\all}^i$, we define
\[\RedSht_{\eta,\mathscr{E},\all}^{d,i,\I}:=\varprojlim_{D}\RedSht_{\eta,\mathscr{E},D}^{d,i,\I},\]
\[\RedSht_{\eta,\mathscr{E},\all}^{d,i,\II}:=\varprojlim_{D}\RedSht_{\eta,\mathscr{E},D}^{d,i,\II},\]
where $D$ runs through all finite subschemes of $X$.
\end{Definition}

\begin{Definition}
Given $d\ge 2, i,j\ge 1, i+j\le d$ and $\mathscr{A}\in\mathfrak{Vect}_{X,D}^i$, $\mathscr{B}\in\mathfrak{Vect}_{X,D}^j$, we define $\RedSht_{\eta,\mathscr{A},\mathscr{B},D}^{d,i,j,\I\wedge\II}$ to be the moduli stack which to each scheme $S$ over $\eta$ associates the groupoid of the following data:

(i) an exact sequence of shtukas with structures of level $D$
\[\begin{tikzcd}
0 \arrow[r]& \mathscr{A}' \arrow[r]& \mathscr{F} \arrow[r]& \mathscr{N} \arrow[r]& 0,
\end{tikzcd}\]
where $\mathscr{A}'\in\TrSht_{\eta,\mathscr{A},D}(S), \mathscr{F}\in\Sht_{\eta,D}^{d,\chi=0}(S), \mathscr{N}\in\Sht_{\eta,D}^{d-i}(S)$, such that the structures of level $D$ induce the standard exact sequence
\[\begin{tikzcd}
0 \arrow[r]& \mathscr{O}_{D\times S}^i \arrow[r]& \mathscr{O}_{D\times S}^d \arrow[r]& \mathscr{O}_{D\times S}^{d-i} \arrow[r]& 0.
\end{tikzcd}\]

(ii) an exact sequence of shtukas with structures of level $D$
\[\begin{tikzcd}
0 \arrow[r]& \mathscr{M} \arrow[r]& \mathscr{N} \arrow[r]& \mathscr{B}' \arrow[r]& 0,
\end{tikzcd}\]
where $\mathscr{M}\in\Sht_{\eta,D}^{d-i-j}(S), \mathscr{B}'\in\TrSht_{\eta,\mathscr{B},D}(S)$, $\mathscr{N}$ is as in (i), such that the structures of level $D$ induce the standard exact sequence
\[\begin{tikzcd}
0 \arrow[r]& \mathscr{O}_{D\times S}^{d-i-j} \arrow[r]& \mathscr{O}_{D\times S}^{d-i} \arrow[r]& \mathscr{O}_{D\times S}^j \arrow[r]& 0.
\end{tikzcd}\]

The above data (i), (ii) are equivalent to the following data:

(i') an exact sequence of shtukas with structures of level $D$
\[\begin{tikzcd}
0 \arrow[r]& \mathscr{L} \arrow[r]& \mathscr{F} \arrow[r]& \mathscr{B}' \arrow[r]& 0,
\end{tikzcd}\]
where $\mathscr{L}\in\Sht_{\eta,D}^{d-j}(S), \mathscr{F}\in\Sht_{\eta,D}^{d,\chi=0}(S), \mathscr{B}'\in\TrSht_{\eta,\mathscr{B},D}(S)$, such that the structures of level $D$ induce the standard exact sequence
\[\begin{tikzcd}
0 \arrow[r]& \mathscr{O}_{D\times S}^{d-j} \arrow[r]& \mathscr{O}_{D\times S}^d \arrow[r]& \mathscr{O}_{D\times S}^j \arrow[r]& 0.
\end{tikzcd}\]

(ii') an exact sequence of shtukas with structures of level $D$
\[\begin{tikzcd}
0 \arrow[r]& \mathscr{A}' \arrow[r]& \mathscr{L} \arrow[r]& \mathscr{M} \arrow[r]& 0,
\end{tikzcd}\]
where $\mathscr{A}'\in\TrSht_{\eta,\mathscr{A},D}(S), \mathscr{M}\in\Sht_{\eta,D}^{d-i-j}(S)$, $\mathscr{L}$ is as in (i'), such that the structures of level $D$ induce the standard exact sequence
\[\begin{tikzcd}
0 \arrow[r]& \mathscr{O}_{D\times S}^i \arrow[r]& \mathscr{O}_{D\times S}^{d-j} \arrow[r]& \mathscr{O}_{D\times S}^{d-i-j} \arrow[r]& 0.
\end{tikzcd}\]
\end{Definition}

\begin{Remark}
We have two Cartesian diagrams
\begin{equation}\label{Cartesian diagrams for reducible shtukas of type I&II}
\begin{tikzcd}
  \RedSht_{\eta,\mathscr{A},\mathscr{B},D}^{d,i,j,\I\wedge\II}, \arrow[r]\arrow[d]& \RedSht_{\eta,\mathscr{B},D}^{d-i,j,\I} \arrow[d]\\
  \RedSht_{\eta,\mathscr{A},D}^{d,i,\II} \arrow[r]& \Sht_{\eta,D}^{d-i}
\end{tikzcd}
\qquad
\begin{tikzcd}
  \RedSht_{\eta,\mathscr{A},\mathscr{B},D}^{d,i,j,\I\wedge\II}, \arrow[r]\arrow[d]& \RedSht_{\eta,\mathscr{A},D}^{d-j,i,\II} \arrow[d]\\
  \RedSht_{\eta,\mathscr{B},D}^{d,j,\I} \arrow[r]& \Sht_{\eta,D}^{d-j}
\end{tikzcd}
\end{equation}
\end{Remark}

\begin{Definition}
  For $\mathscr{A}\in\mathfrak{Vect}_{X,\all}^i, \mathscr{B}\in\mathfrak{Vect}_{X,\all}^j$, we define
  \[\RedSht_{\eta,\mathscr{A},\mathscr{B},\all}^{d,i,j,\I\wedge\II}:=
  \varprojlim_{D}\RedSht_{\eta,\mathscr{A},\mathscr{B},D}^{d,i,j,\I\wedge\II},\]
  where $D$ runs through all finite subschemes of $X$.
\end{Definition}
\subsection{Basic properties of horospherical cycles}
The following statement follows from Corollary 6 of Section 3 of Chapter I of~\cite{L.Lafforgue97}.
\begin{Proposition}
  Let $D$ be a finite subscheme of $X$. Then $\Sht_{\eta,D}^d$ is a Deligne-Mumford stack and it is separated over $\eta$.
\end{Proposition}

The following statement follows from Theorem 9 of Section 2 of Chapter I of~\cite{L.Lafforgue97}.
\begin{Proposition}\label{dimension of the stack of shtukas}
  Let $D$ be a finite subscheme of $X$. The natural morphism $\Sht_{\eta,D}^d\to\eta$ is smooth of pure relative dimension $(2d-2)$. \qed
\end{Proposition}

The following statement is Proposition 5 of Section 3 of Chapter I of~\cite{L.Lafforgue97}.
\begin{Proposition}\label{properties of transition maps between stacks of shtukas}
  For two finite subschemes $D_1\subset D_2\subset X$, the natural morphism $\Sht_{\eta,D_2}^d\to\Sht_{\eta,D_1}^d$ is representable, finite, \'etale and Galois.
\end{Proposition}

The following statement is a consequence of Proposition 2.16(a) of~\cite{Var}.
\begin{Proposition}\label{quasi-compact substack of shtukas representable by a scheme}
  For every finite subscheme $D_1\subset X$ and every quasi-compact open substack $U\subset\Sht_{\eta,D_1}^d$, there exists an integer $N$ such that $U\times_{\Sht_{\eta,D_1}^d}\Sht_{\eta,D_2}^d$ is a scheme for all finite subschemes $D_2\subset X$ satisfying $D_2\ge D_1$ and $\deg D_2\ge N$.
\end{Proposition}

\cref{dimension of the stack of shtukas,properties of transition maps between stacks of shtukas} imply the following statement.
\begin{Proposition}\label{dimension of the schemes of shtukas with structures of all levels}
  The scheme $\Sht_{\eta,\all}^d$ has pure dimension $(2d-2)$. \qed
\end{Proposition}

The following statement follows from Corollary 10 of Section 1 of Chapter II of~\cite{L.Lafforgue97}.
\begin{Proposition}
  Let $D$ be a finite subscheme of $X$ and let $\mathscr{E}\in\mathfrak{Vect}_{X,D}^i$. Then $\RedSht_{\eta,\mathscr{E},D}^{d,i,\I}$ (resp. $\RedSht_{\eta,\mathscr{E},D}^{d,i,\II}$) is a Deligne-Mumford stack and it is separated and locally of finite type over $\eta$.
\end{Proposition}

The following statement follows from Theorem 5 of Section 1 of Chapter II of~\cite{L.Lafforgue97}.
\begin{Proposition}\label{properties of the morphism from RedSht to Sht}
  Let $D$ be a finite subscheme of $X$ and let $\mathscr{E}\in\mathfrak{Vect}_{X,D}^i$. Then the natural morphism
  \[\RedSht_{\eta,\mathscr{E},D}^{d,i,\I}\to\Sht_{\eta,D}^d\]
  \[(\text{resp. $\RedSht_{\eta,\mathscr{E},D}^{d,i,\II}\to\Sht_{\eta,D}^d$}) \]
  is representable, quasi-finite, G-unramified{}\footnote{We use the definition from Stack Project. A morphism is unramified (resp. G-unramified) if and only if it is locally of finite type (resp. locally of finite presentation) and formally unramified.} and separated. \qed
\end{Proposition}

The following statement is Theorem 11 of Section 1 of Chapter II of~\cite{L.Lafforgue97}.
\begin{Proposition}\label{properties of the morphism from RedSht to sub and quotient}
  Let $D$ be a finite subscheme of $X$ and let $\mathscr{E}\in\mathfrak{Vect}_{X,D}^i$. Then the natural morphism
  \[\RedSht_{\eta,\mathscr{E},D}^{d,i,\I}\to\Sht_{\eta,D}^{d-i}\times\TrSht_{\eta,\mathscr{E},D}\]
  \[(\text{resp. $\RedSht_{\eta,\mathscr{E},D}^{d,i,\II}\to\TrSht_{\eta,\mathscr{E},D}\times\Sht_{\eta,D}^{d-i}$}) \]
  is of finite type and smooth of pure relative dimension $i$. \qed
\end{Proposition}

The following statement follows from \cref{dimension of the stack of shtukas,properties of the morphism from RedSht to sub and quotient}.
\begin{Proposition}\label{dimension of the stack of reducible shtukas}
  Let $D$ be a finite subscheme of $X$ and let $\mathscr{E}\in\mathfrak{Vect}_{X,D}^i$. Then the natural morphisms $\RedSht_{\eta,\mathscr{E},D}^{d,i,\I}\to\eta$ and $\RedSht_{\eta,\mathscr{E},D}^{d,i,\II}\to\eta$ are smooth of pure dimension $(2d-i)$. \qed
\end{Proposition}

The following statement is Proposition 4 of Section 1 of Chapter II of~\cite{L.Lafforgue97}.
\begin{Proposition}\label{properties of transition maps between stacks of reducible shtukas}
  Let $D_1\subset D_2$ be two finite subschemes of $X$. Suppose $\mathscr{E}_2\in\mathfrak{Vect}_{X,D_2}^i$ and let $\mathscr{E}_1\in\mathfrak{Vect}_{X,D_1}^i$ be the image of $\mathscr{E}_2$ under the natural map $\mathfrak{Vect}_{X,D_2}\to\mathfrak{Vect}_{X,D_1}$. Then the natural morphism $\RedSht_{\eta,\mathscr{E}_2,D_2}^{d,i,\I}\to\RedSht_{\eta,\mathscr{E}_1,D_1}^{d,i,\I}$ (resp. $\RedSht_{\eta,\mathscr{E}_2,D_2}^{d,i,\II}\to\RedSht_{\eta,\mathscr{E}_1,D_1}^{d,i,\II}$) is representable, finite, \'etale and Galois. \qed
\end{Proposition}

\cref{properties of transition maps between stacks of reducible shtukas,dimension of the stack of reducible shtukas} imply the following statement.
\begin{Proposition}\label{dimension of the scheme of reducible shtukas with structures of all levels}
  For $\mathscr{E}\in\mathfrak{Vect}_{X,\all}^i$, the schemes $\RedSht_{\eta,\mathscr{E},\all}^{d,i,\I}$ and $\RedSht_{\eta,\mathscr{E},\all}^{d,i,\II}$ are reduced and of pure dimension $(2d-i)$. \qed
\end{Proposition}

\cref{properties of the morphism from RedSht to Sht,dimension of the stack of reducible shtukas} imply the following statement.

\begin{Proposition}\label{dimension of horocycles in the stack of shtukas with finite level}
  Let $D$ be a finite subscheme of $X$ and let $\mathscr{E}\in\mathfrak{Vect}_{X,D}^i$. Then the closure of the image of the morphism $\RedSht_{\eta,\mathscr{E},D}^{d,i,\I}\to\Sht_{\eta,D}^d$ (resp. $\RedSht_{\eta,\mathscr{E},D}^{d,i,\I}\to\Sht_{\eta,D}^d$) is reduced and has pure dimension $(2d-i-2)$. \qed
\end{Proposition}

\begin{Lemma}\label{projective limit of closrue of image}
  Let $I$ be a directed set. Let $(X_i)_{i\in I},(Y_i)_{i\in I}$ be two projective systems of schemes with affine surjective transition maps. Let $(X_i\to Y_i)_{i\in I}$ be morphisms compatible with transition maps. Denote $Z_i$ to be the closure of the image of the morphism $X_i\to Y_i$. Then $\varprojlim_{i}Z_i$ equals the closure of the image of the morphism $\varprojlim_{i}X_i\to\varprojlim_{i}Y_i$. \qed
\end{Lemma}

\begin{Proposition}\label{dimension of horocycles in the scheme of shtukas with structures of all levels}
  For $\mathscr{E}\in\mathfrak{Vect}_{X,\all}^i$, the closure of the image of the morphism $\RedSht_{\eta,\mathscr{E},\all}^{d,i,\I}\to\Sht_{\eta,\all}^d$ (resp. $\RedSht_{\eta,\mathscr{E},\all}^{d,i,\II}\to\Sht_{\eta,\all}^d$) is reduced and has pure dimension $(2d-i-2)$. \qed
\end{Proposition}
\begin{proof}
  Let $\mathcal{Y}_{\all}$ denote the closure of the image of the morphism $\RedSht_{\eta,\mathscr{E},\all}^{d,i,\I}\to\Sht_{\eta,\all}^d$. Pick an irreducible component $\mathcal{W}_{\all}$ of $\mathcal{Y}_{\all}$.

  For a finite subscheme $D\subset X$, let $\mathcal{Y}_D$ denote the closure of the image of the morphism $\RedSht_{\eta,\mathscr{E},D}^{d,i,\I}\to\Sht_{\eta,D}^d$, and let $S_D$ denote the set of irreducible components of $\mathcal{Y}_D$ that contain the image of $\mathcal{W}_{\all}$. We see that the transition map $S_{D_2}\to S_{D_1}$ is surjective for any $D_1\subset D_2$. Since the set of finite subschemes of $X$ is countable, $\varprojlim_{D}S_D$ is nonempty. Thus we can find an irreducible component $\mathcal{W}_D\subset\mathcal{Y}_D$ for each $D$ such that $\mathcal{W}_{\all}\subset\varprojlim_D\mathcal{W}_D$. Since each $\mathcal{W}_D$ is irreducible, $\varprojlim_D\mathcal{W}_D$ is also irreducible. \cref{quasi-compact substack of shtukas representable by a scheme,projective limit of closrue of image} show that $\varprojlim_{D}\mathcal{Y}_D=\mathcal{Y}_{\all}$. So $\varprojlim_D\mathcal{W}_D$ is an irreducible component of $\mathcal{Y}_{\all}$. Hence we have $\mathcal{W}_{\all}=\varprojlim_D\mathcal{W}_D$.

  For two finite subschemes $D_1\subset D_2\subset X$, the transition map $\mathcal{W}_{D_2}\to\mathcal{W}_{D_1}$ is finite by \cref{properties of transition maps between stacks of shtukas}. Both $\mathcal{W}_{D_2}$ and $\mathcal{W}_{D_1}$ have dimension $(2d-i-2)$ by \cref{dimension of horocycles in the stack of shtukas with finite level}. Thus the transition map $\mathcal{W}_{D_2}\to\mathcal{W}_{D_1}$ is finite and surjective. Then we see that the morphism $\mathcal{W}_{\all}\to\mathcal{W}_{\emptyset}$ is integral and surjective. Therefore, $\dim\mathcal{W}_{\all}=\dim\mathcal{W}_{\emptyset}=2d-i-2$. This shows that $\mathcal{Y}_{\all}$ has pure dimension $(2d-i-2)$.

  Reducedness of $\mathcal{Y}_{\all}$ follows from reducedness of $\RedSht_{\eta,\mathscr{E},\all}^{d,i,\I}$ by \cref{dimension of the scheme of reducible shtukas with structures of all levels}.

  The statement for the morphism $\RedSht_{\eta,\mathscr{E},\all}^{d,i,\II}\to\Sht_{\eta,\all}^d$ follows from duality.
\end{proof}

The following statement follows from Cartesian diagrams (\ref{Cartesian diagrams for reducible shtukas of type I&II}) and \cref{dimension of the stack of shtukas,dimension of the stack of reducible shtukas,properties of the morphism from RedSht to sub and quotient}.
\begin{Proposition}\label{dimension of the stack of reducible shtukas of type I&II}
  For $\mathscr{A}\in\mathfrak{Vect}_{X,D}^i,\mathscr{B}\in\mathfrak{Vect}_{X,D}^j$, the morphism $\RedSht_{\eta,\mathscr{E}_1,\mathscr{E}_2,D}^{d,i,j,\I\wedge\II}\to\eta$ is smooth of pure dimension $(2d-i-j-2)$. \qed
\end{Proposition}

The following statement follows from Cartesian diagrams (\ref{Cartesian diagrams for reducible shtukas of type I&II}) and \cref{properties of transition maps between stacks of reducible shtukas,properties of transition maps between stacks of shtukas}.

\begin{Proposition}\label{properties of transition maps between stacks of reducible shtukas of type I&II}
  For two finite subschemes $D_1\subset D_2\subset X$ and $\mathscr{A}\in\mathfrak{Vect}_{X,D_2}^i$, $\mathscr{B}\in\mathfrak{Vect}_{X,D_2}^j$, the natural morphism $\RedSht_{\eta,\mathscr{A},\mathscr{B},D_2}^{d,i,j,\I\wedge\II} \to \RedSht_{\eta,\mathscr{A},\mathscr{B},D_1}^{d,i,j,\I\wedge\II}$ is representable, finite, \'etale and Galois.
\end{Proposition}

\cref{dimension of the stack of reducible shtukas of type I&II,properties of transition maps between stacks of reducible shtukas of type I&II} imply the following statement.

\begin{Proposition}\label{dimension of the scheme of reducible shtukas of type I&II with structures of all levels}
  For $\mathscr{A}\in\mathfrak{Vect}_{X,\all}^i$ and $\mathscr{B}\in\mathfrak{Vect}_{X,\all}^j$, the scheme $\RedSht_{\eta,\mathscr{A},\mathscr{B},\all}^{d,i,j,\I\wedge\II}$ has pure dimension $(2d-i-j-2)$. \qed
\end{Proposition}

\subsection{Irreducibility of the scheme of reducible shtukas}
\begin{Theorem}\label{irreducibility of the stack of reducible shtukas with structures of finite level}
  Given $d\ge 2$, a finite subscheme $D\subset X$ and $\mathscr{E}\in\mathfrak{Vect}_{X,D}^1$, the stacks $\RedSht_{\eta,\mathscr{E},D}^{d,1,\I}$ and $\RedSht_{\eta,\mathscr{E},D}^{d,1,\II}$ are irreducible.
\end{Theorem}

\begin{Theorem}\label{irreducibility of the scheme of reducible shtukas with structures of all levels}
  For $d\ge 2$ and $\mathscr{E}\in\mathfrak{Vect}_{X,\all}^1$, the schemes $\RedSht_{\eta,\mathscr{E},\all}^{d,1,\I}$ and $\RedSht_{\eta,\mathscr{E},\all}^{d,1,\II}$ are irreducible.
\end{Theorem}

\begin{Proposition}
  Let $\mathscr{E}\in\mathfrak{Vect}_{X,\all}^1$. Denote $\mathcal{Y}_{\all}^{\I}$ (resp. $\mathcal{Y}_{\all}^{\II}$) to be the closure of the image of the morphism $\RedSht_{\eta,\mathscr{E},\all}^{d,1,\I}\to\Sht_{\eta,\all}^d$ (resp. $\RedSht_{\eta,\mathscr{E},\all}^{d,1,\II}\to\Sht_{\eta,\all}^d$). Then the local ring of $\mathcal{Y}_{\all}^{\I}$ (resp. $\mathcal{Y}_{\all}^{\II}$) in $\Sht_{\eta,\all}^d$ is a discrete valuation ring.
\end{Proposition}
\begin{proof}
  For a finite subscheme $D\subset X$, denote $\mathcal{Y}_{D}^{\I}$ to be the closure of the image of the morphism $\RedSht_{\eta,\mathscr{E},D}^{d,1,\I}\to\Sht_{\eta,D}^d$. \cref{irreducibility of the stack of reducible shtukas with structures of finite level} implies that $\mathcal{Y}_{D}^{\I}$ is irreducible. Choose a quasi-compact substack $U_{\emptyset}\subset\Sht_{\eta,\emptyset}^d$ such that $U_{\emptyset}\cap\mathcal{Y}_{\emptyset}^{\I}$ is dense in $\mathcal{Y}_{\emptyset}^{\I}$. Denote $U_D=U_{\emptyset}\times_{\Sht_{\eta,\emptyset}^d}\Sht_{\eta,D}^d$. \cref{properties of transition maps between stacks of reducible shtukas} shows that the morphism $\mathcal{Y}_D^{\I}\to\mathcal{Y}_{\emptyset}^{\I}$ is dominant. Thus $U_D\cap\mathcal{Y}_D^{\I}$ is dense in $\mathcal{Y}_D^{\I}$ for all $D$. By \cref{quasi-compact substack of shtukas representable by a scheme}, there exists an integer $N$ such that $U_D$ is a scheme when $\deg D\ge N$. Let $A_D$ be the local ring of $U_D\cap\mathcal{Y}_D^{\I}$ in $U_D$ for those $D$ satisfying $\deg D\ge N$. Let $A_{\all}$ be the local ring of $\mathcal{Y}_{\all}^{\I}$ in $\Sht_{\eta,\all}^d$. Then \cref{projective limit of closrue of image} implies that $A_{\all}=\varinjlim_{\deg D\ge N}A_D$. For each $D$,  $\Sht_{\eta,D}^d$ is smooth over $\eta$ by \cref{dimension of the stack of shtukas}, and $\mathcal{Y}_{D}^{\I}$ has codimension 1 in $\Sht_{\eta,D}^d$ by \cref{dimension of horocycles in the stack of shtukas with finite level,dimension of the stack of shtukas}. Hence $A_D$ is a discrete valuation. For two finite subschemes $D_1\subset D_2\subset X$, the transition map $U_{D_2}\to U_{D_1}$ is smooth by \cref{properties of transition maps between stacks of shtukas}. Hence the transition homomorphism $A_{D_1}\to A_{D_2}$ sends uniformizer to uniformizer.  Now \cref{inductive limit of discrete valuation rings} implies that $A_{\all}$ is a discrete valuation ring.

  The statement for $\mathcal{Y}_{\all}^{\II}$ follows from duality.
\end{proof}

\subsection{Criteria for a Tate toy horospherical subscheme to contain the image of a horospherical divisor}
Fix an integer $d\ge 2$. For a shtuka $\mathscr{G}$ over a perfect field with zero and pole satisfying condition (\ref{zero and pole not related by Frobenius}), the notation $\mathscr{G}^{\I}$ and $\mathscr{G}^{\II}$ is defined in \cref{maximal trivial sub and maximal trivial quotient}.

\begin{Definition}\label{notation of horospherical divisors}
  Let $\mathcal{Z}_{\eta,1}^{d,\I}$ (resp. $\mathcal{Z}_{\eta,1}^{d,\II}$) denote the closure of the image of the morphism $\RedSht_{\eta,\mathscr{O}_X,\all}^{d,1,\I}\to\Sht_{\eta,\all}^{d,\chi=0}$ (resp. $\RedSht_{\eta,\mathscr{O}_X,\all}^{d,1,\II}\to\Sht_{\eta,\all}^{d,\chi=0}$), where $\mathscr{O}_X\in\mathfrak{Vect}_{X,\all}^1$ is equipped with the standard structures of all levels.
\end{Definition}

The following statement follows from \cref{irreducibility of the scheme of reducible shtukas with structures of all levels,dimension of the scheme of reducible shtukas with structures of all levels}.
\begin{Theorem}\label{horospherical divisors are reduced and irreducible}
  $\mathcal{Z}_{\eta,1}^{d,\I}$ and $\mathcal{Z}_{\eta,1}^{d,\II}$ are reduced and irreducible.
\end{Theorem}

Let $\xi_{\eta,1}^{d,\I}$ (resp. $\xi_{\eta,1}^{d,\II}$) denote the generic point of $\mathcal{Z}_{\eta,1}^{d,\I}$ (resp. $\mathcal{Z}_{\eta,1}^{d,\II}$).

\begin{Lemma}\label{generic reducible shtuka of type I does not contain trivial sub}
  Let $\xi'$ be a geometric generic point of $\mathcal{Z}_{\eta,1}^{d,\I}$. Let $\mathscr{F}$ be the shtuka over $\xi'$. Then $\mathscr{F}^{\II}=0$.
\end{Lemma}
\begin{proof}
  Denote $\mathscr{A}'=\mathscr{F}^{\II}$ and $i=\rank\mathscr{A}'$. By \cref{descent of coherent sheaves for projective scheme over finite field}, we can find an isomorphism $\mathscr{A}'\cong\mathscr{A}\otimes\xi'$ for some $\mathscr{A}\in\mathfrak{Vect}_{X,\emptyset}^{\rank\mathscr{A}'}$. Suppose $\mathscr{A}\ne 0$.

  For a finite subscheme $D\subset X$, denote $\mathscr{L}'_D$ to be the image of the composition $\mathscr{A}\otimes\mathscr{O}_D\otimes \xi'\to\mathscr{F}\otimes\mathscr{O}_{D\otimes\xi'}$. We see that there exists an $\mathscr{O}_D$-submodule $\mathscr{L}_D\subset\mathscr{O}_D^d$ such that $\mathscr{L}'_D=\mathscr{L}_D\otimes\xi'$. Moreover, for two subschemes $D_1\subset D_2\subset X$, we have $\mathscr{L}_{D_1}=\mathscr{L}_{D_2}\otimes\mathscr{O}_{D_1}$. Thus we can find $g\in GL_d(O)$ such that $g(\mathscr{O}_{D}^i\oplus0)=\mathscr{L}_D$ for all  finite subschemes $D\subset X$. We equip $\mathscr{A}$ with structures of all levels compatible with each other using the standard structures of all levels on $\mathscr{O}_X^i$. Then the image of $g\cdot\xi'$ is contained in the image of the morphism $\RedSht_{\eta,\mathscr{A},\all}^{d,i,\II}\to\Sht_{\eta,\all}^d$. \cref{dimension of horocycles in the scheme of shtukas with structures of all levels} implies that $i=1$ for dimensional reasons.

  From the definition of $\RedSht_{\eta,\mathscr{O}_X,\all}^{d,1,\I}$ we obtain an exact sequence of shtukas
  \[\begin{tikzcd}
  0 \arrow[r]& \mathscr{S} \arrow[r]& \mathscr{F} \arrow[r]& \mathscr{B}' \arrow[r]& 0
  \end{tikzcd}\]
  where $\mathscr{B}'\in\TrSht_{\eta,\mathscr{O}_X,\all}(\xi')$,  $\mathscr{S}\in\Sht_{\eta,\all}^{d-1}(\xi')$

  From the definition of $\mathscr{F}^{\II}$ in \cref{maximal trivial sub and maximal trivial quotient} we know that $\mathscr{A}'$ is saturated in $\mathscr{F}$. If the composition $\mathscr{A}'\to\mathscr{F}\to\mathscr{B}'$ is zero, we have an exact sequence of shtukas
  \[\begin{tikzcd}
  0 \arrow[r]& \mathscr{A}' \arrow[r]& \mathscr{S} \arrow[r]& \mathscr{M} \arrow[r]& 0.
  \end{tikzcd}\]
  Then one can find $h\in GL_d(O)$ such that the image of $h\cdot\xi'$ is contained in the image of the morphism $\RedSht_{\eta,\mathscr{A},\mathscr{O}_X,\all}^{d,1,1,\I\wedge\II}\to\Sht_{\eta,\all}^d$. This is a contradiction to \cref{dimension of horocycles in the scheme of shtukas with structures of all levels,dimension of the scheme of reducible shtukas of type I&II with structures of all levels} for dimensional reasons. Thus the composition $\mathscr{A}'\to\mathscr{F}\to\mathscr{B}'$ is nonzero, and it is injective since $\rank\mathscr{A}'=1$. Therefore, the morphism $\mathscr{A}'\oplus\mathscr{S}\to\mathscr{F}$ is injective.

  \cref{subshtuka of the same rank} implies that $\mathscr{A}'\oplus\mathscr{S}$ is obtained from $\mathscr{F}$ by Construction E. Hence there exists $w\in GL_d(\mathbb{A})$ such that the image of $w\cdot\xi'$ is contained in the image of the morphism $\Sht_{\eta,\all}^{d-1}\to\Sht_{\eta,\all}^d$ which sends a shtuka $\mathscr{S}$ of rank $(d-1)$ to the shtuka $\mathscr{A}'\oplus\mathscr{S}$ of rank $d$. By \cref{dimension of the schemes of shtukas with structures of all levels,dimension of horocycles in the scheme of shtukas with structures of all levels} we get a contradiction for dimensional reasons.
\end{proof}

As in \cref{(Section)basic properties of Tate toy horospherical subschemes}, we denote $\Delta_{\mathbb{A}^d,J}^n=\oToySht_{\mathbb{A}^d}^n\cap\ToySht_{\mathbb{A}^d/J}$ for $J\in\mathbf{P}_{\mathbb{A}^d}$, and we denote $\Delta_{\mathbb{A}^d,H}^n=\oToySht_{\mathbb{A}^d}^n\cap\ToySht_{H}$ for $H\in\mathbf{P}_{(\mathbb{A}^d)^*}$.

Denote $GL_d(\mathbb{A})_0=\{g\in GL_d(\mathbb{A})|\deg g=0\}$.

Recall that for $\chi\in\mathbb{Z}$, we denote $n_\chi$ to be the element of $\Dim_{\mathbb{A}^d}$ such that $n_\chi(O^d)=\chi$. In particular, we have $n_0\in\Dim_{\mathbb{A}^d}$ corresponding to $\chi=0$.

Let $\theta_\eta^d:\Sht_{\eta,\all}^d\to\oToySht_{\mathbb{A}^d}$ be the morphism defined in \cref{construction of the morphism theta}.

\begin{Lemma}\label{image of type I horospherical divisor not contained in any type II toy horospherical divisor}
  For all $g\in GL_d(\mathbb{A})_0$ and $J\in\mathbf{P}_{\mathbb{A}^d}$, we have $\theta_\eta^d(g\cdot\xi_{\eta,1}^{d,\I})\notin\Delta_{\mathbb{A}^d,J}^{n_0}$.
\end{Lemma}
\begin{proof}
  Since $\theta_\eta^d$ is $GL_d(\mathbb{A})_0$-equivariant, it suffices to prove the statement in the case $g=1$. Let $\xi'$ be a geometric generic point of $\mathcal{Z}_{\eta,1}^{d,\I}$. Let $\mathscr{F}$ be the shtuka over $\xi'$. Suppose $\theta_\eta^d(\xi_{\eta,1}^{d,\I})\in\Delta_{\mathbb{A}^d,J}^{n_0}$ for some $J\in\mathbf{P}_{\mathbb{A}^d}$. Then there exists a divisor $D''$ of $X$ and a nonzero element $z\in H^0(X\times \xi',\mathscr{F}(D''\times\xi'))$ such that $\Phi_{\xi'}^*z=z$. Let $\mathscr{G}$ be the subsheaf of $\mathscr{F}(D''\times\xi')$ generated by $z$. We see that $\mathscr{G}\ne 0$ and $\Phi_{\xi'}^*\mathscr{G}=\mathscr{G}$. Thus $(\mathscr{F}(D''\times\xi'))^{\II}\ne 0$. Then \cref{twist of maximal trivial sub and maximal trivial quotient} shows that $\mathscr{F}^{\II}\ne 0$. This is a contradiction to \cref{generic reducible shtuka of type I does not contain trivial sub}.
\end{proof}

\begin{Lemma}\label{generic reducible shtuka of type II has maximal trivial sub of rank 1}
  Let $E$ be a separably closed field over $\eta$. Let $\xi'':\Spec E\to\mathcal{Z}_{\eta,1}^{d,\II}$ be a morphism over $\eta$ whose image lands in the generic point of $\mathcal{Z}_{\eta,1}^{d,\II}$. Let $\mathscr{F}$ be the shtuka over $\xi''$. Then there exists an isomorphism $\mathscr{O}_{X\times\xi''}\xrightarrow{\sim}\mathscr{F}^{\II}$, such that for every finite subscheme $D\subset X$, the composition
  \[\mathscr{O}_{X\times\xi''}\xrightarrow{\sim}\mathscr{F}^{\II}\to\mathscr{F}\to\mathscr{F}\otimes\mathscr{O}_{D\times\xi''} \xrightarrow{\sim}\mathscr{O}_{D\times\xi''}^d\]
  is the following standard morphism
  \[\mathscr{O}_{X\times\xi''}\twoheadrightarrow\mathscr{O}_{D\times\xi''}\hookrightarrow
  \mathscr{O}_{D\times\xi''}\oplus\mathscr{O}_{D\times\xi''}^{d-1}=\mathscr{O}_{X\times\xi''}^d.\]
\end{Lemma}
\begin{proof}
  From the definition of $\RedSht_{\eta,\all}^{d,1,\II}$ and \cref{descent of locally free sheaves with level structure} we see that there exists an injective morphism $\mathscr{O}_{X\times\xi''}\hookrightarrow\mathscr{F}^{\II}$ satisfying the required conditions. Moreover, $\mathscr{O}_{X\times\xi''}$ is saturated in $\mathscr{F}$. Thus it suffices to show that $\rank\mathscr{F}^{\II}=1$.

  Suppose $\rank\mathscr{F}^{\II}\ge 2$. Then we can find $g\in GL_d(O)$ such that the image of $g\cdot\xi''$ is contained in the image of the morphism $\RedSht_{\eta,\mathscr{A},\all}^{d,2,\II}\to\Sht_{\eta,\all}^d$ for some $\mathscr{A}\in\mathfrak{Vect}_{X,\all}^2$. We get a contradiction by \cref{dimension of horocycles in the scheme of shtukas with structures of all levels} for dimensional reasons.
\end{proof}

Recall that $k$ denotes the field of rational functions on $X$.
\begin{Proposition}\label{Tate toy horosphrical divisors containing the image of Z_1^II}
  For $J\in\mathbf{P}_{\mathbb{A}^d}$, $\theta_\eta^d(\mathcal{Z}_{\eta,1}^{d,\II})\subset\Delta_{\mathbb{A}^d,J}^{n_0}$ if and only if
  $J=\mathbb{F}_q\cdot(a,0,\dots,0)^t$ for some $a\in k^\times$.
\end{Proposition}
\begin{proof}
  Recall that $\mathcal{Z}_{\eta,1}^{d,\II}$ is reduced and irreducible by \cref{horospherical divisors are reduced and irreducible}. Let $\xi''$ be a geometric generic point of  $\mathcal{Z}_{\eta,1}^{d,\II}$. Thus $\theta_\eta^d(\mathcal{Z}_{\eta,1}^{d,\II})\subset\Delta_{\mathbb{A}^d,J}^{n_0}$ if and only if the image of $\theta_{\eta}^d\bcirc\xi''$ is contained in $\Delta_{\mathbb{A}^d,J}^{n_0}$.

  Let $\mathscr{F}$ be the shtuka over $\xi''$ and let $L\subset\mathbb{A}^d\otimes\xi''$ be the corresponding Tate toy shtuka. \cref{definition of the stack of reducible shtukas of type II} gives a subshtuka $\mathscr{A}\subset\mathscr{F}$. Here $\mathscr{A}\in \TrSht_{\eta,\mathscr{O}_X,\all}^1(\xi'')$ and $\mathscr{O}_X$ is equipped with the standard structure of all levels. Moreover, the morphism $\mathscr{A}\otimes\mathscr{O}_{D\times\xi''}\to\mathscr{F}\otimes\mathscr{O}_{D\times\xi''}$ induces the standard inclusion of the first entry $\mathscr{O}_{D\times\xi''}\to\mathscr{O}_{D\times\xi''}^d$ for all finite subscheme $D\subset X$. \cref{descent of coherent sheaves for projective scheme over finite field} shows that $\mathscr{A}\cong\mathscr{O}_{X\times\xi''}$. Thus the constant function $1\in H^0(X\times\xi'',\mathscr{O}_{X\times\xi''})$ gives an element $(1,0,\dots,0)^t\in L$. Hence the image of $\theta_{\eta}^d\bcirc\xi''$ is contained in $\Delta_{\mathbb{A}^d,J_1}^{n_0}$, where $J_1=\mathbb{F}_q\cdot(1,0,\dots,0)^t\in\mathbf{P}_{\mathbb{A}^d}$. Since the $k^\times$-action on $\Sht_{\eta,\all}^{d,\chi=0}$ is trivial and the morphism $\theta_\eta^d$ is $k^\times$-equivariant, we deduce that $\theta_{\eta}^d(\mathcal{Z}_{\eta,1}^{d,\II})\subset\Delta_{\mathbb{A}^d,J_a}^{n_0}$ for all $a\in k^\times$, where $J_a=\mathbb{F}_q\cdot(a,0,\dots,0)^t\in\mathbf{P}_{\mathbb{A}^d}$.

  Suppose that the image of $\theta_{\eta}^d\bcirc\xi''$ is contained in $\Delta_{\mathbb{A}^d,J}^{n_0}$ for $J\in\mathbf{P}_{\mathbb{A}^d}$. From the definition of $\theta_\eta^d$ we see that
  \[J\otimes\xi''\subset H^0(X\times\xi'',\mathscr{F}(D''\times\xi''))\]
  for some divisor $D''$ of $X$. Let $\mathscr{F}^{\II}$ be as in \cref{maximal trivial sub and maximal trivial quotient}. Now \cref{twist of maximal trivial sub and maximal trivial quotient} shows that
  \[J\otimes\xi''\subset H^0(X\times\xi'',\mathscr{F}^{\II}(D''\times\xi'')).\]
  Applying \cref{generic reducible shtuka of type II has maximal trivial sub of rank 1}, we deduce that $J=\mathbb{F}_q\cdot(a,0,\dots,0)^t$ for some $a\in k^\times$.
\end{proof}

\subsection{Irreducible components of horospherical divisors}\label{(Section)irreducible components of horospherical divisors}
We introduce two subgroups $P_d^{\I}=\begin{psmallmatrix}GL_{d-1} & * \\ 0 & 1 \end{psmallmatrix}\subset GL_d$ and $P_d^{\II}=\begin{psmallmatrix}1 & * \\ 0 & GL_{d-1} \end{psmallmatrix}\subset GL_d$.

Let $P_d^{\I}(\mathbb{A})_0=\{g\in P_d^{\I}(\mathbb{A})|\deg g=0\}$, $P_d^{\II}(\mathbb{A})_0=\{g\in P_d^{\II}(\mathbb{A})|\deg g=0\}$.

\begin{Lemma}\label{stabilizer of the group action on the irreducible components of horospherical divisors}
  For $g\in GL_d(\mathbb{A})$, $g\cdot\mathcal{Z}_{\eta,1}^{d,\I}=\mathcal{Z}_{\eta,1}^{d,\I}$ if and only if $g\in k^\times\cdot P_d^{\I}(\mathbb{A})_0$, $g\cdot\mathcal{Z}_{\eta,1}^{d,\II}=\mathcal{Z}_{\eta,1}^{d,\II}$ if and only if $g\in k^\times\cdot P_d^{\II}(\mathbb{A})_0$.
\end{Lemma}
\begin{proof}
  Recall that the $k^\times$-action on $\Sht_{\eta,\all}^d$ is trivial.

  Let $g\in P_d^{\II}(\mathbb{A})_0$. Let $\mathscr{O}_X\in\mathfrak{Vect}_{X,\all}^1$ be equipped with the standard structures of all levels. Recall that $\mathcal{Z}_{\eta,1}^{d,\II}$ is irreducible. Let $\xi$ be the generic point of $\mathcal{Z}_{\eta,1}^{d,\II}$. Then $\xi$ is contained in the image of the morphism $\RedSht_{\eta,\mathscr{O}_X,\all}^{d,1,\II}\to\Sht_{\eta,\all}^d$. Let $\mathscr{F}$ be the shtuka over $\xi$. We have an inclusion of shtukas $\mathscr{A}\subset\mathscr{F}$, where $\mathscr{A}\in\TrSht_{\eta,\mathscr{O}_X,\all}$. Moreover, for all finite subscheme $D\subset X$, the induced morphism $\mathscr{A}\otimes\mathscr{O}_{D\times\xi}\to\mathscr{F}\otimes\mathscr{O}_{D\times\xi}$ gives the standard inclusion of the first entry $\mathscr{O}_{D\times\xi}\to\mathscr{O}_{D\times\xi}^d$.  From the construction of the $GL_d(\mathbb{A})$-action on $\Sht_{\eta,\all}^d$, we see that $\mathscr{A}\subset g\cdot\mathscr{F}$, and for all finite subscheme $D\subset X$, the induced morphism $\mathscr{A}\otimes\mathscr{O}_{D\times\xi}\to(g\cdot\mathscr{F})\otimes\mathscr{O}_{D\times\xi}$ gives the standard inclusion of the first entry $\mathscr{O}_{D\times\xi}\to\mathscr{O}_{D\times\xi}^d$. This shows that $g\cdot\xi$ is contained in the image of the morphism $\RedSht_{\eta,\mathscr{O}_X,\all}^{d,1,\II}\to\Sht_{\eta,\all}^d$, hence contained in $\mathcal{Z}_{\eta,1}^{d,\I}$.

  Suppose $g\cdot\mathcal{Z}_{\eta,1}^{d,\II}=\mathcal{Z}_{\eta,1}^{d,\II}$. \cref{Euler characteristic of shtukas under group action} shows that $g\in GL_d(\mathbb{A})_0$. Let $\mathbf{J}_1=\{J\in\mathbf{P}_{\mathbb{A}^d}|J=\mathbb{F}_q\cdot (a,0,\dots,0)^t,a\in k^\times\}$. For $J\in\mathbf{P}_{\mathbb{A}^d}$, \cref{Tate toy horosphrical divisors containing the image of Z_1^II} shows that $\theta_\eta^d(\mathcal{Z}_{\eta,1}^{d,\II})\in\Delta_{\mathbb{A}^d,J}^{n_0}$ if and only if $J\in\mathbf{J}_1$. Now $g\cdot\mathcal{Z}_{\eta,1}^{d,\II}=\mathcal{Z}_{\eta,1}^{d,\II}$ implies that $g\cdot \mathbf{J}_1=\mathbf{J}_1$. Hence $g\in k^\times\cdot P_d^{\II}(\mathbb{A})_0$.

  The second statement follows from duality.
\end{proof}

Denote
\[\ooSht_{\eta,\all}^d:=\Sht_{\eta,\all}^d-\bigcup_{g\in GL_d(\mathbb{A})}g\cdot\mathcal{Z}_{\eta,1}^{d,\I}
-\bigcup_{g\in GL_d(\mathbb{A})}g\cdot\mathcal{Z}_{\eta,1}^{d,\II}.\]

We denote $\Omega_\eta^{d,\I}=GL_d(\mathbb{A})_0/(k^\times\cdot P_d^{\I}(\mathbb{A})_0)$, $\Omega_\eta^{d,\II}=GL_d(\mathbb{A})_0/(k^\times\cdot P_d^{\II}(\mathbb{A})_0)$.

\begin{Definition}
  A Cartier divisor of $\Sht_{\eta,\all}^{d}$ is called a \emph{horospherical divisor} if its support has empty intersection with $\ooSht_{\eta,\all}^d$.
\end{Definition}

The following  statement is a corollary of \cref{stabilizer of the group action on the irreducible components of horospherical divisors}.
\begin{Proposition}\label{irreducible components of horospherical divisors}
  The set of irreducible components of horospherical divisors of $\Sht_{\eta,\all}^{d,\chi=0}$ is $\Omega_{\eta}^{d,\I}\coprod\Omega_{\eta}^{d,\II}$.
\end{Proposition}

\begin{Lemma}
  The $GL_d(\mathbb{A})_0$-action on $\Omega_{\eta}^{d,\I}\coprod\Omega_{\eta}^{d,\II}$ is continuous.
\end{Lemma}
\begin{proof}
  The statement follows the construction of the $GL_d(\mathbb{A})$-action on $\Sht_{\eta,\all}^d$.
\end{proof}

\begin{Lemma}\label{reducible shtukas of type II are contained in horospherical divisors of type II}
  Let $E$ be an algebraically closed field over $\eta$. Suppose $\mathscr{F}\in\Sht_{\eta,\all}^{d}(E)$ satisfies $\mathscr{F}^{\II}\ne 0$. Then there exists $g\in GL_d(\mathbb{A})$ such that $g\cdot\mathscr{F}\in\mathcal{Z}_{\eta,1}^{d,\II}$.
\end{Lemma}
\begin{proof}
  By \cref{descent of coherent sheaves for projective scheme over finite field}, we have an isomorphism $\mathscr{F}^{\II}\cong\mathscr{M}\otimes E$ for some $\mathscr{M}\in\mathfrak{Vect}_{X,\emptyset}^{\rank\mathscr{F}}$. Choose an invertible subsheaf $\mathscr{A}\subset\mathscr{M}$ on $X$ which is saturated in $\mathscr{M}$. Then we have an exact sequence of shtukas
  \[\begin{tikzcd}
  0 \arrow[r]& \mathscr{A}\otimes E \arrow[r]& \mathscr{F} \arrow[r]& \mathscr{Q} \arrow[r]& 0.
  \end{tikzcd}\]

  For a finite subscheme $D\subset X$, Let $\mathscr{L}'_{D}$ be the image of the composition
  $\mathscr{A}\otimes\mathscr{O}_{D\otimes E} \to \mathscr{F}\otimes\mathscr{O}_{D\otimes E}
  \simto \mathscr{O}_{D\otimes E}^d$
  We see that $\mathscr{L}'_D$ is invariant under the natural isomorphism $\Phi_E^*\mathscr{O}_{D\otimes E}^d\simto\mathscr{O}_{D\otimes E}^d$. By \cref{descent of coherent sheaves for projective scheme over finite field}, we have $\mathscr{L}'_D=\mathscr{L}_D\otimes E$ for some $\mathscr{O}_D$-submodule $\mathscr{L}_D\subset\mathscr{O}_{D}^d$. Moreover, for two finite subschemes $D_1\subset D_2\subset X$, we have $\mathscr{L}_{D_1}=\mathscr{L}_{D_2}\otimes\mathscr{O}_{D_1}$. Thus we can choose $g_1\in GL_d(O)$ such that $g_1(\mathscr{O}_D\oplus0)=\mathscr{L}_D$ for all finite subschemes $D\subset X$. We equip $\mathscr{A}$ with structures of all levels compatible with each other using the standard structures of all levels on $\mathscr{O}_X$.

  Since $\mathbb{A}^\times$ acts transitively on $\mathfrak{Vect}_{X,\all}^1$, we can choose $a\in\mathbb{A}^\times$ such that $a\cdot\mathscr{A}=\mathscr{O}_X$, where $\mathscr{O}_X$ is equipped with the standard structures of all levels. Pick $B\in GL_{d-1}(\mathbb{A})$ such that $\deg a+\deg\det B+\chi(\mathscr{F})=0$. Let
  $g_2=\begin{psmallmatrix}a & 0 \\ 0 & B\end{psmallmatrix}$.

  Then we see that $g_2g_1\cdot\mathscr{F}\in\mathcal{Z}_{\eta,1}^{d,\II}$.
\end{proof}

The following statement is dual to \cref{reducible shtukas of type II are contained in horospherical divisors of type II}.
\begin{Lemma}\label{reducible shtukas of type I are contained in horospherical divisors of type I}
  Let $E$ be an algebraically closed field over $\eta$. Suppose $\mathscr{F}\in\Sht_{\eta,\all}^{d}(E)$ satisfies $\mathscr{F}^{\I}\ne 0$. Then there exists $g\in GL_d(\mathbb{A})$ such that $g\cdot\mathscr{F}\in\mathcal{Z}_{\eta,1}^{d,\I}$. \qed
\end{Lemma}

\begin{Lemma}\label{reducible shtukas are contained in horospherical divisors}
  Let $s$ be a geometric point of $\Sht_{\eta,\all}^d$ such that $s\notin g\cdot\mathcal{Z}_{\eta,1}^{d,\I},s\notin g\cdot\mathcal{Z}_{\eta,1}^{d,\II}$ for all $g\in GL_d(\mathbb{A})$. Then the shtuka over $s$ is irreducible.
\end{Lemma}
\begin{proof}
  The statement follows from \cref{reducible shtukas of type II are contained in horospherical divisors of type II,reducible shtukas of type I are contained in horospherical divisors of type I}.
\end{proof}

\begin{Lemma}\label{image of irreducible shtukas under theta}
  For any point $s\in\ooSht_{\eta,\all}^{d,\chi}$, we have $\theta_{\eta}^d(s)\in\ooToySht_{\mathbb{A}^d}^{n_\chi}$.
\end{Lemma}
\begin{proof}
\cref{construction of the morphism theta} shows that $\theta_{\eta}^d(s)\in\oToySht_{\mathbb{A}^d}^{n_\chi}$.

We may assume that $s$ is a geometric point of $\ooSht_{\eta,\all}^{d,\chi}$. Let $\mathscr{F}$ be the shtuka over $s$. Then $\mathscr{F}$ is irreducible by \cref{reducible shtukas are contained in horospherical divisors}. Hence $\mathscr{F}^{\II}=0$ by \cref{maximal trivial sub and maximal trivial quotient}. Now \cref{twist of maximal trivial sub and maximal trivial quotient} shows that $\mathscr{F}(D''\times s)^{\II}=0$ for any divisor $D''$ of $X$. Thus $\theta_{\eta}^d(s)\notin\Delta_{\mathbb{A}^d,J}^{n_\chi}$ for any $J\in\mathbf{P}_{\mathbb{A}^d}$. By duality, $\theta_{\eta}^d(s)\notin\Delta_{\mathbb{A}^d,H}^{n_\chi}$ for any $H\in\mathbf{P}_{(\mathbb{A}^d)^*}$. The statement follows.
\end{proof}

\section{Pullback of Tate toy horospherical divisors under $\theta$}\label{(Section)pullback of toy horospherical divisors under theta}
We use the notation and conventions of \cref{(Section)notation and conventions for shtukas}.

Fix an integer $d\ge 2$.

Let $\theta_\eta^d:\Sht_{\eta,\all}^d\to\oToySht_{\mathbb{A}^d}$ be the morphism defined in \cref{construction of the morphism theta}.

For $\chi\in\mathbb{Z}$, we denote $n_\chi$ to be the element of $\Dim_{\mathbb{A}^d}$ such that $n_\chi(O^d)=\chi$. In particular, we have $n_0\in\Dim_{\mathbb{A}^d}$ corresponding to $\chi=0$.

Denote $\theta_{\eta}^{d,0}:\Sht_{\eta,\all}^{d,\chi=0}\to\oToySht_{\mathbb{A}^d}^{n_0}$.
\subsection{$(\theta_{\eta}^{d,0})^*$ as a direct sum}\label{(Section)theta as a direct sum}
Recall that a Cartier divisor of $\oToySht_{\mathbb{A}^d}$ is called a \emph{(Tate) toy horospherical divisor} if its restriction to $\ooToySht_{\mathbb{A}^d}$ is zero. \cref{image of irreducible shtukas under theta} shows that it makes sense to pullback Tate toy horospherical divisors under $\theta_{\eta}^d$, and the pullback is a horospherical divisor of $\Sht_{\eta,\all}^d$.

Recall that $\mathfrak{O}_{\mathbb{A}^d}^{n_0}$ denotes the space of Tate toy horospherical divisors on $\oToySht_{\mathbb{A}^d}^{n_0}$. \cref{space of Tate toy horospherical subschemes} gives an isomorphism
\[\mathfrak{O}_{\mathbb{A}^d}^{n_0} \cong C_+((\mathbb{A}^d)^*-\{0\})^{\mathbb{F}_q^\times}\oplus C_+(\mathbb{A}^d-\{0\})^{\mathbb{F}_q^\times}.\]

\cref{irreducible components of horospherical divisors} shows that the set of irreducible components of horospherical divisors of $\Sht_{\eta,\all}^{d,\chi=0}$ is $\Omega_{\eta}^{d,\I}\coprod\Omega_{\eta}^{d,\II}$

Considering the multiplicity of pullback of Tate toy horospherical divisors under the morphism $\theta_{\eta}^{d,0}:\Sht_{\eta,\all}^{d,\chi=0}\to\oToySht_{\mathbb{A}^d}^{n_0}$, we get a homomorphism of (partially) ordered abelian groups
\[(\theta_\eta^{d,0})^*:
C_+((\mathbb{A}^d)^*-\{0\})^{\mathbb{F}_q^\times}\oplus C_+(\mathbb{A}^d-\{0\})^{\mathbb{F}_q^\times} \to
C^\infty(\Omega_\eta^{d,\I})\oplus C^\infty(\Omega_\eta^{d,\II}).\]

\begin{Lemma}\label{no interference between I&II for pullback of toy horospherical divisors}
  The homomorphisms
  \[C_+((\mathbb{A}^d)^*-\{0\})^{\mathbb{F}_q^\times}\to C^\infty(\Omega_{\eta}^{d,\II}),\]
  \[C_+(\mathbb{A}^d-\{0\})^{\mathbb{F}_q^\times}\to C^\infty(\Omega_{\eta}^{d,\I})\]
  induced by $(\theta_{\eta}^{d,0})^*$ are zero.
\end{Lemma}
\begin{proof}
  The first homomorphism is zero by \cref{image of type I horospherical divisor not contained in any type II toy horospherical divisor}. By duality, the second one is also zero.
\end{proof}

The homomorphism $(\theta_{\eta}^{d,0})^*$ induces two maps
\[(\theta_{\eta}^{d,0,\I})^*:C_+((\mathbb{A}^d)^*-\{0\})^{\mathbb{F}_q^\times}\to C^\infty(\Omega_{\eta}^{d,\I}),\]
\[(\theta_{\eta}^{d,0,\II})^*:C_+(\mathbb{A}^d-\{0\})^{\mathbb{F}_q^\times}\to C^\infty(\Omega_{\eta}^{d,\II}).\]

The following statement is a corollary of \cref{no interference between I&II for pullback of toy horospherical divisors}.
\begin{Lemma}\label{theta as a direct sum}
  $(\theta_{\eta}^{d,0})^*=(\theta_{\eta}^{d,0,\I})^*\oplus(\theta_{\eta}^{d,0,\II})^*$. \qed
\end{Lemma}

\subsection{Definition of the averaging maps}
Recall that $k$ denotes the field of rational functions on $X$. We introduced two subgroups $P_d^{\I}=\begin{psmallmatrix}GL_{d-1} & * \\ 0 & 1 \end{psmallmatrix}\subset GL_d$ and $P_d^{\II}=\begin{psmallmatrix}1 & * \\ 0 & GL_{d-1} \end{psmallmatrix}\subset GL_d$ in \cref{(Section)irreducible components of horospherical divisors}. We also denoted $\Omega_\eta^{d,\I}=GL_d(\mathbb{A})_0/(k^\times\cdot P_d^{\I}(\mathbb{A})_0)$, $\Omega_\eta^{d,\II}=GL_d(\mathbb{A})_0/(k^\times\cdot P_d^{\II}(\mathbb{A})_0)$ .

We introduce two varieties over $\mathbb{F}_q$.
\[Y_d^{\I}:= GL_d/P_d^{\I},\qquad Y_d^{\II}:= GL_d/P_d^{\II}.\]
We equip $Y_d^{\I}(\mathbb{A})$ and $Y_d^{\II}(\mathbb{A})$ with topologies as homogeneous spaces of $GL_d(\mathbb{A})$.

We identify $\Omega_\eta^{d,\I}$ (resp. $\Omega_\eta^{d,\II}$) with $Y_d^{\I}(\mathbb{A})/k^\times$ (resp. $Y_d^{\II}(\mathbb{A})/k^\times$).

Define a map
\begin{align*}
\iota_d^{\II}:Y_d^{\II}(\mathbb{A}) &\to \mathbb{A}^d-\{0\} \\
g &\mapsto g\cdot(1,0,\dots,0)^t
\end{align*}
We see that $\iota_d^{\II}$ is well-defined, injective and continuous, but the topology on $Y_d^{\II}(\mathbb{A})$ is different from the subset topology.

Define an averaging map
\[\Av_d^{\II}: C_+(\mathbb{A}^d-\{0\})^{\mathbb{F}_q^\times}\to C^\infty(\Omega_{\eta}^{d,\II})\]
by the composition
\[C_+(\mathbb{A}^d-\{0\})^{\mathbb{F}_q^\times}=
C_+((\mathbb{A}^d-\{0\})/\mathbb{F}_q^\times) \xrightarrow{f_1}
C_+(Y_d^{\II}(\mathbb{A})/\mathbb{F}_q^\times) \xrightarrow{f_2}
C^\infty(Y_d^{\II}(\mathbb{A})/k^\times) = C^\infty(\Omega_{\eta}^{d,\II}) \]
where $f_1$ and $f_2$ are as follows. We identify $C^\infty(\mathbb{A}^d-\{0\})^{\mathbb{F}_q^\times}$ with $C^\infty((\mathbb{A}^d-\{0\})/\mathbb{F}_q^\times)$, and $C_+((\mathbb{A}^d-\{0\})/\mathbb{F}_q^\times)$ is the subgroup of $C^\infty((\mathbb{A}^d-\{0\})/\mathbb{F}_q^\times)$ corresponding to $C_+(\mathbb{A}^d-\{0\})^{\mathbb{F}_q^\times}$. The inclusion $\iota_d^{\II}$ gives a pullback $C^\infty((\mathbb{A}^d-\{0\})/\mathbb{F}_q^\times) \to C^\infty(Y_d^{\II}(\mathbb{A})/\mathbb{F}_q^\times)$. Let $C_+(Y_d^{\II}(\mathbb{A})/\mathbb{F}_q^\times)$ be the image of $C_+((\mathbb{A}^d-\{0\})/\mathbb{F}_q^\times)$. This gives $f_1$. The map $f_2$ is the pushforward, i.e., for $\varphi\in C_+(Y_d^{\II}(\mathbb{A})/\mathbb{F}_q^\times)$, $f_2(\varphi)$ is defined by
\[(f_2(\varphi))(x)=\sum_{a\in k^\times/\mathbb{F}_q^\times}\varphi(ax), \qquad (x\in Y_d^{\II}(\mathbb{A}))\]
From the definition of $C_+$ in \cref{(Section)formulation of the main result for the group of Tate toy horospherical divisors} we see that all but finitely many summands are zero.

For a rational differential $\omega$ on $X$ and an element $a\in\mathbb{A}$, we denote $\Res_\omega(a)$ to be the sum of residues of $a\omega$ at all closed points of $X$.

For a rational differential $\omega$ on $X$, we define a map
\begin{align*}
\iota_{d,\omega}^{\I}:Y_d^{\I}(\mathbb{A}) &\to (\mathbb{A}^d)^*-\{0\} \\
g &\mapsto (v\mapsto\Res_\omega((0,\dots,0,1)\cdot g^{-1}v)), \qquad (v\in \mathbb{A}^d)
\end{align*}
We see that $\iota_{d,\omega}^{\I}$ is well-defined, injective and continuous, but the topology on $Y_d^{\I}(\mathbb{A})$ is different from the subset topology.

Similarly to $\Av_d^{\II}$, we define an averaging map
\[\Av_d^{\I}: C_+((\mathbb{A}^d)^*-\{0\})^{\mathbb{F}_q^\times}\to C^\infty(\Omega_{\eta}^{d,\I})\]
by the composition
\[C_+((\mathbb{A}^d)^*-\{0\})^{\mathbb{F}_q^\times} =
C_+(((\mathbb{A}^d)^*-\{0\})/\mathbb{F}_q^\times) \xrightarrow{h_1}
C_+(Y_d^{\I}(\mathbb{A})/\mathbb{F}_q^\times) \xrightarrow{h_2}
C^\infty(Y_d^{\I}(\mathbb{A})/k^\times) = C^\infty(\Omega_{\eta}^{d,\I}). \]
Here $h_2$ is the pullback along $\iota_{d,\omega}^{\I}$ for a nonzero rational differential $\omega$ on $X$. We see that the composition does not depend on the choice of $\omega$.

\subsection{Formula for pullback of horospherical divisors under $\theta$}
Recall that the homomorphism
\[(\theta_\eta^{d,0})^*:
C_+((\mathbb{A}^d)^*-\{0\})^{\mathbb{F}_q^\times}\oplus C_+(\mathbb{A}^d-\{0\})^{\mathbb{F}_q^\times} \to
C^\infty(\Omega_\eta^{d,\I})\oplus C^\infty(\Omega_\eta^{d,\II})\]
is defined in \cref{(Section)theta as a direct sum}, and we have $(\theta_\eta^{d,0})^*=(\theta_\eta^{d,0,\I})^*\oplus(\theta_\eta^{d,0,\II})^*$
by \cref{theta as a direct sum}.

\begin{Proposition}\label{formula for pullback of toy horospherical divisors of type II under theta}
  $(\theta_{\eta}^{d,0,\II})^*=\Av_d^{\II}$. \qed
\end{Proposition}
The proof of \cref{formula for pullback of toy horospherical divisors of type II under theta} will be given in \cref{(Section)description of pullback of toy horospherical divisors of type II under theta}.

\begin{Lemma}\label{adeles with all residues being zero}
  An element $a\in\mathbb{A}$ is contained in $k$ if and only if $\Res_\omega(a)=0$ for all rational differentials $\omega$ of $X$. \qed
\end{Lemma}

In view of \cref{adeles with all residues being zero}, the following statement is dual to \cref{formula for pullback of toy horospherical divisors of type II under theta}.

\begin{Proposition}\label{formula for pullback of toy horospherical divisors of type I under theta}
  $(\theta_{\eta}^{d,0,\I})^*=\Av_d^{\I}$. \qed
\end{Proposition}

Now we obtain the main theorem of this section.
\begin{Theorem}\label{formula for pullback of toy horospherical divisors under theta}
  $(\theta_{\eta}^{d,0})^*=\Av_d^{\I}\oplus\Av_d^{\II}$.
\end{Theorem}
\begin{proof}
  The statement follows from \cref{formula for pullback of toy horospherical divisors of type II under theta,formula for pullback of toy horospherical divisors of type I under theta}.
\end{proof}

\subsection{A subspace of principal horospherical $\mathbb{Z}[\frac{1}{p}]$-divisors}
We normalize the Haar measure on $\mathbb{A}^d$ such that $O^d$ has measure 1.

Fix a nontrivial additive character $\psi:\mathbb{F}_q\to\mathbb{C}^\times$. We define the Fourier transform
\[\Four_\psi:C_c^\infty(\mathbb{A}^d;\mathbb{C})\to C_c^\infty((\mathbb{A}^d)^*;\mathbb{C})\]
such that for any $f\in C_c^\infty(\mathbb{A}^d;\mathbb{C}), \omega\in (\mathbb{A}^d)^*$, we have
\[\Four_\psi(f)(\omega)=\int_{\mathbb{A}^d} f(v)\psi(\omega(v))dv.\]

When $f\in C_c^\infty(\mathbb{A}^d;\mathbb{Z}[\frac{1}{p}])^{\mathbb{F}_q^\times}$, we have $\Four_\psi(f)\in C_c^\infty((\mathbb{A}^d)^*;\mathbb{Z}[\frac{1}{p}])^{\mathbb{F}_q^\times}$, and $\Four_\psi(f)$ does not depend on the choice of $\psi$.

Denote $C_0^\infty(\mathbb{A}^d)=\{f\in C_c^\infty(\mathbb{A}^d)|f(0)=\int_{\mathbb{A}^d} fdv=0\}$.

Combining \cref{space of principal rational Tate toy horospherical divisors,formula for pullback of toy horospherical divisors under theta}, we get the following statement.

\begin{Theorem}\label{a subspace of principal horospherical Q-divisors}
  For $d\ge 2$, the (partially) ordered abelian group of principal horospherical $\mathbb{Z}[\frac{1}{p}]$-divisors of $\Sht_{\eta,\all}^{d,\chi=0}$ contains the following subgroup
  \[\{(\Av_d^{\I}(\Four_\psi(f)),\Av_d^{\II}(f))\in C^\infty(\Omega_{\eta}^{d,\I};\mathbb{Z}[\tfrac{1}{p}])\oplus C^\infty(\Omega_{\eta}^{d,\II};\mathbb{Z}[\tfrac{1}{p}])|f\in C_0^\infty(\mathbb{A}^d;\mathbb{Z}[\tfrac{1}{p}])^{\mathbb{F}_q^\times}\}.\]
\end{Theorem}

\subsection{Multiplicity one for pullback of toy horospherical divisors}
Recall that the notation $O_D$ for a divisor $D$ of $X$ is defined in \cref{(Section)from a shtuka with all level structures to a Tate toy shtuka}. It is the c-lattice of $\mathbb{A}$ which consists of those adeles with poles bounded by $D$.

Let $\chi\in\mathbb{Z}$. Let $D', D''$ be two divisors of $X$ such that $(O_{D''}^d,O_{D'}^d)\in AP_{n_\chi}(\mathbb{A}^d)$. (See \cref{(Section)admissible pairs of c-lattices} for the definition of the notation $AP$.)

Let $E$ be a separably closed field over $\mathbb{F}_q$ and fix two morphisms $\alpha,\beta:\Spec E\to X$ satisfying condition $(*)$.

Let $\Sht_{E,D}^{d,\chi,D',D''}$ be the moduli stack which to each scheme $S$ over $\Spec E$ associates the groupoid of shtukas $\mathscr{F}$ over $S$ of rank $d$ with zero $\alpha\bcirc\mathsf{p}_S$ and pole $\beta\bcirc\mathsf{p}_S$ equipped with a structure of level $D$, satisfying the following conditions:

(i) $\chi(\mathscr{F})=\chi$.

(ii) For every $s\in S$ one has $H^0(X\times s,\mathscr{F}_s(D'\times s))=0$ and $H^1(X\times s,\mathscr{F}_s(D''\times s))=0$.

\cref{from a shtuka to a toy shtuka} gives a morphism $\Sht_{E,D}^{d,\chi,D',D''}\to\ToySht_{O_{D''}^d/O_{D'}^d}^{\chi+d\cdot\deg D''}\otimes E$. The image of the morphism lands in $\oToySht_{O_{D''}^d/O_{D'}^d}^{\chi+d\cdot\deg D''}$ since we have $\alpha\ne\beta$ by condition $(*)$.

Define $\Art_E$ to be the category of local Artinian rings whose residue fields are identified with $E$. Define $\Art_E^{(n)}$ to be the full subcategory of $\Art_E$ whose objects satisfy the condition that $x^{q^n}=0$ for all $x$ in the maximal ideal.

For a shtuka $\mathscr{G}$ over a perfect field, the notation $\mathscr{G}^{\I}$ and $\mathscr{G}^{\II}$ is defined in \cref{maximal trivial sub and maximal trivial quotient}.

\begin{Lemma}\label{deformation for pullback of toy horospherical divisors}
  Let $\mathscr{G}\in\Sht_{E,D}^{d,\chi,D',D''}(E)$. Suppose that $\mathscr{G}^{\II}=\mathscr{A}\otimes E$, where $\mathscr{A}$ is an invertible sheaf on $X$. Let $A\in\Art_E^{(1)}$. Let $\widetilde{\mathscr{G}}$ be a shtuka over $\Spec A$ extending the one over $\Spec E$. Let $\widetilde{L}={\pi_A}_*(\widetilde{\mathscr{G}}(D''\otimes A))\in\oToySht_{O_{D''}^d/O_{D'}^d}^{\chi+d\cdot\deg D''}(A)$ be the associated toy shtuka over $\Spec A$. Assume that $\widetilde{L}\supset J\otimes A$ for some $J\in\mathbf{P}_{O_{D''}^d/O_{D'}^d}$. Then $\widetilde{\mathscr{G}}$ contains $\mathscr{A}\otimes A$, and $\widetilde{\mathscr{G}}/(\mathscr{A}\otimes A)$ is locally free.
\end{Lemma}
\begin{Remark}
  Since $A\in\Art_E^{(1)}$, we have $\Phi_A^*\widetilde{\mathscr{G}}=(\Phi_E^*\mathscr{G})\otimes_E A$. We view $\widetilde{\mathscr{G}}$ as a subsheaf of $(\Phi_A^*\widetilde{\mathscr{G}})(\Gamma_{\widetilde{\alpha}})$, where $\widetilde{\alpha}:\Spec A\to X$ is the zero of $\mathscr{G}$. We also have $\Phi_E^*\mathscr{G}\supset\mathscr{A}\otimes E$. So $\mathscr{A}\otimes A$ in the above lemma is viewed as a subsheaf of $\Phi_A^*\widetilde{\mathscr{G}}$.
\end{Remark}
\begin{proof}[Proof of \cref{deformation for pullback of toy horospherical divisors}]
  We have $J\otimes E\subset H^0(X\otimes E, \mathscr{G}(D''\otimes E))$. Since $\Fr_E^*(J\otimes E)=J\otimes E$, \cref{maximal trivial sub and maximal trivial quotient} shows that $J\otimes E\subset H^0(X\otimes E, (\mathscr{G}(D''\otimes E))^{\II})$. Since $\mathscr{G}^{\II}=\mathscr{A}\otimes E$, \cref{twist of maximal trivial sub and maximal trivial quotient} implies that $(\mathscr{G}(D''\otimes E))^{\II}=\mathscr{A}(D'')\otimes E$. Together with the fact that $J$ is defined over $\mathbb{F}_q$, we deduce that $J\subset H^0(X,\mathscr{A}(D''))$. Let $\mathscr{J}_0$ be the subsheaf of $\mathscr{A}(D'')$ generated by $J$. Then $\mathscr{J}_0=\mathscr{A}(D''-D_0)$, where $D_0$ is an effective divisor of $X$. The assumption $\widetilde{L}\supset J\otimes A$ implies that $\widetilde{\mathscr{G}}(D''\otimes A)\supset\mathscr{J}_0\otimes A$. Hence
  \begin{equation}\label{deformation of shtuka contains trivial subbundle up to a finite subscheme of X}
  \widetilde{\mathscr{G}}\supset(\mathscr{A}\otimes A)(-D_0\otimes A).
  \end{equation}

  Since $A\in\Art_E^{(1)}$, we have $\Phi_A^*\widetilde{\mathscr{G}}=(\Phi_E^*\mathscr{G})\otimes_E A$. From the definition of shtuka we see that $\widetilde{\mathscr{G}}\supset((\Phi_E^*\mathscr{G})\otimes_E A)(-\Gamma_{\widetilde{\beta}})$, where $\widetilde{\beta}:\Spec A\to X$ is the pole of $\widetilde{\mathscr{G}}$. We also have $(\Phi_E^*\mathscr{G})\otimes_E A\supset \mathscr{A}\otimes A$. Thus
  \begin{equation}\label{deformation of shtuka contains trivial subbundle up to the pole}
  \widetilde{\mathscr{G}}\supset(\mathscr{A}\otimes A)(-\Gamma_{\widetilde{\beta}}).
  \end{equation}

  Condition $(*)$ implies that image of $\beta:\Spec E\to X$ is the generic point of $X$. Hence $D_0\otimes A$ and $\Gamma_{\widetilde{\beta}}$ are disjoint. Now (\ref{deformation of shtuka contains trivial subbundle up to a finite subscheme of X}) and (\ref{deformation of shtuka contains trivial subbundle up to the pole}) imply that $\widetilde{\mathscr{G}}\supset\mathscr{A}\otimes A$.

  The base change of the inclusion $\mathscr{A}\otimes A\hookrightarrow\widetilde{\mathscr{G}}$ to $X\otimes E$ is the inclusion $\mathscr{A}\otimes E\hookrightarrow\mathscr{G}$, which is injective at each point according to the definition of $\mathscr{G}^{\II}=\mathscr{A}\otimes E$ in \cref{maximal trivial sub and maximal trivial quotient}. Since $A\in \Art_E^{(1)}$, the morphism $X\otimes E\to X\otimes A$ is bijective on points. Thus the inclusion $\mathscr{A}\otimes A\hookrightarrow\widetilde{\mathscr{G}}$ has strictly constant rank 1. Now \cref{cokernel is free for constant rank morphism between locally free sheaves} implies that $\widetilde{\mathscr{G}}/(\mathscr{A}\otimes A)$ is locally free.
\end{proof}

Recall that $\mathcal{Z}_{\eta,1}^{d,\II}$ is defined in \cref{notation of horospherical divisors}.

Now we consider the case $\chi=0$. Denote $\mu:\Sht_{\eta,\all}^{d,\chi=0,D',D''}\to\oToySht_{O_{D''}^d/O_{D'}^d}^{d\cdot\deg D''}$ to be the morphism given by \cref{from a shtuka to a toy shtuka}.

\begin{Proposition}\label{multiplicity one for pullback of toy horospherical divisors}
  Let $J\in\mathbf{P}_{O_{D''}^d/O_{D'}^d}$. Let $\Delta=\Delta_{O_{D''}^d/O_{D'}^d,J}^{d\cdot\deg D''}$ be a toy horospherical divisor of $\oToySht_{O_{D''}^d/O_{D'}^d}^{d\cdot\deg D''}$. Then the multiplicity of $\mu^{-1}(\Delta)$ at $(g\cdot\mathcal{Z}_{\eta,1}^{d,\II})\cap\Sht_{\eta,\all}^{d,D',D''}$ is~$\le 1$ for all $g\in GL_d(\mathbb{A})_0$.
\end{Proposition}

\begin{proof}
  Let $\xi=\Spec L$ be the generic point of $(g\cdot\mathcal{Z}_{\eta,1}^{d,\II})\cap\Sht_{\eta,\all}^{d,D',D''}$. Let $R$ be the local ring of $\xi$ in $\Sht_{\eta,\all}^{d,D',D''}$. Let $S=\mu^{-1}(\Delta)\cap\Spec R$.

  Let $E$ be a separable closure of $L$. Let $\mathscr{F}$ be the shtuka over $\Spec E$. \cref{generic reducible shtuka of type II has maximal trivial sub of rank 1} shows that $(g^{-1}\cdot\mathscr{F})^{\II}\cong\mathscr{O}_X$. Hence  $\mathscr{F}^{\II}=\mathscr{A}\otimes E$ for some invertible sheaf $\mathscr{A}$ on $X$.

  Let $\mathcal{Y}$ be the closure of the image of the morphism $\RedSht_{\eta,\mathscr{A},\emptyset}^{d,1,\II}\to\Sht_{\eta,\emptyset}^d$. \cref{dimension of the stack of reducible shtukas} shows that $\mathcal{Y}$ is reduced. \cref{properties of transition maps between stacks of reducible shtukas} shows that the inverse image of $\mathcal{Y}$ under the morphism $\Spec R\to\Sht_{\eta,\emptyset}^d$ is equal to $\xi$.

  \cref{deformation for pullback of toy horospherical divisors} shows that for any $A\in\Art_E^{(1)}$, a morphism $\Spec A\to S$ factors through $\xi$ as long as the composition $\Spec E\to\Spec A\to S$ equals the composition $\Spec E\to\Spec L\to S$. Hence $S=\xi$.
\end{proof}

\subsection{Description of $(\theta_{\eta}^{d,0,\II})^*$}\label{(Section)description of pullback of toy horospherical divisors of type II under theta}
In this subsection, we describe $(\theta_{\eta}^{d,0,\II})^*$ in terms of whether a Tate toy horospherical divisor of $\oToySht_{\mathbb{A}^d}^{n_0}$ contains the image of a horospherical divisor of $\Sht_{\eta,\all}^{d,\chi=0}$ (See \cref{description of pullback of toy horospherical divisors under theta}). Then we finish the proof of \cref{formula for pullback of toy horospherical divisors of type II under theta} to obtain a formula for $(\theta_{\eta}^{d,0,\II})^*$.

Let $x\in \Omega_\eta^{d,\II}$. Denote $\mathcal{Z}_{\eta,x}^{d,\II}=x\cdot\mathcal{Z}_{\eta,1}^{d,\II}$. Let $\xi_{\eta,x}^{d,\II}$ denote the generic point of $\mathcal{Z}_{\eta,x}^{d,\II}$.

We denote
\[\mathbf{J}_x=\{J\in\mathbf{P}_{\mathbb{A}^d}| \theta_\eta^d(\xi_{\eta,x}^{d,\II})\in\Delta_{\mathbb{A}^d,J}^{n_0}\}.\]

\begin{Lemma}\label{description of pullback of toy horospherical divisors under theta}
  Let $f\in C_+(\mathbb{A}^d-\{0\})^{\mathbb{F}_q^\times}$. Then the set
  \[\mathbf{J}_{x,f}:=\{J\in\mathbf{J}_{x}|f(J)\ne 0\}\]
  is finite, and we have
  \[((\theta_\eta^{d,0,\II})^*f)(x)=\sum_{J\in \mathbf{J}_x}f(J).\]\qed
\end{Lemma}

\begin{proof}[Proof of \cref{formula for pullback of toy horospherical divisors of type II under theta}]
  The group $GL_d(\mathbb{A})_0$ acts transitively on $\Omega_{\eta}^{d,\II}$. The morphism $\theta_{\eta}^{d,0}$ is $GL_d(\mathbb{A})_0$-equivariant. Hence it suffices to show that $((\theta_{\eta}^{d,0,\II})^*f)(1)=(\Av_d^{\II}f)(1)$ for any $f\in C_+(\mathbb{A}^d-\{0\})^{\mathbb{F}_q^\times}$. This follows from \cref{description of pullback of toy horospherical divisors under theta,Tate toy horosphrical divisors containing the image of Z_1^II}.
\end{proof}

Recall that the notation $O_D$ for a divisor $D$ of $X$ is defined in \cref{(Section)from a shtuka with all level structures to a Tate toy shtuka}. Choose two divisors $D',D''$ of $X$ such that $(O_{D'}^d,O_{D''}^d)\in AP_{n_0}(\mathbb{A}^d)$ and $\xi_{\eta,x}^{d,\II}\in\Sht_{\eta,\all}^{d,\chi=0,D',D''}$.

For divisors $\widetilde{D'}\le D'\le D''\le \widetilde{D''}$ of $X$, we denote
\[U_{\widetilde{D'},\widetilde{D''}}^{n_0,D',D''}
=\oToySht_{O_{\widetilde{D''}}^d/O_{\widetilde{D'}}^d}^{ d\cdot\deg\widetilde{D''},O_{D'}^d/O_{\widetilde{D'}}^d,O_{D''}^d/O_{\widetilde{D'}}^d},\]
\[\mathbf{J}_{x,\widetilde{D'},\widetilde{D''}}^{D',D''}=
\{J\in \mathbf{J}_x|J\not\subset O_{D'}^d,J\subset O_{\widetilde{D''}}^d\}.\]
For $\overline{J}\in\mathbf{P}_{O_{\widetilde{D''}}^d/O_{\widetilde{D'}}}$, we denote $Z_{\widetilde{D'},\widetilde{D''},\overline{J}}^{D',D''}
=(\Delta_{O_{\widetilde{D''}}^d/O_{\widetilde{D'}}^d,\overline{J}}^{d\cdot\deg\widetilde{D''}}\cap U_{\widetilde{D'},\widetilde{D''}}^{n_0,D',D''})$. We denote
\[\overline{\mathbf{J}}_{x,\widetilde{D'},\widetilde{D''}}^{D',D''}=
\{\overline{J}\in\mathbf{P}_{O_{\widetilde{D''}}^d/O_{\widetilde{D'}}^d}|
(u_{\widetilde{D'},\widetilde{D''}}^{n_0,D',D''}\bcirc\theta_\eta^d)(\xi_{\eta,x}^{d,\II})\in Z_{\widetilde{D'},\widetilde{D''},\overline{J}}^{D',D''}\}\]
where $u_{\widetilde{D'},\widetilde{D''}}^{n_0,D',D''}:\oToySht_{\mathbb{A}^d}^{n_0,D',D''}\to U_{\widetilde{D'},\widetilde{D''}}^{n_0,D',D''}$ is the projection.

For $J\in\mathbf{J}_{x,\widetilde{D'},\widetilde{D''}}^{D',D''}$, we have $u_{\widetilde{D'},\widetilde{D''}}^{n_0,D',D''}(\Delta_{\mathbb{A}^d,J}^{n_0,D',D''})\subset Z_{\widetilde{D'},\widetilde{D''},\overline{J}}^{D',D''}$, where $\overline{J}=\im(J\cap O_{\widetilde{D''}}^d\to O_{\widetilde{D''}}^d/O_{\widetilde{D'}}^d)$. Thus we get a map
\begin{align*}
\mathbf{J}_{x,\widetilde{D'},\widetilde{D''}}^{D',D''} &\to \overline{\mathbf{J}}_{x,\widetilde{D'},\widetilde{D''}}^{D',D''}\\
J &\mapsto \overline{J}
\end{align*}

\begin{Lemma}\label{bijection between toy horospherical divisors containing the image of a horospherical divisor}
  The map $\mathbf{J}_{x,\widetilde{D'},\widetilde{D''}}^{D',D''}\to \overline{\mathbf{J}}_{x,\widetilde{D'},\widetilde{D''}}^{D',D''}$ is bijective.
\end{Lemma}
\begin{proof}

  Let $\overline{J}\in \overline{\mathbf{J}}_{x,\widetilde{D'},\widetilde{D''}}^{D',D''}$. By \cref{Tate toy horospherical divisor as a projective limit}, it suffices to show that for any two divisors $\mathring{D'},\mathring{D''}$ of $X$ such that $\mathring{D'}\le\widetilde{D'}\le\widetilde{D''}\le\mathring{D''}$, there is a unique $\mathring{J}\in\mathring{\mathfrak{J}}$ such that $(u_{\mathring{D'},\mathring{D''}}^{n_0,D',D''}\bcirc\theta_\eta^d)(\xi_{\eta,x}^{d,\II})\in Z_{\mathring{D'},\mathring{D''},\mathring{J}}^{D',D''}$, where the set $\mathring{\mathfrak{J}}$ is defined by
  \[\mathring{\mathfrak{J}}=\{J\in\mathbf{P}_{O_{\mathring{D''}}^d/O_{\mathring{D'}}^d}|
    \im(J\cap O_{\widetilde{D''}}^d\to O_{\widetilde{D''}}^d/O_{\widetilde{D'}}^d)=\overline{J}\}.\]

  Let $\overline{\Delta}=Z_{\widetilde{D'},\widetilde{D''},\overline{J}}^{D',D''}$. Denote $u: U_{\mathring{D'},\mathring{D''}}^{n_0,D',D''} \to U_{\widetilde{D'},\widetilde{D''}}^{n_0,D',D''}$ to be the projection. Consider the composition of morphisms
  \[\Sht_{\eta,\all}^{d,\chi=0,D',D''}\xrightarrow{\theta_{\eta}^d}\oToySht_{\mathbb{A}^d}^{n_0,D',D''} \xrightarrow{u_{\mathring{D'},\mathring{D''}}^{n_0,D',D''}}
  U_{\mathring{D'},\mathring{D''}}^{n_0,D',D''}\xrightarrow{u} U_{\widetilde{D'},\widetilde{D''}}^{n_0,D',D''}.\]

  \cref{pullback of toy horospherical divisors under transition maps} implies that
  \begin{equation}\label{reference 1 of pullback of toy horospherical divisors under transition maps}
  u^*(\overline{\Delta})=\sum_{J\in \mathring{\mathfrak{J}}}Z_{\mathring{D'},\mathring{D''},J}^{D',D''}.
  \end{equation}
  \cref{transition maps of toy shtukas are smooth} shows that the morphism $u$ is smooth. Hence
  \[u^{-1}(\overline{\Delta})=\bigcup_{J\in \mathring{\mathfrak{J}}}Z_{\mathring{D'},\mathring{D''},J}^{D',D''}.\]
  The assumption $\overline{J}\in \overline{\mathbf{J}}_{x,\widetilde{D'},\widetilde{D''}}^{D',D''}$ implies that $(u\bcirc u_{\mathring{D'},\mathring{D''}}^{n_0,D',D''}\bcirc\theta_{\eta}^d)(\xi_{\eta,x}^{d,\II})\in \overline{\Delta}$. Hence $(u_{\mathring{D'},\mathring{D''}}^{n_0,D',D''}\bcirc\theta_{\eta}^d)(\xi_{\eta,x}^{d,\II})\in Z_{\mathring{D'},\mathring{D''},\mathring{J}}^{D',D''}$ for some $\mathring{J}\in\mathring{\mathfrak{J}}$. This shows the existence of $\mathring{J}$.

  \cref{multiplicity one for pullback of toy horospherical divisors} shows that the pullback of $\overline{\Delta}$ to $\Sht_{\eta,\all}^{d,\chi=0,D',D''}$ has multiplicity~$\le 1$ at $\mathcal{Z}_{\eta,x}^{d,\II}$. From formula (\ref{reference 1 of pullback of toy horospherical divisors under transition maps}), we see that there is at most one $\mathring{J}\in\mathring{\mathfrak{J}}$ such that $(u_{\mathring{D'},\mathring{D''}}^{n_0,D',D''}\bcirc\theta_\eta^d)(\xi_{\eta,x}^{d,\II})\in Z_{\mathring{D'},\mathring{D''},\mathring{J}}^{D',D''}$. We obtain the uniqueness of $\mathring{J}$.
\end{proof}

\begin{proof}[Proof of \cref{description of pullback of toy horospherical divisors under theta}]
  Let $Q\in\mathfrak{O}_{\mathbb{A}^d}^{n_0}$ be the divisor of $\oToySht_{\mathbb{A}^d}^{n_0}$ corresponding to $0\oplus f$. Let $Q^{D',D''}\in\mathfrak{O}_{\mathbb{A}^d}^{n_0,D',D''}$ be the restriction of $Z$ to $\oToySht_{\mathbb{A}^d}^{n_0,D',D''}$. From \cref{open subscheme of Tate toy shtukas as a projective limit} we see that there exist two divisors $\widetilde{D'},\widetilde{D''}$ of $X$ such that $\widetilde{D'}\le D'\le D''\le\widetilde{D''}$ and that $Q^{D',D''}$ equals the pullback of a toy horospherical divisor $\overline{Z}$ of $U_{\widetilde{D'},\widetilde{D''}}^{n_0,D',D''}$. Thus we have an inclusion $\mathbf{J}_{x,f}\subset\mathbf{J}_{x,\widetilde{D'},\widetilde{D''}}^{D',D''}$.
  Since $\overline{\mathbf{J}}_{x,\widetilde{D'},\widetilde{D''}}^{D',D''}$ is finite, $\mathbf{J}_{x,\widetilde{D'},\widetilde{D''}}^{D',D''}$ is also finite by \cref{bijection between toy horospherical divisors containing the image of a horospherical divisor}, hence is $\mathbf{J}_{x,f}$.

  \cref{multiplicity one for pullback of toy horospherical divisors} implies that
  \[((\theta_\eta^{d,0,\II})^*f)(x)=\sum_{\overline{J}\in\overline{\mathbf{J}}_{x,\widetilde{D'},\widetilde{D''}}^{D',D''}}
  m(f,\overline{J})\]
  where $m(f,\overline{J})$ is the multiplicity of $\overline{Z}$ at $Z_{\widetilde{D'},\widetilde{D''},\overline{J}}^{D',D''}$. For any  $\overline{J}\in\overline{\mathbf{J}}_{x,\widetilde{D'},\widetilde{D''}}^{D',D''}$, \cref{bijection between toy horospherical divisors containing the image of a horospherical divisor} shows that there is a unique $J\in \mathbf{J}_{x,\widetilde{D'},\widetilde{D''}}^{D',D''}$ such that $\overline{J}=\im(J\cap O_{\widetilde{D''}}^d\to O_{\widetilde{D''}}^d/O_{\widetilde{D'}}^d)$. \cref{equality of multiplicities at toy horospherical divisors for a projective limit} shows that $m(f,\overline{J})=f(J)$. The statement follows.
\end{proof}

\appendix
\section{Recollections of linear algebra}
\subsection{}
\begin{Definition}
For two coherent locally free sheaves $\mathscr{F},\mathscr{G}$ over a scheme $S$, we say that a morphism $\psi:\mathscr{F}\to\mathscr{G}$ has rank at most $r$ if the induced morphism $\bigwedge^{r+1}\psi:\bigwedge^{r+1}\mathscr{F}\to\bigwedge^{r+1}\mathscr{G}$ is zero. In this case, we say that $\psi$ has strictly constant rank $r$ if the induced morphism $\bigwedge^{r}\psi:\bigwedge^{r}\mathscr{F}\to\bigwedge^{r}\mathscr{G}$ is nonzero at any point $s\in S$.
\end{Definition}

\begin{Lemma}\label{cokernel is free for constant rank morphism between modules}
Let $\psi:M\to N$ be a morphism between finitely generated free modules over a local ring $A$. Put $e=\rank N$. Then $\psi$ has strictly constant rank $r$ if and only if $\coker\psi$ is free of rank $e-r$.
\end{Lemma}
\begin{proof}
If $\psi$ has strictly constant rank $r$, then the matrix for $\psi$ under two bases of $M,N$ has an $(r\times r)$-minor with determinant being a unit of $A$. Hence we can choose bases of $M,N$ so that $\psi$ has matrix
\[\begin{pmatrix}\Id_{r\times r} & 0\\0 & B\end{pmatrix}.\]
Since $\psi$ has strictly constant rank $r$, we get $B=0$. Hence $\coker\psi$ is free of rank $e-r$.

If $\coker\psi$ is free of rank $(e-r)$, then $\im\psi$ is free of rank $r$, so $\psi$ has rank at most $r$. Let $E$ be the residue field of $A$. Then $\coker(\psi\otimes E)=(\coker\psi)\otimes E$ has dimension $(e-r)$ over $E$. Hence $\psi\otimes E$ has rank $r$. Therefore $\psi$ has strictly constant rank $r$.
\end{proof}

\begin{Lemma}\label{cokernel is free for constant rank morphism between locally free sheaves}
Let $S$ be a scheme. Let $\psi:\mathscr{F}\to\mathscr{G}$ be a morphism of coherent locally free sheaves on $S$. Put $g=\rank\mathscr{G}$. Then $\psi$ has strictly constant rank $r$ if and only if $\coker\psi$ is locally free of rank $g-r$.
\end{Lemma}
\begin{proof}
Since $\coker\psi$ is a finitely presented $\mathscr{O}_S$-module, it suffices to prove the statements for each stalks. This follows from \cref{cokernel is free for constant rank morphism between modules}.
\end{proof}

\begin{Corollary}\label{splittings of source and target for morphism of locally constant rank}
Let $S,\psi,\mathscr{F},\mathscr{G}$ be as in \cref{cokernel is free for constant rank morphism between locally free sheaves}, then the following two exact sequences split locally.
\[\begin{tikzcd}0\arrow[r]&\ker\psi\arrow[r]&\mathscr{F}\arrow[r]&\im\psi\arrow[r]&0,\end{tikzcd}\]
\[\begin{tikzcd}0\arrow[r]&\im\psi\arrow[r]&\mathscr{G}\arrow[r]&\coker\psi\arrow[r]&0.\end{tikzcd}\]
In particular, $\ker\psi$ and $\im\psi$ are locally free.
\end{Corollary}

\begin{Corollary}\label{equivalent conditions of the relative position of two subbundles}
Let $\mathscr{L}_1,\mathscr{L}_2$ be subbundles of a coherent locally free sheaf $\mathscr{F}$ over a scheme $S$. Assume $\mathscr{L}_1,\mathscr{L}_2$ have rank $n$. Then the following conditions are equivalent:

(i) the morphism $\mathscr{L}_1\to\mathscr{F}/\mathscr{L}_2$ has strictly constant rank $r$;

(ii) the morphism $\mathscr{L}_2\to\mathscr{F}/\mathscr{L}_1$ has strictly constant rank $r$;

(iii) the morphism $\mathscr{L}_1\oplus\mathscr{L}_2\to\mathscr{F}$ has strictly constant rank $n+r$.

When the above conditions hold, $\mathscr{L}_1+\mathscr{L}_2$ and $\mathscr{L}_1\cap\mathscr{L}_2$ are subbundles of $\mathscr{F}$, and the four sheaves $(\mathscr{L}_1+\mathscr{L}_2)/\mathscr{L}_1,(\mathscr{L}_1+\mathscr{L}_2)/\mathscr{L}_2,\mathscr{L}_1/(\mathscr{L}_1\cap\mathscr{L}_2),
\mathscr{L}_2/(\mathscr{L}_1\cap\mathscr{L}_2)$ are locally free of rank $r$.
\end{Corollary}

\bibliography{../../../Bibliography/Reference}
\bibliographystyle{amsplain}
\end{document}